\documentclass[12pt,twoside]{article}
\usepackage[hmargin=0.8in,vmargin=0.9in]{geometry}
\usepackage{fancyhdr}
\usepackage{graphicx}
\usepackage{amssymb}
\usepackage{amsmath}
\usepackage{amsthm}
\usepackage{amsfonts}
\usepackage{mathrsfs}
\usepackage{mathtools}
\usepackage{bm}
\usepackage{color}
\usepackage{setspace}
\usepackage{exscale}
\usepackage{relsize}
\usepackage{float}
\usepackage{cite}
\usepackage{hyperref}
\usepackage{cleveref}
\usepackage{nicefrac}
\usepackage{placeins}

\usepackage{etoolbox} 
\AtBeginEnvironment{pmatrix}{\everymath{\displaystyle}}
\AtBeginEnvironment{bmatrix}{\everymath{\displaystyle}}

\usepackage{booktabs,multirow} 
\usepackage{array} 
\usepackage{paralist} 

\usepackage{sectsty}

\sectionfont{\fontsize{14}{15}\selectfont}
\subsectionfont{\fontsize{13}{15}\selectfont}

\theoremstyle{plain}                    
\newtheorem{thm}{Theorem}[section]

\newtheorem{rmk}[thm]{Remark}

\usepackage{tikz}
\usetikzlibrary{positioning}
\usepackage{xcolor}

\numberwithin{equation}{section}
\numberwithin{figure}{section}
\numberwithin{table}{section}

\newcommand*\xbar[1]{%
  \hbox{%
    \vbox{%
      \hrule height 0.5pt 
      \kern0.4ex
      \hbox{%
        \kern-0.05em
        \ensuremath{#1}%
        \kern-0.00em
      }%
    }%
  }%
}

\usepackage[utf8]{inputenc}
\usepackage[english]{babel}
\usepackage{empheq}
\usepackage{subcaption} 
\usepackage{epsfig}
\usepackage{bbm}
\usepackage{sectsty}

\newcommand{\bb}[1]{{\bf{#1}}} 
\newcommand{\mc}[1]{{\mathcal{#1}}} 

\newcommand{\p}[1]{{\left( #1 \right)}}
\newcommand{\set}[1]{{\left\{ #1 \right\}}}
\newcommand{\br}[1]{{\left[ #1 \right]}}
\newcommand{\abs}[1]{{\left| #1 \right|}}

\def\d{\partial}

\def\eps{\varepsilon}

\newcommand{\pd}[3][]{\frac{\partial^{#1} #2}{\partial #3^{#1}}}

\newcommand{\kph}{{k+\frac{1}{2}}}
\newcommand{\kmh}{{k-\frac{1}{2}}}
\newcommand{\jph}{{j+\frac{1}{2}}}
\newcommand{\jmh}{{j-\frac{1}{2}}}

\newcommand{\dy}{\Delta y}
\newcommand{\dz}{\Delta \zeta}

\allowdisplaybreaks[1]

\graphicspath{{Figures/}}

\title{Moment Approximations to Magnetic Rotating Shallow Flows}
\author{ %
\hyperlink{https://orcid.org/0000-0002-8822-461X}{Julian Koellermeier}%
\thanks{Bernoulli Institute for Mathematics, Computer Science and Artificial Intelligence, University of Groningen, 9747 AG Groningen,
The Netherlands;}$\ ^,$%
\thanks{Department of Mathematics, Computer Science and Statistics, Ghent University,
9000 Ghent, Belgium; {\tt julian.koellermeier@ugent.be}}, 
\hyperlink{https://orcid.org/0000-0002-4099-2978}{Michael Redle}%
\thanks{Department of Mathematics, RWTH Aachen University, Schinkelstraße 2, 52062 Aachen, Germany;}$\ ^,$%
\thanks{Institute of Mathematics, TU Clausthal, Erzstraße 1, 38678 Clausthal-Zellerfeld, Germany; {\tt michael.redle@tu-clausthal.de}}, and %
\hyperlink{https://orcid.org/0000-0003-0008-2061}{Manuel Torrilhon}%
\thanks{Department of Mathematics, RWTH Aachen University, Schinkelstraße 2, 52062 Aachen, Germany; {\tt mt@acom.rwth-aachen.de}}
}

\begin{document}
\date{\today}
\maketitle

\begin{abstract}
Originally introduced to describe a transition region in stars, the magnetic rotating shallow water (MRSW) model is now used in many solar physics and geophysical applications. 
Derived from the 3-D incompressible magnetohydrodynamic system, the shallow nature of these applications motivates depth-averaging of both the velocities and magnetic fields. 
This is advantageous in terms of computational efficiency---but at the loss of vertical information, thus limiting the predictive power of the MRSW model. 
To overcome this problem, we employ higher-order vertical moments, but now in the context of conductive fluids. 
In doing so, the new approximation maintains non-constant vertical profiles of both the horizontal magnetic fields and horizontal velocities, while still remaining in the simplified 2-D framework corresponding to depth integration. 
In this work, we extend the derivation of the shallow water moment equations to derive the MRSW moment system of arbitrary order; i.e., we represent the vertical profiles of the velocities---and now additionally the magnetic fields---by arbitrary-order polynomial expansions, and close the new expanded 2-D system with evolution equations for these polynomial coefficients, found via Galerkin projection. 
Through numerical experiments for MRSW moment systems up to third-order, we demonstrate that these moment approximations reduce model error without significantly sacrificing computational efficiency.

\end{abstract}



\newpage

\section{Introduction}\label{sec1}

The rotating shallow water magnetohydrodynamic (MHD) equations, also referred as magnetic rotating shallow water (MRSW) equations, describe the behavior of a thin layer, conductive
fluid with a free surface.
Since their introduction in  \cite{Gilman_2000} in the context of solar tachocline (the thin
layer of the Sun between the convective and radiative zones), the number of applications
has vastly grown; some examples include neutron star atmosphere dynamics \cite{Heng_2009,Zaqarashvili2011Rossby} and accreting
matter flows \cite{inogamov2010spread}, and exoplanets with magnetically active atmospheres \cite{cho2008atmospheric,heng2015atmospheric,Heng_2014,Showman_2009}. 
While these applications
could alternatively be described by the full, spatially 3-D incompressible MHD system, the flow in the vertical direction is much smaller than that in the horizontal directions---motivating the introduction of the ‘shallow fluid’ assumption.

Through a suitable vertical averaging, the shallow fluid assumption ultimately leads to a reduced complexity, two-dimensional depth-averaged model (e.g., the MRSW system).
Consequently, computational costs of the MRSW model are significantly lower than its fully 3-D incompressible MHD counterpart---where the difficulties lie not only in the extra vertical dimension, but also in the computational complexities that arise with a free surface.
Unfortunately, however, depth-averaging the velocity and the magnetic field comes at the cost of losing vertical information. 
In some cases, this may not be critical; for example, in cases without a magnetic field, so called `plug-flow' arises often in geophysical contexts; see, e.g., \cite{SLEATH1999Plug}. 
Many other applications, though, have significant vertical profiles of velocity and magnetic fields; see, e.g., \cite{Dikpati_1999}, where results in the context of a solar tachocline are presented.

In contexts without magnetic fields, such as the standard shallow water equations, the accuracy loss of depth-averaged models when dealing with varying velocity profiles is already well-known \cite{schaefer2013velocity,christen2010ramms}.
This has led to a variety of more advanced models for free-surface flows.
Two typical directions can be distinguished:
(i) the consideration of multiple layers with piecewise constant fluid velocities---the so-called multi-layer models \cite{Fernandez-Nieto2016}, which unfortunately lose important mathematical properties for the interesting case of more than two layers; and
(ii) using the Green-Naghdi approach to model vertical variation with additional shape parameters; see, e.g.,\cite{GreenNaghdi1976,Parisot2019}; leading to additional equations that have to be derived on a case-by-case basis.
Note that both of these methodologies to improve vertical profile approximations have also been extended to the MRSW system; see, e.g., \cite{alonso2021asymptotic,hunter2015waves,Zeitlin2013} for examples of multilayer model extensions and, e.g., \cite{Dellar2003} for the extension of the Green-Naghdi approach to MRSW.

A promising new method has been recently employed in \cite{Kowalski2019Moment} in which the authors (i) represent the vertical profile of the velocity by an $M$th-degree polynomial; and (ii) derive additional evolution equations for the $M$ polynomial coefficients using a Galerkin projection of (a rescaled version of) the incompressible Navier-Stokes equations. 
Consequently, one obtains the hierarchical model known as the Shallow Water Moment Equations (SWME).
Since their derivation, several studies have subsequently analyzed hyperbolicity \cite{KoellermeierRominger2020}, steady state solutions \cite{Koellermeier2020steady}, and the behavior near equilibrium manifolds \cite{Koellermeier2020Eq} of the SWME, further revealing several beneficial analytical properties of the models.
In addition, versions of the SWME have been applied to sediment transport \cite{garres2021}, and extended to 2-D \cite{bauerle2024rotationalinvariancehyperbolicityshallow,VerbiestKoellermeier}, including open channel flows \cite{Steldermann10112023} and non-hydrostatic cases \cite{scholz2023}.

In this work, we extend this SMWE approach, i.e., employing higher-order vertical moments, to magnetic rotating flows. 
This new approximation of the 3-D incompressible MHD system maintains non-constant vertical profiles of the magnetic field in addition to those of the horizontal velocities, while still remaining in the computationally efficient 2-D framework corresponding to depth integration. 
We present the derivation of this extended approach, resulting in the MRSW moment system of arbitrary order.
The obtained system naturally gives rise to the so-called Godunov-Powell source---a nonconservative term often used to numerically treat the ever-important divergence-free constraint of the magnetic field.
Utilizing this term, we discretize both the moment and the (rescaled) incompressible Navier-Stokes systems using a direct extension of that introduced in \cite{Chertock2024Locally}; namely, we use the Godunov-Powell source to help treat the divergence of the magnetic field, and discretize the nonconservative system using the second-order path-conservative central-upwind (PCCU) scheme---an unstaggered Riemann-problem-solver-free finite volume scheme. 
Through 1-D numerical results of MRSW moment systems of up to third-order (i.e., third-degree polynomial representations of the vertical profiles), we show that the proposed model hierarchy better approximates the derived 2-D reference system than the constant-profile MRSW model. 
After demonstrating that these models can increase the predictive power of conductive shallow flow models without significantly sacrificing computational efficiency,
we present an exploratory example on magneto-geostrophic adjustment.
We show that an initial velocity vertical profile induces a vertical profile of the magnetic field and alters wave behavior (particularly in the Alfv\'en waves), further demonstrating the importance of additionally representing magnetic field vertical profiles.

The rest of this paper is organized as follows. 
In \S\ref{sec3}, we derive the vertically resolved reference MRSW model, which is used both for error analysis and further derivation. 
We follow with the derivation of the arbitrary-order MRSW moment equations in \S\ref{sec4}.
Section \ref{sec5} then describes the numerical methods needed to numerically represent both the moment system (\S\ref{sec5.1}) and the reference system (\S\ref{sec5.2}), as each requires its own discretization.
We follow by presenting numerical tests of the MRSW moment system in \S\ref{sec6}. Finally, we make our concluding remarks in \S\ref{sec7}.


\section{Vertically resolved magnetic rotating shallow water reference model}\label{sec3}
To derive the MRSW moment approximation, we first map the vertical variable $z \in [Z,h+Z]$ in the ideal incompressible magnetohydrodynamic (MHD) equations ($Z$ and $h$ being the bathymetry and fluid thickness, respectively) to a normalized vertical variable $\zeta \in [0,1]$.
The resulting mapped equations will be referred to as the `vertically resolved MRSW reference system,' and forming these equations benefits twofold.
First, the 3-D mapping obtained in this section is directly used to derive the reduced-dimension 2-D MRSW moment model, which uses polynomial expansions to form higher-order vertical profiles and additional evolution equations for the polynomial coefficients; see \S\ref{sec4} for further details.  
Second, the vertically resolved MRSW reference system provides a crucial validation aspect of the model hierarchy we form in \S\ref{sec4}.
This is because the reference model still fully resolves the $z$-profiles of the velocity and magnetic field, in addition to matching the shallow flow assumption used to derive the MRSW system; i.e., it additionally assumes magneto-hydrostatic equilibrium.
Furthermore, this MRSW reference system renders a full 3-D representation of the problem of interest---thus providing a precise model to evaluate the error of the MRSW moment approximation, which reduces the dimension through vertical averaging. 

%

The original equations used to derive the vertically resolved MRSW reference model are the incompressible ideal MHD equations with rotation.
In addition, we assume an inviscid fluid, zero magnetic diffusivity, and a rotation term described by the Coriolis force. 
Doing so, the incompressible ideal MHD system reads
\begin{equation}
    \begin{aligned}
    &\bm u_t + \nabla \cdot \br{\bm u \bm u^\top + \p{\frac{1}{\rho}p + \frac{1}{2}\abs{\bm b}^2}I 
    - \bm b \bm b^\top} = -f \bm{e}_z\times \bm{u}-g \bb{e}_z,	\\[4pt]
    & \bm b_t - \nabla \times \p{\bm u \times \bm b} = 0,	\\[4pt]
    &\nabla \cdot \bm u = 0, \\[4pt]
    &\nabla \cdot \bm b = 0,
    \end{aligned}
    \label{eq:inc_mhd}
\end{equation}
%
where $(x,y,z)$ denotes the spatial coordinate, $t$ is time, $\bb u = \p{u,v,w}^\top$ denotes velocity, $p$ is the thermal pressure, $\bm b = \p{a,b,c}^\top$ denotes the magnetic field scaled by $1/\sqrt{\mu_0 \rho}$ (thus having units of velocity), $f = f(y)$ is the Coriolis parameter, $g$ is the constant acceleration due to gravity, $\rho$ denotes the (assumed constant) density, $I$ denotes the $3\times3$ identity matrix, and $\bb e_z = \p{0,0,1}^\top$ is the vertical unit vector. 

Since we are interested in magnetic free-surface shallow flow, however, we may start with a simpler version of system \eqref{eq:inc_mhd}---Appendix \ref{secA1} presents this standard dimensionless scaling, resulting in the equations 
\eqref{eq:start}--\eqref{eq:mhp_sc}. 
For sake of readability, we will use the analogous dimensional formulation of \eqref{eq:start}--\eqref{eq:mhp_sc} as the starting point in deriving the MRSW moment framework; that is, we start with the system
\begin{align}
	&u_t + \br{u^2 +  \p{\frac{1}{\rho} p + \frac{1}{2} \abs{\bm b}^2} - a^2}_x 
	+ \br{uv - ab}_y + \br{uw - ac}_z = fv, \label{eq:mo_x}\\[4pt]
	&v_t + \br{uv - ab}_x + \br{v^2 +  \p{\frac{1}{\rho}p + \frac{1}{2} \abs{\bm b}^2} - b^2}_y 
	+ \br{vw - bc}_z = -fu, \label{eq:mo_y}\\[4pt]
	& a_t + \p{av-bu}_y + \p{aw-cu}_z = 0, \label{eq:mag_x}\\[4pt]
	& b_t + \p{bu-av}_x + \p{bw-cv}_z = 0, \label{eq:mag_y}\\[4pt]
	&u_x + v_y + w_z = 0,\label{eq:divu}\\[4pt]
	&a_x + b_y + c_z = 0, \label{eq:divb}
\end{align}
%
paired with the magneto-hydrostatic pressure relation
\begin{equation}\label{eq:MHpressure}
p + \frac{1}{2} \rho \abs{\bm b}^2 = \p{h + Z - z}\rho g,
\end{equation}
and the free surface boundary conditions; which, for the vertical velocity, read
\begin{equation}
\begin{aligned}
	(h+Z)_t ={\bm u} \cdot \bm n \qquad \implies \qquad
	&(h+Z)_t + (u,v)\Big|_{z = h+Z} \cdot \nabla (h + Z) = w \Big|_{z = h+Z}\ ,\\[4pt]
	Z_t ={\bm u} \cdot \bm n \qquad \implies \qquad
	&Z_t + (u,v)\Big|_{z = Z} \cdot \nabla Z =  w \Big|_{z = Z}\ ,
\end{aligned}
\label{eq:ubc}
\end{equation}
and for the vertical magnetic field, we use
\begin{equation}
\begin{aligned}
	{\bm b} \cdot \bm n = 0 \qquad \implies \qquad 
	& c \Big|_{z = h+Z} = (a,b)\Big|_{z = h+Z} \cdot \nabla (h + Z), \\[4pt]
	\qquad \implies \qquad & c \Big|_{z = Z} = (a,b)\Big|_{z = Z} \cdot \nabla Z.
\end{aligned}
\label{eq:bbc}
\end{equation}
Note again that here, $h(t,x,y)$ denotes the fluid thickness, and $Z(t,x,y)$ denotes the bottom topography. 

\subsection{Mapping}\label{sec3.1}

To render the reference system \eqref{eq:mo_x}--\eqref{eq:bbc} more accessible in the derivation of the MRSW moment approximation, we follow \cite{Kowalski2019Moment} and first map the vertical variable $z \in [Z,h+Z]$ to the normalized vertical variable $\zeta \in [0,1]$ using 
\begin{equation*}\label{eq:zeta}
\zeta = \frac{z - Z(t,x,y)}{h(t,x,y)},
\end{equation*}
or equivalently, $z = \zeta h + Z$.
Let us consider the arbitrary function $\psi(t,x,y,z)$ dependent on space and time. 
Then its mapped counterpart reads
\begin{equation*}\label{eq:arbmap}
	\tilde{\psi}(t,x,y,\zeta) = \psi(t,x,y,\zeta h + Z ).
\end{equation*}
To transform the system \eqref{eq:mo_x}--\eqref{eq:bbc} to also be dependent on $\zeta$ instead of $z$, we additionally require the differential operators of the mapping $\psi$, which read
\begin{equation}
\begin{aligned}
	&h\psi_s = \p{h \tilde \psi}_s - \br{\p{\zeta h + Z}_s \tilde \psi}_\zeta, \qquad s \in (t, x, y),\\[4pt]
	&h\psi_z = \tilde \psi_\zeta.
\end{aligned}
\label{eq:mappings}
\end{equation}
The remainder of \S\ref{sec3.1} describes the implementation of this mapping on the starting system \eqref{eq:mo_x}--\eqref{eq:bbc}.

\subsubsection{Mass balance mapping}\label{sec3.1.1}

Starting from the velocity divergence in \eqref{eq:divu}, we multiply by fluid thickness $h$, apply the transformation rules in \eqref{eq:mappings}, and rearrange to recover the mapped mass balance:
\begin{equation}\label{eq:w_tilde}
\tilde w_\zeta = \br{\p{\zeta h + Z}_x \tilde u + \p{\zeta h + Z}_y \tilde v}_\zeta - (h \tilde u )_x - (h \tilde v )_y.
\end{equation}
Integrating \eqref{eq:w_tilde} over $\zeta$, we obtain 
$$
\tilde w \Big|_{\zeta = 1} -  \tilde w \Big|_{\zeta = 0} =  \tilde u \Big|_{\zeta = 1}(h+Z)_x - \tilde u \Big|_{\zeta = 0}Z_x + \tilde v \Big|_{\zeta = 1}(h + Z)_y - \tilde v \Big|_{\zeta = 0} Z_y - \p{h \int_0^1 \tilde u \ d\zeta}_x - \p{h \int_0^1 \tilde v \ d\zeta}_y.
$$
Applying the free surface boundary conditions for velocity in \eqref{eq:ubc}, we obtain the standard depth-averaged fluid height equation of both the shallow water and MRSW systems:
\begin{equation}\label{eq:sw_mass}
	h_t + \p{hu_m}_x + \p{hv_m}_y = 0,
\end{equation}
where $\displaystyle u_m = \int_0^1 \tilde u \, d\zeta$ and $\displaystyle v_m = \int_0^1 \tilde v \,d\zeta$ are the depth-averaged mean velocities in the $x$- and $y$-directions, respectively.

\subsubsection{Mapping of the divergence-free condition}\label{sec3.1.2}

Recovering the mapped version of the magnetic field divergence-free condition follows similarly to that of the mass balance.
Starting with the magnetic divergence-free constraint in \eqref{eq:divb}, multiplying by fluid thickness $h$, and using the differential operator transformations in \eqref{eq:mappings}, we find
\begin{equation}\label{eq:c_tilde}
	\tilde c_\zeta = \br{\p{\zeta h + Z}_x \tilde a + \p{\zeta h + Z}_y \tilde b}_\zeta 
	- (h \tilde a )_x - (h \tilde b )_y.
\end{equation}
This is again integrated over $\zeta$ and evaluated, yielding
$$
\tilde c \Big|_{\zeta = 1} -  \tilde c \Big|_{\zeta = 0} =  \tilde a \Big|_{\zeta = 1}(h+Z)_x - \tilde a \Big|_{\zeta = 0}Z_x + \tilde b \Big|_{\zeta = 1}(h + Z)_y - \tilde b \Big|_{\zeta = 0} Z_y - \p{h \int_0^1 \tilde a \ d\zeta}_x - \p{h \int_0^1 \tilde b \ d\zeta}_y.
$$
Applying the zero outward normal boundary conditions in \eqref{eq:bbc}, we recover the divergence-free condition of the magnetic field in the MRSW equations:
\begin{equation}\label{eq:sw_divb}
	\p{ha_m}_x + \p{hb_m}_y = 0,
\end{equation}
where $\displaystyle a_m = \int_0^1 \tilde a \ d\zeta$ and $\displaystyle b_m = \int_0^1 \tilde b \ d\zeta$ are the mean magnetic fields in the $x$- and $y$-directions, respectively.

\subsubsection{Momentum balance mapping}\label{sec3.1.3}
Starting with the momentum balance \eqref{eq:mo_x}, we again multiply by $h$ and apply the differential operator mappings in \eqref{eq:mappings} to get
\begin{equation}\label{eq:map_mo_x}
\begin{aligned}
	&\p{h \tilde u}_t +\br{h\tilde u^2 +h\p{ \frac{1}{\rho}\tilde p 
	+ \frac{1}{2}\abs{\tilde{\bm b}}^2}- h \tilde a^2}_x +\br{h \tilde u \tilde v - h \tilde a \tilde b}_y 
	+ \Bigg[-\p{\frac{1}{\rho}\tilde p + \frac{1}{2}\abs{\tilde{\bm b}}^2}\p{\zeta h + Z}_x\\[4pt]
	&+\tilde u \br{\tilde w - \p{\zeta h + Z}_x \tilde u - \p{\zeta h + Z}_y \tilde v  - \p{\zeta h + Z}_t}
	- \tilde a \br{\tilde c - \p{\zeta h + Z}_x \tilde a - \p{\zeta h + Z}_y \tilde b}\Bigg]_\zeta
	= fh\tilde v.
\end{aligned}
\end{equation}
Similarly, we also transform the magneto-hydrostatic pressure relation \eqref{eq:MHpressure} to obtain
\begin{equation*}\label{eq:mapped_mhp}
	{ \tilde p + \frac{1}{2}\rho\abs{{\tilde{\bm b}}}^2} = (1-\zeta) \rho g h,
\end{equation*}
which allows for the following simplification of the pressure dependent terms in \eqref{eq:map_mo_x}:
\begin{equation}\label{eq:sw_mhs}
	\br{h\p{\frac{1}{\rho}\tilde p + \frac{1}{2}\abs{\tilde{\bm b}}^2}}_x 
	- \br{\p{\frac{1}{\rho}\tilde p + \frac{1}{2}\abs{\tilde{\bm b}}^2}\p{\zeta h + Z}_x}_\zeta 
	= \p{\frac{1}{2}gh^2}_x + ghZ_x.
\end{equation}
Combining \eqref{eq:map_mo_x} with \eqref{eq:sw_mhs}, and defining the vertical coupling operators for the velocity and magnetic field as 
\begin{align}
	h\omega [h, \tilde u, \tilde v] & = \tilde w - \p{\zeta h + Z}_x \tilde u 
	- \p{\zeta h + Z}_y \tilde v  - \p{\zeta h + Z}_t, \label{eq:omega}\\[4pt]
	h \mc C [h, \tilde a, \tilde b] & = \tilde c - \p{\zeta h + Z}_x \tilde a 
	- \p{\zeta h + Z}_y \tilde b, \label{eq:mcC}
\end{align}
respectively, we recover the mapped horizontal momentum balance in the $x$-direction:
\begin{equation}\label{eq:sw_hu}
	\p{h \tilde u}_t + \p{h \tilde u^2 + \frac{1}{2}gh^2 - h \tilde a^2}_x 
	+\p{h \tilde u \tilde v - h \tilde a \tilde b}_y 
	+\br{h \tilde u \omega[h,\tilde u, \tilde v] - h\tilde a\mc C [h, \tilde a, \tilde b]}_\zeta
	= fhv-ghZ_x.
\end{equation}
The $y$-direction momentum balance is formulated following the same rationale; see equation \eqref{eq:newhv}. 

\begin{rmk}\label{rmk3.1}
Extending that in \cite{Kowalski2019Moment}, we can combine the integral forms of \eqref{eq:w_tilde} and \eqref{eq:c_tilde} with \eqref{eq:omega} and \eqref{eq:mcC}, respectively, to obtain vertical coupling terms that no longer have dependence on the $z$-direction velocity or magnetic field:
\begin{align}
	h \omega [h, \tilde u, \tilde v] &= - \p{h \int_0^\zeta \tilde u - u_m \ d\zeta}_x 
	- \p{h \int_0^\zeta \tilde v - v_m \ d\zeta}_y, \label{eq:omega_new}\\[4pt]
	h \mc C [h, \tilde a, \tilde b] &= - \p{h \int_0^\zeta \tilde a \ d\zeta}_x 
	- \p{h \int_0^\zeta \tilde b\ d\zeta}_y. \label{eq:mcC_new}
\end{align}
Notice that \eqref{eq:omega_new} and \eqref{eq:mcC_new} also imply the following divergence relations:
\begin{align}
	&(h \tilde u)_x + (h \tilde v)_y + (h \omega)_\zeta = (h u_m)_x + (h v_m)_y, 
	\label{eq:omega_rel}\\[4pt]
	&(h \tilde a)_x + (h \tilde b)_y + (h \mc C)_\zeta = 0, \label{eq:mcC_rel}
\end{align}
which are listed here for future reference; see \S\ref{sec5.2} for their appearance within the numerical method.
In addition, notice that through the mapped divergence-free condition of the magnetic field in \eqref{eq:sw_divb}, we may equivalently write the magnetic vertical coupling term $h\mc C$ using the subtraction of the mapped and mean magnetic fields---providing an identical form to that of the velocity vertical coupling \eqref{eq:omega_new}.
We, however, will not use this alternative setup of \eqref{eq:mcC_new}, since keeping it as seen will benefit us in the discretization; see \S\ref{sec5.2} for further details.
\end{rmk}

\subsubsection{Mapping of Faraday's equation}\label{sec3.1.4}
We follow the same steps for Faraday's equation in \eqref{eq:mag_x}; that is, we multiply by $h$, apply the differential operator transformations of \eqref{eq:mappings}, and rearrange to obtain

$$
	(h\tilde a)_t + \p{h \tilde a \tilde v - h \tilde b \tilde u}_y 
	+ \p{\tilde a \br{\tilde w - (\zeta h + Z)_t - (\zeta h +Z)_y\tilde v} 
	- \tilde u \br{\tilde c - (\zeta h + Z)_y \tilde b}}_\zeta = 0.
$$
This can be simplified by substituting in the definitions of the vertical coupling terms $h \omega [h, \tilde u, \tilde v]$ and $h \mc C [h, \tilde a, \tilde b]$ from \eqref{eq:omega} and \eqref{eq:mcC}, respectively, to get
$$
	(h\tilde a)_t + \p{h \tilde a \tilde v - h \tilde b \tilde u}_y 
	+ \p{\tilde a \big[h \omega [h, \tilde u, \tilde v] +(\zeta h +Z)_x\tilde u\big] 
	- \tilde u \br{h \mc C [h, \tilde a, \tilde b] + (\zeta h + Z)_x \tilde a}}_\zeta = 0. 
$$
This further reduces to the mapping of the magnetic field equation in the $x$-direction:
\begin{equation}\label{eq:sw_ha}
	(h\tilde a)_t + \p{h \tilde a \tilde v - h \tilde b \tilde u}_y 
	+ \p{h\tilde a \omega [h, \tilde u, \tilde v] -  h \tilde u \mc C [h, \tilde a, \tilde b] }_\zeta = 0.
\end{equation}
The formulation for the $y$-direction magnetic field follows analogously; see equation \eqref{eq:newhb}. 


\subsection{Complete reference system}\label{sec3.2}
Combining the resulting mapped equations from \S \ref{sec3.1.1}--\ref{sec3.1.4}, we recover the complete, vertically resolved MRSW reference system, which reads:
\begin{align}
	& h_t + \p{hu_m}_x + \p{hv_m}_y = 0, \label{eq:newh}\\[4pt]
	& \p{ha_m}_x + \p{hb_m}_y = 0, \label{eq:newdivb}\\[4pt]
	& \p{h \tilde u}_t + \p{h \tilde u^2 + \frac{1}{2}gh^2 - h \tilde a^2}_x 
	+\p{h \tilde u \tilde v - h \tilde a \tilde b}_y 
	+\br{h \tilde u \omega - h\tilde a\mc C}_\zeta
	= fh\tilde v-ghZ_x, \label{eq:newhu}\\[4pt]
	& \p{h \tilde v}_t +\p{h \tilde u \tilde v - h \tilde a \tilde b}_x +\p{h \tilde v^2 
	+ \frac{1}{2}gh^2 - h \tilde b^2}_y
	+\br{h \tilde v \omega - h\tilde b\mc C }_\zeta
	= -fh\tilde u-ghZ_y, \label{eq:newhv}\\[4pt]
	& (h\tilde a)_t + \p{h \tilde a \tilde v - h \tilde b \tilde u}_y + \p{h\tilde a \omega
	 -  h \tilde u \mc C  }_\zeta = 0,\label{eq:newha}\\[4pt]
	& (h\tilde b)_t + \p{h \tilde b \tilde u - h \tilde a \tilde v}_x 
	+ \p{h\tilde b \omega -  h \tilde v \mc C }_\zeta = 0, 
	\label{eq:newhb}
\end{align}
where $\omega$ and $\mc C$ are the vertical coupling terms defined in \eqref{eq:omega_new} and \eqref{eq:mcC_new}, respectively, and 
$$
    \psi_m(t,x,y) = \int_0^1 \tilde{\psi}\ d\zeta,
$$
for $\psi = u, v, a,$ and $b$.

Notice that the system \eqref{eq:newh}--\eqref{eq:newhb} reduces to the MRSW equations (see, e.g., \S\ref{sec4.4.1}) in the case of constant vertical profiles of velocity and magnetic field. 
In addition, this system no longer has dependence on the vertical velocity $w$ and magnetic field $c$ from \eqref{eq:mo_x}--\eqref{eq:bbc} thanks to the definitions of the vertical coupling in \eqref{eq:omega_new}--\eqref{eq:mcC_new}. 
Furthermore, this reduction of total number of variables within the MRSW reference system implies \eqref{eq:newh}--\eqref{eq:newhb} is more computationally friendly in comparison to the original system \eqref{eq:mo_x}--\eqref{eq:bbc}, but still fully 3-D; details on the discretization of the MRSW reference system can be found in \S\ref{sec5.2}.

\section{A moment closure for magnetic rotating shallow flow}\label{sec4}
Standard shallow flow models---in this paper the most relevant being the MRSW equations---spatially reduce the mapped reference system in \eqref{eq:newh}--\eqref{eq:newhb} from 3-D to 2-D by replacing the velocities and magnetic fields with their depth-averaged mean values $u_m$, $v_m$, $a_m$, and $b_m$. 
This comes with the potential loss of information regarding the corresponding vertical profiles. 
In this section, we also reduce the reference system to 2-D, but to mitigate this vertical information loss, we follow the work of \cite{Kowalski2019Moment} and express variables with vertical profiles through a polynomial expansion:
\begin{equation}\label{eq:expan}
	\widetilde{\bm{W}}(t, x, y, \zeta) = \bm{W}_m(t,x,y) + \sum_{\ell = 1}^M \bm{\Psi}_\ell(t,x,y)\phi_\ell\left(\zeta\right),
\end{equation}
where $\widetilde{\bm{W}} = (\tilde u,\tilde v,\tilde a,\tilde b)^\top$, $\bm{W}_m = (u_m,v_m,a_m,b_m)^\top$, $\bm{\Psi}_\ell = (\alpha_\ell,\beta_\ell,\gamma_\ell,\eta_\ell)^\top$ are the vertical profile polynomial coefficients, and $\phi_\ell$ are basis functions on the domain $[0,1]$ for $\ell = 1,\dots, M$. 
Following \cite{Kowalski2019Moment}, we define $\phi_\ell$ using the scaled (orthogonal) Legendre polynomials on the domain $[0,1]$ normalized such that $\phi_\ell(0) = 1$.
For example, the first three basis polynomials read 
\begin{equation}\label{eq:basis}
	\phi_1(\zeta) = 1-2\zeta, \quad 
	\phi_2(\zeta) = 1-6\zeta+6\zeta^2, \quad 
	\phi_1(\zeta) = 1-12\zeta + 30\zeta^2 - 20\zeta^3.
\end{equation}
The expansions in \eqref{eq:expan} imply that the vertical profiles of velocity and magnetic field are now modeled by the spatially 2-D functions in $\bm{W}_m(t,x,y)$, and their corresponding 2-D polynomial coefficients within $\bm{\Psi}_\ell$ for $\ell = 1,\dots, M$.
Furthermore, it is clear that as the number $M$ increases, the approximation of the vertical profile will improve. 
Additionally, note that the depth-averaged quantities $\bm W_m$ are essentially an expansion of the zeroth-order Legendre polynomial $\phi_0(\zeta) = 1$, thus matching the single-variable depth-averaged vertical profiles of standard shallow flow approximations.

By \eqref{eq:expan}, ${\bm{W}}_m$ and $\bm \Psi_\ell$ together provide an approximate 3-D rendering of $\widetilde{\bm W}$. 
Thus, the remainder of this section presents the derivation of the evolution equations for the 2-D variables within ${\bm{W}}_m$ and $\bm \Psi_\ell$.
In doing so, this will achieve the overall goal of reducing \eqref{eq:newh}--\eqref{eq:newhb} from 3-D to an approximate system dependent on only two spatial dimensions.
The presented derivation closely parallels that of \cite{Kowalski2019Moment}, where the vertically resolved reference system for the shallow water equations is derived together with their corresponding moment approximations---except we must now consider contributions of the magnetic field. 
Even so, the magnetic contributions interestingly parallel that of the velocity, as we will show in the remainder of Section \ref{sec4}. 
For the remainder of this paper, the tildes seen in \eqref{eq:newh}--\eqref{eq:newhb} will be dropped for readability.

We will first successively derive the evolution equations for the mean variables $u_m,v_m,a_m,b_m$, before deriving the remaining evolution equations for the additional coefficients $\alpha_\ell,\beta_\ell,\gamma_\ell,\eta_\ell$ for $\ell = 1,\dots, M$.

\subsection{Mean momentum and magnetic field evolution derivation}\label{sec4.1}
To derive the mean momentum and magnetic field evolution equations, we take the equations \eqref{eq:newhu}--\eqref{eq:newhb} and (i) integrate over $\zeta$ from zero to one, and (ii) substitute in the polynomial expansions of \eqref{eq:expan} and simplify.
Starting with the $x$-direction momentum equation \eqref{eq:newhu}, and integrating over $\zeta$, we obtain
\begin{align*}
	\p{hu_m}_t + \br{h \int_0^1 u^2 \ d \zeta + \frac{1}{2}g h^2 
	- h \int_0^1 a^2 \ d \zeta}_x
	&+  \br{h \int_0^1 uv \ d \zeta - h \int_0^1 ab \ d \zeta}_y \\[4pt]
	& + \int_0^1h  u \omega[h, u, v] - h a\mc C [h, a, b]\ d\zeta = fhv_m-gh Z_x.
\end{align*}
We then wish to implement the polynomial expansions from \eqref{eq:expan}. To do so, we must simplify the higher order integrals seen above. 
Using the orthogonality and properties of Legendre polynomials, we find 
\begin{equation}\label{eq:ints}
	\int_0^1 W^{(j)}W^{(k)} \ d\zeta = W_m^{(j)} W^{(k)}_m + \sum_{\ell = 1}^M 
	\frac{\Psi_\ell^{(j)} \Psi_\ell^{(k)}}{2\ell + 1}, 
\end{equation}
where $W^{(j)}$, $W_m^{(j)}$, and $\Psi_\ell^{(j)}$ denote the $j$th component of $\bm W$, $\bm W_m$, and $\bm\Psi_\ell$ in \eqref{eq:expan}, respectively. 
Then, by using the definitions of the mean velocities and magnetic fields, and the corresponding vertical coupling terms in \eqref{eq:omega_new}--\eqref{eq:mcC_new}, we can see that for an arbitrary function $\psi$, we have 
\begin{equation}\label{eq:int_coupling}
\begin{aligned}
	&\int_0^1 \p{\psi h\omega[h,u,v]}_\zeta \ d \zeta = 0, \\[4pt]
	&\int_0^1 \p{\psi h\mc C[h,a,b]}_\zeta \ d \zeta 
	= \psi \Big|_{\zeta = 1}\br{(ha_m)_x + (hb_m)_y}.
\end{aligned}
\end{equation}
Note that the integral of the magnetic vertical coupling term multiplied by an arbitrary function $\psi$ also equals zero through the divergence-free condition \eqref{eq:newdivb}, but keeping it as written above proves beneficial in the discretization; see further discussion in \S\ref{sec4.3}.

Combining the simplifications in \eqref{eq:ints}--\eqref{eq:int_coupling} with the $\zeta$-integrated $x$-direction momentum equation, we recover the evolution equation for $hu_m$ purely in 2-D: 
\begin{equation}
\begin{aligned}
	\p{hu_m}_t +& \br{h \p{u_m^2 + \sum_{\ell = 1}^M\frac{\alpha_\ell^2}{2 \ell +1}} 
	+ \frac{1}{2}g h^2 - h \p{a_m^2 + \sum_{\ell = 1}^M
	\frac{\gamma_\ell^2}{2 \ell +1}}}_x\\[4pt]
	+&  \br{h \p{u_mv_m + \sum_{\ell = 1}^M\frac{\alpha_\ell\beta_\ell}{2 \ell +1}} 
	- h \p{a_mb_m + \sum_{\ell = 1}^M\frac{\gamma_\ell\eta_\ell}{2 \ell +1}}}_y \\[4pt]
	& \hspace{6.5cm}= fhv_m-ghZ_x-a\Big|_{\zeta = 1} \br{\p{ha_m}_x + \p{hb_m}_y}.
\end{aligned}
\label{eq:hu_tmp}
\end{equation}
The derivations of the equations for $hv_m,\ ha_m,$ and $hb_m$ follow the same procedure; for the sake of brevity, we refrain from presenting these evolution equations until the final system in \S\ref{sec4.3}.

\subsection{Vertical coefficient evolution equations}\label{sec4.2}
To obtain the full 3-D approximation to the velocity and magnetic fields as seen in \eqref{eq:expan}, we need expressions for the vertical profile coefficients $\bm{\Psi}_i = (\alpha_i,\beta_i,\gamma_i,\eta_i)^\top$; in this work, these come in the form of evolution equations.
To recover these evolution equations, we again follow \cite{Kowalski2019Moment} and take moments of the velocity and magnetic fields with respect to the Legendre polynomials.
By orthogonality of the selected basis functions $\phi_i(\zeta)$, doing so results in 
\begin{equation}\label{eq:moment}
	\int_0^1 \phi_i \bm W \ d \zeta = \frac{\bm \Psi_i}{2i+1}, 
\end{equation}
where $\bm W$ is defined in \eqref{eq:expan}.
Therefore, by multiplying the reference velocity and magnetic field equations \eqref{eq:newhu}--\eqref{eq:newhb} by $\phi_i(\zeta)$ and integrating over $\zeta \in [0,1]$, we will obtain evolution equations for the vertical profile coefficients $\bm{\Psi}_i = (\alpha_i,\beta_i,\gamma_i,\eta_i)^\top$.
Due to the similarities of this process with the computation of moments, we will also refer to these as `moment' equations. 
We present this derivation for the $x$-velocity moments $\alpha_i$ in this section, and the moment equations for $\beta_i,\ \gamma_i,$ and $\eta_i$ are derived analogously. 

Multiplying the $x$-momentum reference equation \eqref{eq:newhu} by $\phi_i(\zeta)$ and integrating over $[0,1]$, we obtain
\begin{equation}\label{eq:int_hu}
\begin{aligned}
	&\p{h\int_0^1 \phi_i u\ d \zeta}_t 
	+ \br{h\int_0^1 \phi_i u^2 \ d \zeta - h\int_0^1 \phi_i a^2 \ d \zeta}_x
	+\br{h\int_0^1 \phi_i u v\ d \zeta - h\int_0^1 \phi_i
	ab\ d\zeta}_y\\[4pt]
	&\hspace{8.5cm}+\int_0^1 \phi_i \p{hu\omega - ha\mc C}_\zeta \ d \zeta 
	= \int_0^1 fh\phi_iv\ d\zeta.
\end{aligned}
\end{equation}
Note that here we already applied orthogonality to the $\frac{1}{2}gh^2$ term in the $x$-direction flux and to the bottom topography source term.
The term inside of the time derivative and the remaining source term can be simplified using \eqref{eq:moment}. 
In addition, all integrals residing in the flux terms of \eqref{eq:int_hu} have an integrand of the form $ \phi_iW^{(j)}W^{(k)}$, where $W^{(j)}$ again denotes the $j$th component of $\bm W$ in \eqref{eq:expan}. Therefore, all flux term integrals can be simplified by substituting in the vertical expansions of \eqref{eq:expan} and evaluating, resulting in
\begin{equation*}\label{eq:int_uv}
	\int_0^1 \phi_i W^{(j)}W^{(k)} \ d \zeta = \frac{1}{2i+1}\p{W^{(j)}_m \Psi^{(k)}_i 
	+ W^{(k)}_m \Psi^{(j)}_i 
	+\sum_{\ell,n = 1}^M A_{i\ell n}\Psi^{(j)}_\ell \Psi^{(k)}_n},
\end{equation*}
where 
$$
	A_{i\ell n} = (2i+1)\int_0^1 \phi_i \phi_\ell\phi_n\ d \zeta.
$$
Lastly, we must evaluate the integrals in \eqref{eq:int_hu} associated with the vertical coupling. 
We present these computations in Appendix \ref{secA2}. 
Using these results presented in \eqref{eq:int_omega} and \eqref{eq:int_mcC}, we obtain the following evolution equation for the $x$-momentum vertical profile coefficients $\alpha_i$, $i = 1,\dots, M$:
\begin{equation}
\begin{aligned}
	&\p{h \alpha_i}_t  + \br{h \p{2u_m \alpha_i+ \sum_{\ell,n=1}^M
	A_{i\ell n}\alpha_\ell \alpha_n}
	-h\p{2a_m \gamma_i+\sum_{\ell,n=1}^M A_{i\ell n}
	\gamma_\ell \gamma_n}}_x\\[4pt]
	&\hspace{1cm} +\br{h \p{u_m \beta_i + v_m \alpha_i+ \sum_{\ell,n=1}^M 
	A_{i\ell n}\alpha_\ell \beta_n}
	-h\p{a_m \eta_i+b_m\gamma_i+\sum_{\ell,n=1}^M 
	A_{i\ell n}\gamma_\ell \eta_n}}_y\\[4pt]
	&= fh\beta_i + u_m \br{(h\alpha_i)_x + (h\beta_i)_y} - \sum_{\ell,n=1}^M
	B_{i\ell n} \alpha_n \br{(h\alpha_\ell)_x + (h\beta_\ell)_y} \\[4pt]
	&\ \  -a_m \br{(h\gamma_i)_x + (h\eta_i)_y} 
	+ \sum_{\ell,n=1}^M  B_{i \ell n}\gamma_n 
	\br{(h\gamma_\ell)_x + (h\eta_\ell)_y} -\br{\p{ha_m}_x 
	+ \p{hb_m}_y} \br{\sum_{\ell=1}^M \gamma_\ell \p{1 + \Gamma_{i \ell}}},
\end{aligned}
\label{eq:h_alpha1}
\end{equation}
where $B_{i\ell n}$ is defined in \eqref{eq:B_iln}, and 
$$
\Gamma_{i\ell}  = (2i+1) \int_0^1\zeta\phi_i \phi_\ell'\ d \zeta.
$$
The moment equations for $\beta_i$, $\gamma_i$, and $\eta_i$, $i = 1, \dots, M$, are derived in an analogous manner, and their evolution equations are presented in \S\ref{sec4.3}.

\subsection{The arbitrary-order magnetic rotating shallow water moment equations}\label{sec4.3}
Using the derivation techniques presented in \S\ref{sec4.1} and \S\ref{sec4.2}, we obtain the full, arbitrary $M$th-order Magnetic Rotating Shallow Water Moment Equations (MRSWME):
\begin{align}
	& h_t + \p{hu_m}_x + \p{hv_m}_y = 0, \label{eq:fin_h}\\[4pt]
	& \p{ha_m}_x + \p{hb_m}_y = 0, \label{eq:fin_divb}\\[4pt]
	& \begin{aligned}
		\p{hu_m}_t +& \br{h u_m^2 + \frac{1}{2}g h^2 - ha_m^2 + 
		h\sum_{\ell = 1}^M\frac{\alpha_\ell^2-\gamma_\ell^2}{2 \ell +1} }_x
		+ \br{h u_mv_m - ha_mb_m + h\sum_{\ell = 1}^M
		\frac{\alpha_\ell\beta_\ell-\gamma_\ell\eta_\ell}{2 \ell +1} }_y \\[4pt]
		& \hspace{7.5cm}= fhv_m-ghZ_x-a\Big|_{\zeta = 1} \br{\p{ha_m}_x + \p{hb_m}_y}, 
	\end{aligned}\label{eq:hu_poly}\\[4pt]
	& \begin{aligned}
		\p{hv_m}_t +&  \br{h u_mv_m - ha_mb_m + h\sum_{\ell = 1}^M
		\frac{\alpha_\ell\beta_\ell-\gamma_\ell\eta_\ell}{2 \ell +1} }_x
		+ \br{h v_m^2 + \frac{1}{2}g h^2 - hb_m^2 + 
		h\sum_{\ell = 1}^M\frac{\beta_\ell^2-\eta_\ell^2}{2 \ell +1}}_y\\[4pt]
		& \hspace{7.5cm}= -fhu_m- ghZ_y-b\Big|_{\zeta = 1} \br{\p{ha_m}_x + \p{hb_m}_y},
	\end{aligned}\label{eq:hv_poly}\\[4pt]
	& \p{ha_m}_t + \br{ha_mv_m - h b_mu_m + h\sum_{\ell = 1}^M
	\frac{\beta_\ell\gamma_\ell-\alpha_\ell\eta_\ell}{2 \ell +1} }_y 
	=  -u\Big|_{\zeta = 1} \br{\p{ha_m}_x + \p{hb_m}_y}, \label{eq:poly_ha}\\[4pt]
	& \p{hb_m}_t + \br{h b_mu_m - ha_mv_m + h\sum_{\ell = 1}^M
	\frac{\alpha_\ell\eta_\ell-\beta_\ell\gamma_\ell}{2 \ell +1} }_x 
	=  -v\Big|_{\zeta = 1} \br{\p{ha_m}_x + \p{hb_m}_y}, \label{eq:poly_hb}\\[4pt]
	& \begin{aligned}
		\p{h \alpha_i}_t  &+ \br{2 \p{hu_m \alpha_i-ha_m \gamma_i}+h\sum_{\ell,n=1}^M
		A_{i\ell n}\p{\alpha_\ell \alpha_n-\gamma_\ell \gamma_n} }_x\\[4pt]
		& +\br{h u_m \beta_i + hv_m \alpha_i-ha_m \eta_i-hb_m\gamma_i
		+h\sum_{\ell,n=1}^M A_{i\ell n}\p{\alpha_\ell \beta_n
		-\gamma_\ell \eta_n}}_y\\[4pt]
		&\quad= fh\beta_i + u_m D_i -a_m \mc D_i -\br{\p{ha_m}_x 
		+ \p{hb_m}_y} {\sum_{\ell=1}^M \gamma_\ell \p{1 + \Gamma_{i \ell}}}
		- \sum_{\ell,n=1}^M B_{i\ell n} \p{ \alpha_n D_\ell - \gamma_n\mc D_\ell},
	\end{aligned}\label{eq:h_alpha} \\[4pt]
	& \begin{aligned}
		\p{h \beta_i}_t &+\br{h u_m \beta_i + hv_m \alpha_i-ha_m \eta_i-hb_m\gamma_i
		+h\sum_{\ell,n=1}^M A_{i\ell n}\p{\alpha_\ell \beta_n
		-\gamma_\ell \eta_n}}_x\\[4pt]
		&+ \br{2 \p{hv_m \beta_i-hb_m \eta_i}+h\sum_{\ell,n=1}^M
		A_{i\ell n}\p{\beta_\ell \beta_n-\eta_\ell \eta_n}}_y\\[4pt]
		&\quad= -fh\alpha_i +v_m D_i-b_m\mc D_i 
		-\br{\p{ha_m}_x + \p{hb_m}_y} {\sum_{\ell=1}^M \eta_\ell 
		\p{1 + \Gamma_{i \ell}}}
		- \sum_{\ell,n=1}^M B_{i\ell n} \p{\beta_n D_\ell -\eta_n\mc D_\ell},
	\end{aligned}\label{eq:h_beta} \\[4pt]
	& \begin{aligned}
		\p{h \gamma_i}_t &+\br{h a_m \beta_i + hv_m \gamma_i -h b_m\alpha_i - h u_m\eta_i
		+h \sum_{\ell,n=1}^M A_{i\ell n}\p{\beta_\ell\gamma_n 
		- \alpha_\ell\eta_n}}_y\\[4pt]
		& = a_m D_i -u_m \mc D_i 
		-\br{\p{ha_m}_x + \p{hb_m}_y}{\sum_{\ell=1}^M \alpha_\ell \p{1 + \Gamma_{i \ell}}}
		- \sum_{\ell,n=1}^M B_{i\ell n}\p{ \gamma_n D_\ell -\alpha_n \mc D_n },
	\end{aligned} \label{eq:h_gamma} \\[4pt]
	& \begin{aligned}
		\p{h \eta_i}_t &+\br{hb_m \alpha_i+hu_m\eta_i-h a_m \beta_i -h v_m \gamma_i
		+h\sum_{\ell,n=1}^M A_{i\ell n}\p{\alpha_\ell\eta_n
		-\beta_\ell\gamma_n}}_x\\[4pt]
		& = b_m D_i -v_m \mc D_i 
		-\br{\p{ha_m}_x + \p{hb_m}_y}{\sum_{\ell=1}^M \beta_\ell \p{1 + \Gamma_{i \ell}}}
		- \sum_{\ell,n=1}^M B_{i\ell n}\p{ \eta_n D_\ell -\beta_n \mc D_n },
	\end{aligned} \label{eq:h_eta}
\end{align}
where $i = 1,\dots M$, the notation $\displaystyle u\Big|_{\zeta = 1} = u_m + \sum_{\ell = 1}^M \alpha_\ell\phi_\ell(1)$, and similarly for $v,\ a,$ and $b$ evaluated at $\zeta = 1$, and
\begin{equation}
\begin{aligned}
	&A_{i\ell n} = (2i+1)\int_0^1 \phi_i \phi_\ell\phi_n\ d \zeta,\qquad 
	B_{i\ell n} = (2i+1) \int_0^1 \phi_i' 
	\p{\int_0^\zeta \phi_\ell \d \hat\zeta}\phi_n \d \zeta,\\[4pt]
	&\Gamma_{i\ell}  = (2i+1) \int_0^1\zeta\phi_i \phi_\ell'\ d \zeta, \qquad
	D_i = \p{h \alpha_i}_x + \p{h \beta_i}_y,\qquad 
	\mc D_i = \p{h \gamma_i}_x + \p{h \eta_i}_y.
\end{aligned}
\label{eq:Coefs}
\end{equation}
The complete $M$th-order MRSWME above is a system of $4M+5$ evolution equations in addition to the divergence-free condition of the magnetic field.
In the case that vertical deviations from the depth-averaged value are negligible, \eqref{eq:fin_h}--\eqref{eq:Coefs} simplifies to the the MRSW system often used in geophysical studies, and it is clear that $\alpha_i,\ \beta_i,\ \gamma_i,$ and $\eta_i$ will remain zero for all time if all moments are zero initially.
In addition, the Coriolis force plays a role in the velocity moment equations for $\alpha_i$ and $\beta_i$, implying the rotational source influences the vertical profile. 
This is unlike the gravitational forces (magneto-hydrostatic pressure and bottom topography source), which do not enter the moment evolution equations. 

Within the derivation of the MRSWME, the structure of the $\gamma_\ell, \eta_\ell$ pieces in the fluxes of \eqref{eq:fin_h}--\eqref{eq:h_eta} significantly parallel those of $\alpha_\ell,\beta_\ell$. 
Even the vertical coupling contribution of the magnetic field in \eqref{eq:int_mcC} is identical in form to that of velocity vertical coupling contribution in \eqref{eq:int_omega}, if the zero-divergence magnetic field equation \eqref{eq:fin_divb} is applied.
However, in opting to keep this term (and all other source terms scaled by $[(ha_m)_x+(hb_m)_y]$), we see the appearance of the so-called Godunov-Powell source; see, e.g.,
\cite{Godunov1972symmetric,
Powell1995Upwind,
Powell1999Sol}; 
which, much like that of MHD or MRSW, helps recover a missing eigenvector of the MRSWME. 
This additional source has proven quite beneficial in numerical approximations of MHD systems; see, e.g.,
\cite{Chandrashekar2016Entropy,
Powell1995Upwind,
Powell1999Sol,
Winters2017Uniquely};
and is a key component of the discretizations presented in \S\ref{sec5.1} and \S\ref{sec5.2}.

As shown for the shallow water moment equations in \cite{Kowalski2019Moment}, one downside of the moment approximation is the loss of global hyperbolicity for $M \geq 2$.
We discuss how this takes effect in the magnetic extension in \S\ref{sec4.4.2}.

\subsection{MRSWME examples}\label{sec4.4}
In this section, we present two examples of the model hierarchy shown in \eqref{eq:fin_h}--\eqref{eq:h_eta}.
For the sake of simplicity, we only present the one-dimensional system in the $y$-direction. 
Note that, due to the Coriolis force further coupling the $u$ and $v$ equations, we may not simply cancel out the $x$-direction velocity and magnetic field---even if $u$ and $a$ are initially zero.

\subsubsection{Zeroth order system: magnetic rotating shallow water equations}\label{sec4.4.1}
Using $M=0$ and thus assuming a constant profile in the vertical direction for velocity and magnetic field, all higher order moment evolutions drop from the system, resulting in the 5-equation model
\begin{equation}\label{eq:zeroorder}
    \begin{pmatrix}
    h \\ 
    hu_m \\ 
    hv_m \\ 
    ha_m \\ 
    hb_m
    \end{pmatrix}_t 
    + 
    \begin{pmatrix}
    hv_m \\ 
    hu_m v_m - ha_m b_m \\ 
    hv_m^2 - hb_m^2 + \frac{1}{2}gh^2 \\ 
    ha_m v_m- hb_m u_m \\ 
    0\end{pmatrix}_y
    = 
    \begin{pmatrix}
    0\\ 
    fhv_m - a_m (hb_m)_y\\ 
    -fhu_m - ghZ_y - b_m (hb_m)_y\\ 
    -u_m (hb_m)_y \\ 
    -v_m (hb_m)_y
    \end{pmatrix},
\end{equation}
paired with the 1-D divergence-free condition $(hb_m)_y = 0$.
These are known as the shallow water MHD or MRSW equations. 
Due to the addition of the Godunov-Powell source term, this system is strongly hyperbolic for all $h >0$. The characteristic wave speeds are given by 
\begin{equation}\label{eq:MRSWspeeds}
	\lambda_1 = v_m, \qquad
	\lambda_{2,3} = v_m \pm \abs{b_m}, \qquad 
	\lambda_{4,5} = v_m \pm \sqrt{b_m^2 + gh},
\end{equation}
and known as the material wave, Alfv\'en waves, and magnetogravity waves, respectively. 
In \S\ref{sec6}, we will see how vertical profiles change the solutions of these traveling waves.

\subsubsection{First order system}\label{sec4.4.2}
Taking $M = 1$ in system \eqref{eq:fin_h}--\eqref{eq:h_eta} implies linear velocity and magnetic field profiles, giving the first-order MRSWME consisting of 9 equations:
\begin{equation}\label{eq:firstorder}
    \begin{pmatrix}
    h \\ 
    hu_m \\ 
    hv_m \\ 
    ha_m \\ 
    hb_m \\ 
    h\alpha_1 \\ 
    h \beta_1 \\ 
    h \gamma_1 \\ 
    h \eta_1
    \end{pmatrix}_t 
    + 
    \begin{pmatrix}
    \textstyle hv_m \\ 
    \textstyle hu_m v_m + \frac{1}{3} h \alpha_1 \beta_1 - ha_m b_m - \frac{1}{3} h \gamma_1 \eta_1 \\ 
    \textstyle hv_m^2 + \frac{1}{3} h \beta_1^2 - hb_m^2 - \frac{1}{3} h \eta_1^2 + \frac{g}{2}h^2 \\ 
    \textstyle ha_m v_m + \frac{1}{3} h \beta_1 \gamma_1 - hb_m u_m - \frac{1}{3} h \alpha_1 \eta_1 \\ 
    0\\ 
    hu_m\beta_1 + h v_m \alpha_1 - ha_m \eta_1 - h b_m \gamma_1 \\ 
    2hv_m\beta_1 - 2h b_m \eta_1 \\ 
    ha_m\beta_1 + h v_m \gamma_1 - hb_m \alpha_1 - h u_m \eta_1 \\ 
    0
    \end{pmatrix}_y
    = 
    \begin{pmatrix}
    0\\ 
    fhv_m - (a_m-\gamma_1) (hb_m)_y\\ 
    -fhu_m - ghZ_y - (b_m-\eta_1) (hb_m)_y\\ 
    -(u_m - \alpha_1) (hb_m)_y \\ 
    -(v_m - \beta_1) (hb_m)_y \\ 
    fh\beta_1 + u_m(h\beta_1)_y - a_m (h \eta_1)_y - 2\gamma_1 (hb_m)_y \\ 
    -fh\alpha_1+v_m(h\beta_1)_y - b_m (h \eta_1)_y - 2\eta_1 (hb_m)_y \\ 
    a_m(h\beta_1)_y - u_m (h \eta_1)_y - 2\alpha_1 (hb_m)_y \\ 
    b_m(h\beta_1)_y - v_m (h \eta_1)_y - 2\beta_1 (hb_m)_y 
    \end{pmatrix},
\end{equation}
which is again paired with the 1-D divergence-free condition $(hb_m)_y = 0$.
Looking into the hyperbolicity of the system, we see that five of the nine wave speeds have a simple closed form:
$$
	\lambda_1 = v_m - \beta_1, \qquad
	\lambda_{2,3} = v_m - b_m \pm \frac{1}{\sqrt{3}}\abs{\beta_1-\eta_1},\qquad
	\lambda_{4,5} = v_m + b_m \pm \frac{1}{\sqrt{3}}\abs{\beta_1+\eta_1}.
$$
The remaining four eigenvalues of the flux Jacobian come from finding the roots $\widetilde{\lambda}$ of the quartic
\begin{equation}\label{eq:quartic}
	\widetilde{\lambda}^4 
	- 3\p{\frac{3b_m^2}{b_m^2 + gh} + 1 + 3\widetilde{\beta}_1^2+\widetilde{\eta}_1^2}\widetilde{\lambda}^2 
	-72b_m\widetilde{\beta}_1\widetilde{\eta}_1\widetilde{\lambda}
	+27b_m^2\p{3-\widetilde{\beta}_1^2-3\widetilde{\eta}_1^2} = 0,
\end{equation}
where 
\begin{equation}\label{eq:q_scalings}
	\widetilde{\lambda} = \lambda_{i} + 3v_m,\qquad
	\widetilde{\beta}_1 = \frac{\beta_1}{\sqrt{b_m^2 + gh}},\qquad
	\widetilde{\eta}_1 = \frac{\eta_1}{\sqrt{b_m^2 + gh}},
\end{equation}
$i = 6,\dots,9$, are defined to increase readability. 

\begin{figure}[H]
\centerline{
\includegraphics[trim=0.9cm 1.3cm 1.4cm 3.7cm, clip, width=8cm]{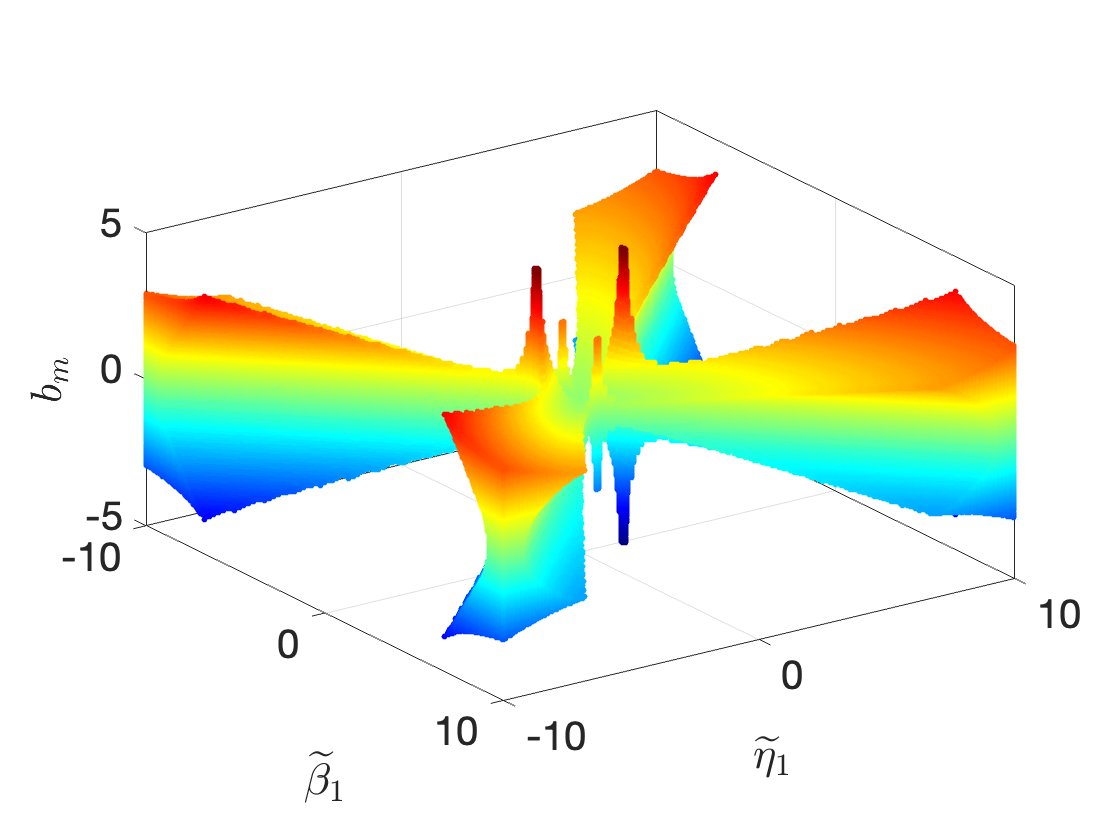}
\includegraphics[trim=0.9cm 1.3cm 1.4cm 3.7cm, clip, width=8cm]{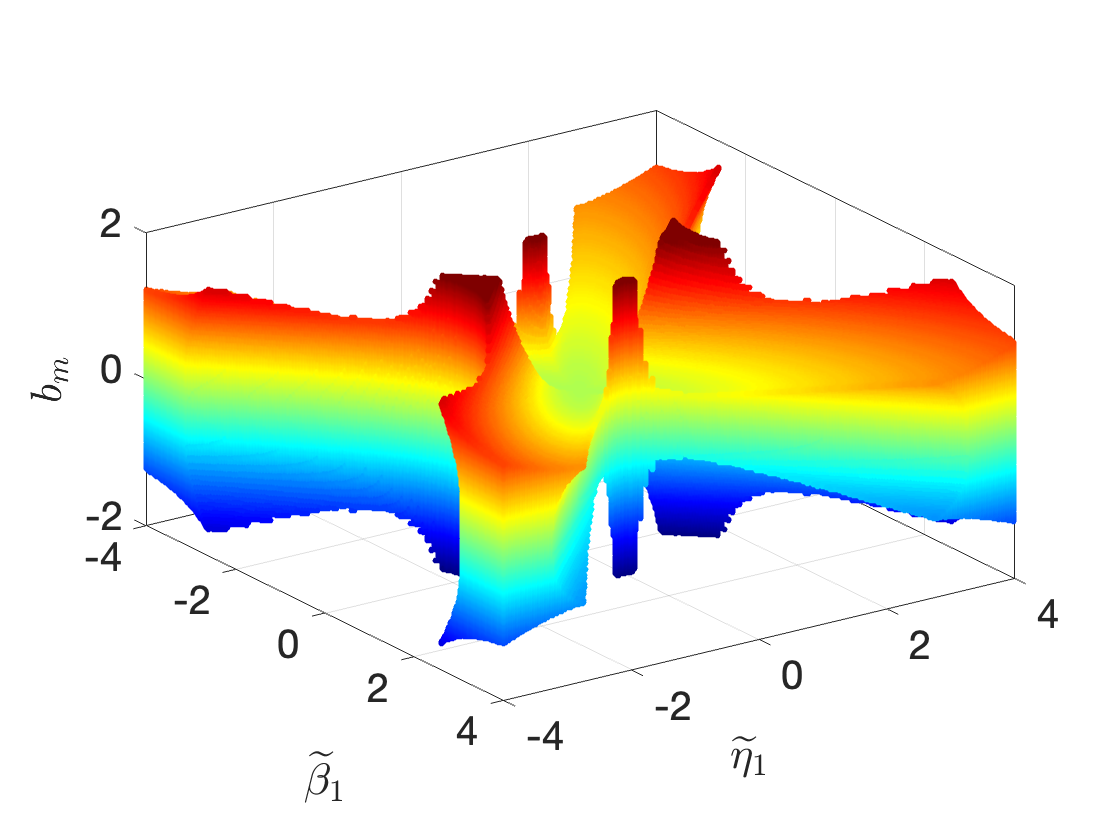}}
\vskip7pt
\centerline{
\includegraphics[trim=1.5cm 3.7cm 3.0cm 6cm, clip, width=8cm]{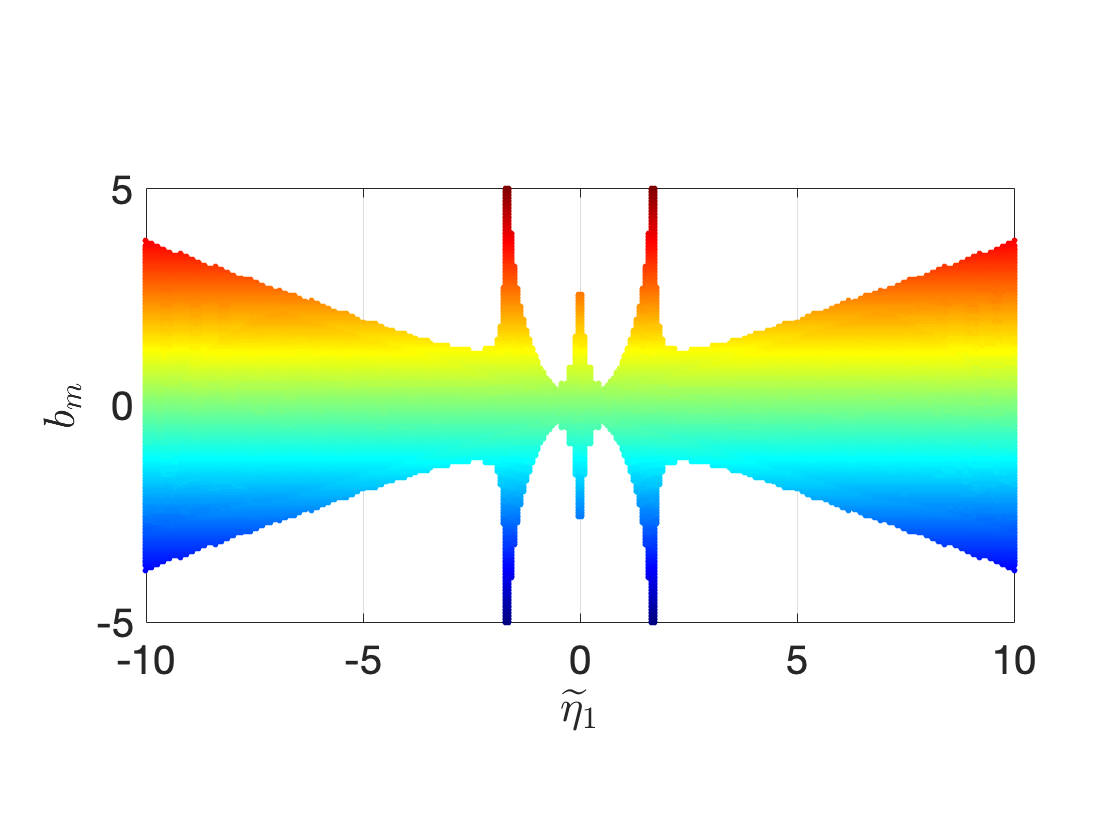}
\includegraphics[trim=1.5cm 3.7cm 3.0cm 6cm, clip, width=8cm]{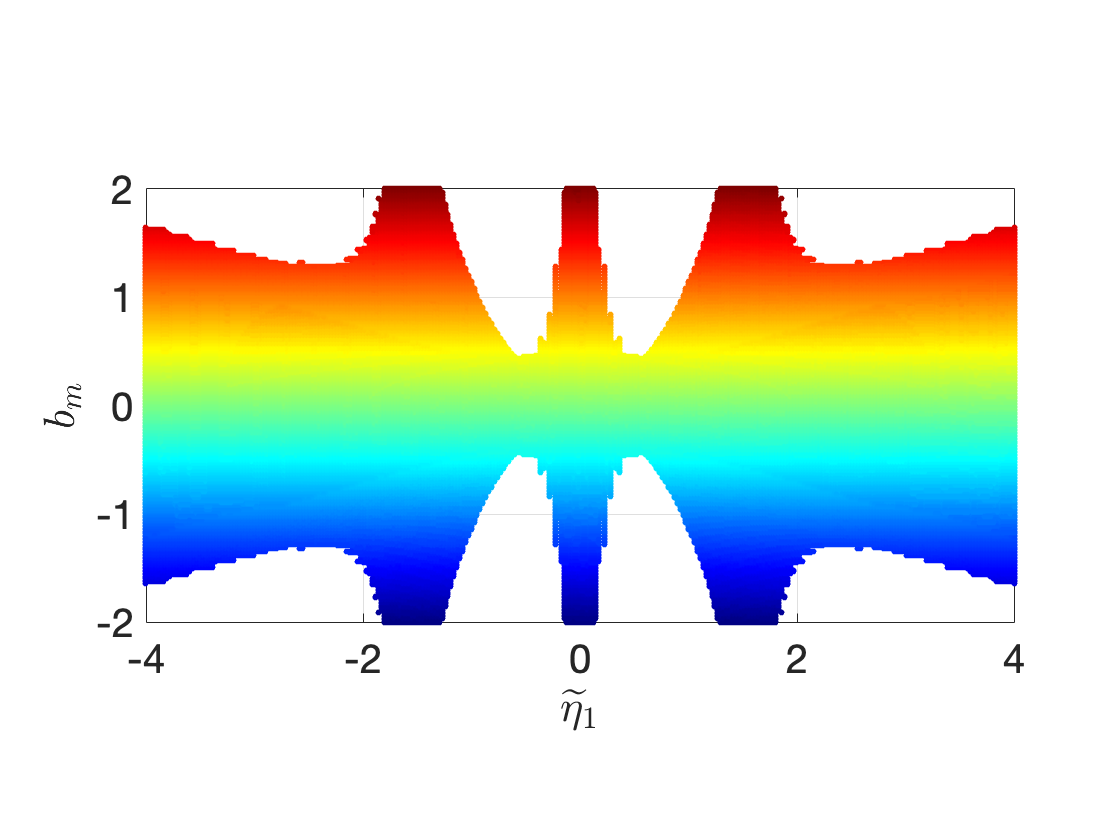}}
\vskip7pt
\centerline{
\includegraphics[trim=4.2cm 0.2cm 6.0cm 1.5cm, clip, width=6cm]{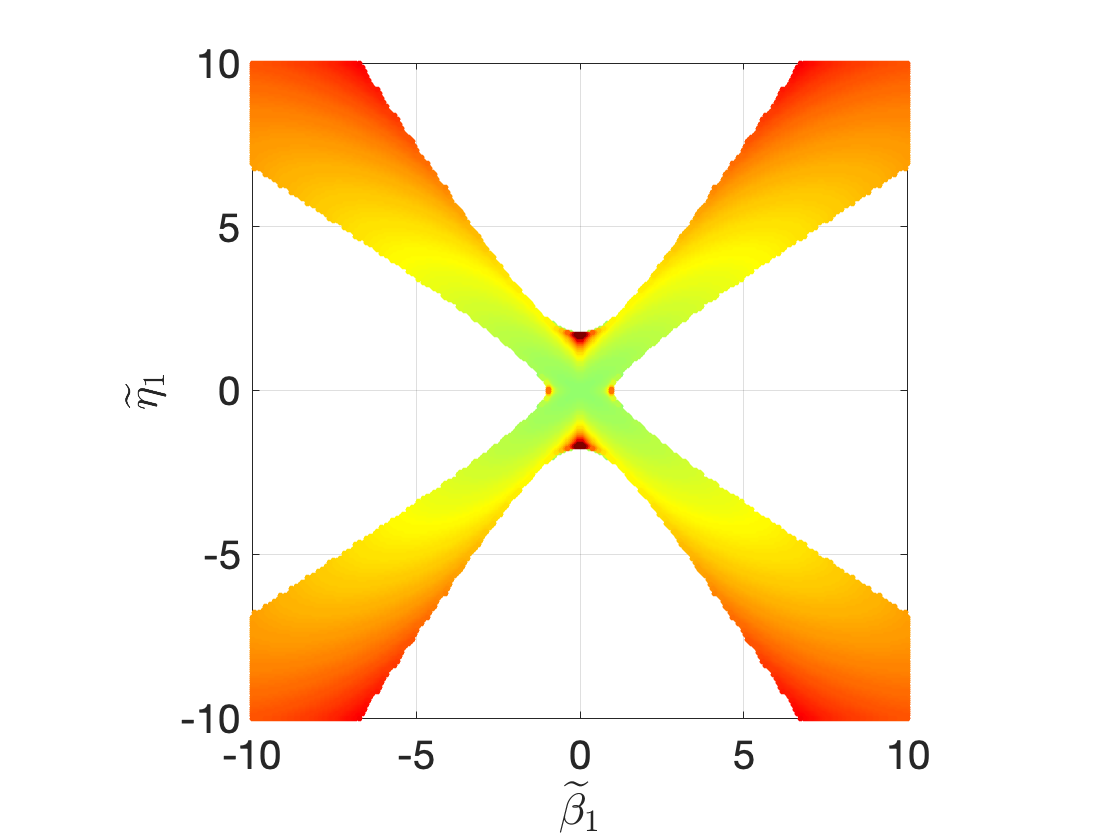}
\hspace{2cm}
\includegraphics[trim=4.2cm 0.2cm 6.0cm 1.5cm, clip, width=6cm]{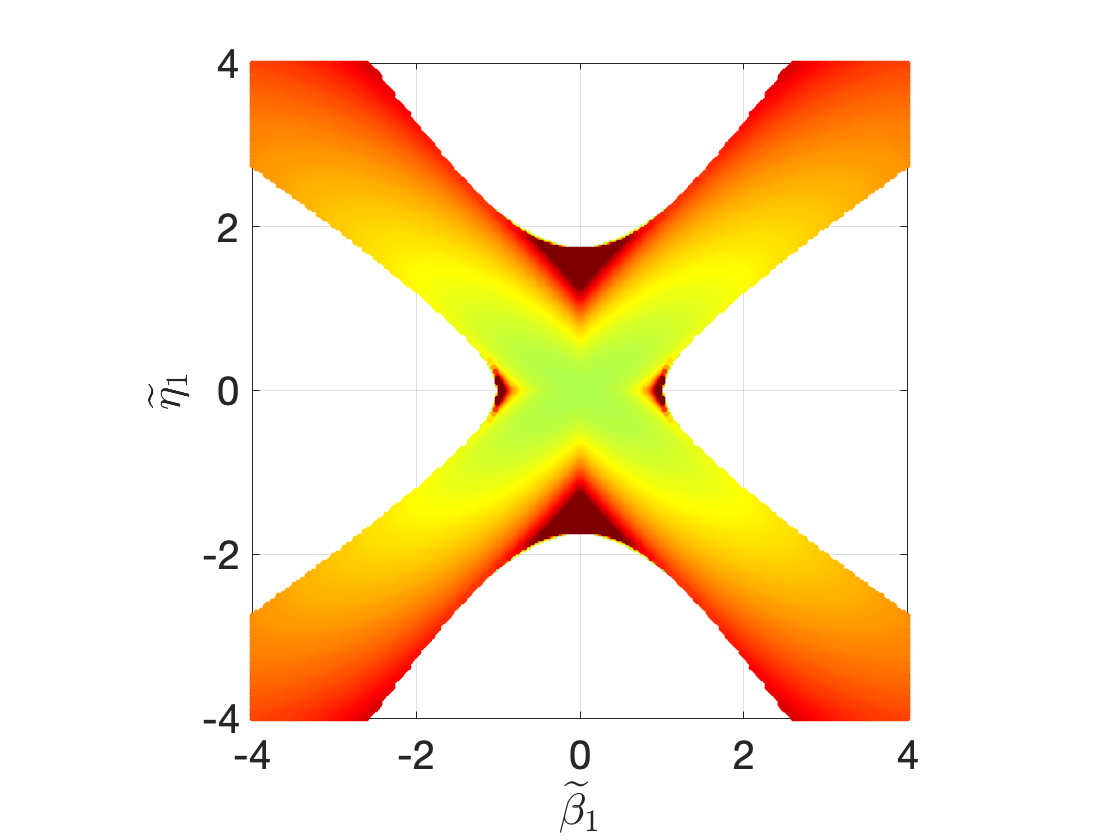}}
\caption{\sf Hyperbolic region of the system \eqref{eq:firstorder} for values $b_m\in[-5,5]$ and $\widetilde{\beta},\ \widetilde{\eta}\in[-10,10]$ (left column) and for the zoomed region $b_m\in[-2,2]$ and $\widetilde{\beta},\ \widetilde{\eta}\in[-4,4]$ (right column), where $\widetilde{\beta},\ \widetilde{\eta}$ are defined in \eqref{eq:q_scalings}. Note the slice $b_m = 0$, which is always hyperbolic, is omitted for visibility of the rest of the region.
Also note that the colors are purely for the purpose of providing depth.}
\label{fig:hyperbolic1}
\end{figure}

The region in which the system \eqref{eq:firstorder} is hyperbolic, i.e., where the quartic in \eqref{eq:quartic} has purely real roots, is presented in Figure \ref{fig:hyperbolic1} for $gh = 1$. 
Note that for the sake of visibility of the rest of the region, the value $b_m = 0$ is omitted, as the system \eqref{eq:firstorder} is globally hyperbolic when $b_m = 0$.
In addition, note that the colors in Figure \ref{fig:hyperbolic1} are present purely to provide depth to the 3-D region.
In comparison to the shallow water moment equations of \cite{Kowalski2019Moment}, which is globally hyperbolic for $M = 1$, the MRSWME loses \textit{global} hyberbolicity for any $M > 0$.
This is due to the introduction of a magnetic field, rather than a magnetic field vertical profile, which can be seen in Figure \ref{fig:hyperbolic1} along the slice $\widetilde{\eta} = 0$.

This loss of global hyperbolicity is potentially problematic, as it can result in unmitigated growth of linear instabilities.
The breakdown of hyperbolicity is a well known challenge in the contexts of moment models; see e.g., \cite{torrilhon2016modeling} where this arises in moment models for the Boltzmann equation.
However, all tests presented in the \S\ref{sec6} are either in the hyperbolic region, or barely outside of it (i.e., $\mathfrak{Im}(\lambda)/ \mathfrak{Re}(\lambda) \ll 1$). 
For those cases that do produce eigenvalues with a relatively small complex portion, we project these eigenvalues back into the real plane, and do not see any instabilities arise. 

For problems with significant breakdown in hyperbolicity (i.e., $\mathfrak{Im}(\lambda)/ \mathfrak{Re}(\lambda) $ larger than some threshold), there are several different techniques proposed to treat potential issues; see, e.g., \cite{fan2016model,levermore1996moment,KoellermeierRominger2020,Koellermeier2020steady}.
The numerical treatment of such cases where hyperbolicity is lost will be considered in future work on the MRSWME.

%
%

\section{Numerical Methods}\label{sec5}
Here we describe the methods to solve the 1-D version of MRSWME \eqref{eq:fin_h}--\eqref{eq:h_eta} in the $y$-direction, and the corresponding 2-D (with dimensions $y\zeta$) reference system.
The methods we describe are variants of that presented in \cite{Chertock2024Locally}, modified to fit the MRSWME \eqref{eq:fin_h}--\eqref{eq:h_eta} and reference system \eqref{eq:newh}--\eqref{eq:newhb}.

\subsection{Scheme for the MRSWME }\label{sec5.1}
Before introducing the numerical method, we write the 1-D MRSWME in the vector form
\begin{equation}\label{eq:gen_form_1D}
	\bm U_t + \bm G(\bm U)_y = Q(\bm U) \bm U_y + \bm S (\bm U),
\end{equation}
where the $4M+5$ conservative variables are $\bm U = \p{h, hu_m, hv_m, ha_m, hb_m, h\bm \Psi_1^\top,\dots, h\bm \Psi_M^\top}^\top$, $\bm\Psi_i = \p{\alpha_i, \beta_i, \gamma_i, \eta_i}^\top$ for $i = 1,\dots, M$,

\begin{equation}\label{eq:fluxes}
\begin{aligned}
	\bm G(\bm U) &= 
	\begin{pmatrix}
		hv_m\\
		h u_mv_m - ha_mb_m + h\sum_{\ell = 1}^M
		\frac{\alpha_\ell\beta_\ell-\gamma_\ell\eta_\ell}{2 \ell +1} \\
		h v_m^2 +\frac{1}{2}gh^2- hb_m^2 + h\sum_{\ell = 1}^M
		\frac{\beta_\ell^2-\eta_\ell^2}{2 \ell +1} \\
		ha_mv_m - h b_mu_m + h\sum_{\ell = 1}^M
		\frac{\beta_\ell\gamma_\ell-\alpha_\ell\eta_\ell}{2 \ell +1}\\
		0 \\
		\bm G ^{\bm \Psi}_1 (\bm U)\\
		\vdots\\
		\bm G ^{\bm \Psi}_M (\bm U)
	\end{pmatrix},\\[4pt]
	\bm G^{\bm \Psi}_i (\bm U) &= 
	\begin{pmatrix}
		h u_m \beta_i + hv_m \alpha_i-ha_m \eta_i-hb_m\gamma_i
		+h\sum_{\ell,n=1}^M A_{i\ell n}\p{\alpha_\ell \beta_n
		-\gamma_\ell \eta_n} \\
		2 \p{hv_m \beta_i-hb_m \eta_i}+h\sum_{\ell,n=1}^M
		A_{i\ell n}\p{\beta_\ell \beta_n-\eta_\ell \eta_n}\\
		h a_m \beta_i + hv_m \gamma_i -h b_m\alpha_i - h u_m\eta_i
		+h \sum_{\ell,n=1}^M A_{i\ell n}\p{\beta_\ell\gamma_n 
		- \alpha_\ell\eta_n}\\
		0
	\end{pmatrix},
\end{aligned}
\end{equation}
are the fluxes of the 1-D system, the source term reads
\begin{equation}\label{eq:sources}
	\bm S(\bm U) = 
	\begin{pmatrix}
		0\\
		fhv_m\\
		-fhu_m-ghZ_y \\
		0\\
		0 \\
		\bm S ^{\bm \Psi}_1 (\bm U)\\
		\vdots\\
		\bm S ^{\bm \Psi}_M (\bm U)
	\end{pmatrix},
	\qquad 
	\bm S ^{\bm \Psi}_i (\bm U)
	\begin{pmatrix}
	fh\beta_i\\
	-fh\alpha_i\\
	0\\
    0
	\end{pmatrix},
\end{equation}
and 
\begin{equation}\label{eq:Qmat}
\begin{aligned}
	&Q(\bm U) = \\
	&\begin{pmatrix}
		0 & 0 & 0 & 0 & 0 & 0\ \cdots\ 0 & 0\quad \cdots\quad\ 0 & 0\ \cdots\ 0 & 0\ \cdots\ 0 \\
		 & & & & -a_m - \sum_{\ell = 1}^M \gamma_\ell \phi_\ell(1) & & & & \\
		 & & & & -b_m - \sum_{\ell = 1}^M \eta_\ell \phi_\ell(1) & & & & \\
		\vdots &\vdots &\vdots &\vdots & -u_m - \sum_{\ell = 1}^M \alpha_\ell \phi_\ell(1) &\vdots \qquad \vdots & \vdots \qquad\qquad\ \vdots&\vdots \qquad \vdots & \vdots \qquad\qquad\ \vdots\\
		 & & & & -v_m - \sum_{\ell = 1}^M \beta_\ell \phi_\ell(1) & & 0\quad\ \cdots\quad\ 0& & 0\quad\ \cdots\quad\ 0\\
		 & & & & \bm Q^{(hb)}_1 & & \bm Q^{(h\beta_1)}_1\ \cdots\ \bm Q^{(h\beta_M)}_1 & &\bm Q^{(h\eta_1)}_1\ \cdots\ \bm Q^{(h\eta_M)}_1 \\
		 & & & & \vdots & &\vdots\quad\  \ddots\quad\ \  \vdots & & \vdots\quad\  \ddots\quad\ \  \vdots \\
		 0& 0& 0& 0& \bm Q^{(hb)}_M & 0\ \cdots \ 0 &  \bm Q^{(h\beta_1)}_M\ \cdots\ \bm Q^{(h\beta_M)}_M&0\ \cdots \ 0 &\bm Q^{(h\eta_1)}_M\ \cdots\ \bm Q^{(h\eta_M)}_M \\
	\end{pmatrix},
\end{aligned}
\end{equation}
is the nonconservative matrix, where 
\begin{equation} \label{eq:Q_Piece}
	\bm Q^{(hb)}_i = 
	\begin{pmatrix}
		-\sum_{\ell = 1}^M \gamma_\ell \p{1 + \Gamma_{i\ell}}\\
		-\sum_{\ell = 1}^M \eta_\ell \p{1 + \Gamma_{i\ell}}\\
		-\sum_{\ell = 1}^M \alpha_\ell \p{1 + \Gamma_{i\ell}}\\
		-\sum_{\ell = 1}^M \beta_\ell \p{1 + \Gamma_{i\ell}}
	\end{pmatrix},\ 
	\bm Q^{(h\beta_n)}_i = 
	\begin{pmatrix}
		u_m \delta_{in} - \sum_{\ell = 1}^M B_{i\ell n} \alpha_\ell \\
		v_m \delta_{in} - \sum_{\ell = 1}^M B_{i\ell n} \beta_\ell \\
		a_m \delta_{in} - \sum_{\ell = 1}^M B_{i\ell n} \gamma_\ell \\
		b_m \delta_{in} - \sum_{\ell = 1}^M B_{i\ell n} \eta_\ell
	\end{pmatrix},\ 
	\bm Q^{(h\beta_n)}_i = 
	\begin{pmatrix}
		-a_m \delta_{in} + \sum_{\ell = 1}^M B_{i\ell n} \gamma_\ell \\
		-b_m \delta_{in} + \sum_{\ell = 1}^M B_{i\ell n} \eta_\ell \\
		-u_m \delta_{in} + \sum_{\ell = 1}^M B_{i\ell n} \alpha_\ell \\
		-v_m \delta_{in} + \sum_{\ell = 1}^M B_{i\ell n} \beta_\ell
	\end{pmatrix},
\end{equation}
with Kronecker delta $\delta_{in}$, and $B_{i\ell n}$ defined in \eqref{eq:Coefs}. 
The system \eqref{eq:gen_form_1D}--\eqref{eq:Q_Piece} is paired with the 1-D divergence-free condition of the magnetic field $(hb_m)_y = 0,$ which, in combination with the 1-D evolution equation for $hb_m$, implies $hb_m$ is constant for all time.

\subsubsection{Path-conservative central upwind scheme}\label{sec5.1.1}

Following that of \cite{Chertock2024Locally}, we assume that the numerical solution, realized in terms of its cell averages 
\begin{equation*}\label{eq:cellavg}
	\xbar{\bm U}_j(t) \approx \frac{1}{\dy}\int\displaylimits_{C_j} \bm U(y,t) \ dy,
\end{equation*}
is available at time $t$. 
Here, we introduce a uniform Cartesian mesh with finite volume cells $C_j = [y_{\jmh}, y_{\jph}]$ with $y_{\jph} - y_{\jmh} = \dy$, $j = 1,\dots,N_y$.
Note we will omit the time dependence in the semi-discretization for the sake of brevity. 

To evolve the cell averages $\xbar{\bm U}_j$ in time, we extend the Path-Conservative Central-Upwind (PCCU) scheme introduced in \cite{Castro2019PCCU} to the system \eqref{eq:gen_form_1D}--\eqref{eq:Q_Piece}, resulting in the following semi-discretization for the MRSWME:
\begin{equation}\label{eq:PCCU}
\begin{aligned}
	\frac{d}{dt}\,\xbar{\bm U}_{j}=
	&-\frac{1}{\dy}\br{\bm{{\cal G}}_{\jph}-\bm{{\cal G}}_{\jmh}-\bm Q_{j}
	-\frac{s_{\jmh}^+}{s_{\jmh}^+-s_{\jmh}^-}\,\bm Q_{\bm\Upsilon,\jmh}
	+\frac{s_{\jph}^-}{s_{\jph}^+-s_{\jph}^-}\,\bm Q_{\bm\Upsilon,\jph}},
\end{aligned}
\end{equation}
where 
\begin{equation}\label{eq:CUflux}
	\bm{{\cal G}}_{\jph}=\frac{s_{\jph}^+\bm G\big(\bm U^{\rm N}_{j}\big)
	-s_{\jph}^-\bm G\big(\bm U^{\rm S}_{j+1}\big)}
	{s_{\jph}^+-s_{\jph}^-}+\frac{s_{\jph}^+s_{\jph}^-}{s_{\jph}^+
	-s_{\jph}^-}\left(\bm U^{\rm S}_{j+1}-\bm U^{\rm N}_{j}\right),
\end{equation}
are the CU numerical fluxes from \cite{Kurganov2002Sol}, $\bm U_j^{\rm N,S}$ denote the reconstructed right and left cell interface point values of cell $C_j$, respectively (see \S\ref{sec5.1.2}), $s_{\jph}^\pm$ denote the one-sided local speeds of propagation, estimated by using the maximum and minimum eigenvalues of the matrix $J(\bm U) = \pd{\bm G}{\bm U}(\bm U) - Q(\bm U)$:
\begin{equation}\label{eq:wavespeeds}
\begin{aligned}
	s_{\jph}^+ = \max\set{\lambda_i\p{J(\bm U_j^{\rm N})},\lambda_i\p{J(\bm U_{j+1}^{\rm S})},0},\\
	s_{\jph}^- = \min\set{\lambda_i\p{J(\bm U_j^{\rm N})},\lambda_i\p{J(\bm U_{j+1}^{\rm S})},0}, 
\end{aligned}
\end{equation}
for $i = 1,\dots, 5+4M$, where $M$ is the number of moments, and $\bm Q_j,\ \bm Q_{\bm \Upsilon,\jph}$ denote the discretization of the nonconservative products appearing on the right hand side of \eqref{eq:gen_form_1D} (see \S\ref{sec5.1.3}).

Since \eqref{eq:PCCU} is a system of ODEs, we must numerically integrate in time via an appropriate ODE solver. 
In the reported numerical results presented in \S\ref{sec6}, we opt for the explicit three-stage third-order strong stability preserving (SSP) Runge-Kutta (RK) method; see, for example, 
\cite{Gottlieb2001SSP,
Gottlieb2011SSP};
in which we meet the CFL condition
\begin{equation}\label{eq:CFL}
	\Delta t \leq \nu \frac{\dy}{\max\limits_j\set{s_{\jph}^+,-s_{\jph}^-}},
\end{equation}
with CFL number $\nu \leq 0.5$, to ensure stability.

\subsubsection{Reconstruction}\label{sec5.1.2}
To compute the cell interface point values $\bm U_j^{\rm N,S}$, we use a conservative piecewise linear reconstruction:
\begin{equation}\label{eq:gen_rec}
	\widetilde{\bm U}(y) = \xbar{\bm U}_j + \p{\bm U_y}_j(y-y_j),
	\qquad y \in C_j,
\end{equation}
resulting in the following right and left cell interface values, respectively:
\begin{equation}\label{eq:1D_rec}
	\bm U_j^{\rm N} = \xbar{\bm U}_j + \frac{\dy}{2}\p{\bm U_y}_j,\qquad 
	\bm U_j^{\rm S} = \xbar{\bm U}_j - \frac{\dy}{2}\p{\bm U_y}_j.
\end{equation}
To avoid oscillations, the slopes $\p{\bm U_y}_j$ are computed using a nonlinear limiter; we opt for the generalized minmod limiter (see, e.g., 
\cite{Lie2003Art,
Nessyahu1990Nonosc,
Sweby1984High}):
\begin{equation}\label{eq:limiter}
    \p{\bm U_y}_j = {\textrm{minmod}}\p{
    \theta \frac{\xbar{\bm U}_{j+1}-\xbar{\bm U}_{j}}{\dy},\ 
    \frac{\xbar{\bm U}_{j+1}-\xbar{\bm U}_{j-1}}{2\dy},\ 
    \theta \frac{\xbar{\bm U}_{j}-\xbar{\bm U}_{j-1}}{\dy} },
\end{equation}
for $\theta \in [1,2]$, and 
\begin{equation}\label{eq:minmod}
	\mbox{minmod}(z_1,z_2,\ldots)=\left
	\{\begin{aligned}
		&\min(z_1,z_2,\cdots)&&\mbox{if}~z_i>0,~\forall i,\\
		&\max(z_1,z_2,\cdots)&&\mbox{if}~z_i<0,~\forall i,\\
		&0&&\mbox{otherwise},
	\end{aligned}\right.
\end{equation}
is applied component-wise. 
Note that, in general, additional treatment must be used to implement the divergence-free condition of the magnetic field. 
However, in the 1-D case, the divergence-free condition implies $hb_m$ is constant for all time.
Therefore, if $hb_m$ is constant initially, then the reconstruction \eqref{eq:1D_rec} with limiter \eqref{eq:limiter} will preserve the divergence-free constraint for all time.

\subsubsection{Discretization of nonconservative products}\label{sec5.1.3}
To evaluate the nonconservative terms $Q(\bm U)\bm U_y$ on the right hand side of \eqref{eq:gen_form_1D}, we follow
\cite{Castro2019PCCU,Chu2023Fifth}
and compute the terms $\bm Q_j$ and $\bm Q_{\bm \Upsilon, \jph}$ in \eqref{eq:PCCU} by evaluating the following integrals exactly:
\begin{equation}\label{eq:noncons}
	\bm Q_{j}=\int\limits_{y_\jmh}^{y_\jph}
	Q\p{\widetilde{\bm U}(y)} \p{\widetilde{\bm U}(y)}_y\, dy, \qquad
	\bm Q_{\bm\Upsilon,\jph}=\int\limits_0^1
	Q\p{\bm U\p{\bm\Upsilon_{\jph}(s)}}\bm\Upsilon_{\jph}'(s)\, ds,
\end{equation}
where $Q(\bm U)$ is defined in \eqref{eq:Qmat}--\eqref{eq:Q_Piece}, $\widetilde{\bm U}(y)$ is the global (in space) interpolant defined in \eqref{eq:gen_rec}, and $\bm\Upsilon_{\jph}(s)$ is a linear path connecting the states $\bm U^{\rm N}_{j}$ and $\bm U^{\rm S}_{j+1}$:
\begin{equation*}\label{eq:Upsilon}
	\bm\Upsilon_{\jph}(s)=\bm U^{\rm N}_{j}+s[\bm U]_{\jph}, \qquad 
	[\bm U]_{\jph} = \bm U_{j+1}^{\rm S} - \bm U_{j}^{\rm N}.
\end{equation*}
Notice, though, that due to the structure of all the non-zero terms in $Q(\bm U)$, the integrals in \eqref{eq:noncons} can be written in the following arbitrary form:
\begin{equation}\label{eq:noncons_gen}
\begin{aligned}
	Q_{j}^{(i)}=\sum_k A_k\int\limits_{y_\jmh}^{y_\jph}
	\frac{\widetilde{\chi}_k(y)}{\widetilde{h}(y)}\p{\widetilde{\psi}_k(y)}_y\, dy,\qquad
	Q_{\bm\Upsilon,\jph}^{(i)}=\sum_k A_k\int\limits_0^1
	\frac{\chi_k\p{\bm\Upsilon_{\jph}(s)}}
 	{h\p{\bm\Upsilon_{\jph}(s)}}
	[\psi_k]_\jph\, ds,
\end{aligned}
\end{equation}
where $A_k$ are the appropriate coefficients appearing in \eqref{eq:Qmat}--\eqref{eq:Q_Piece}, and $\chi_k(y),\psi_k(y) \in \bm U(y)$ denote the appropriate components of $\bm U(y)$ that appear in $Q(\bm U)$ in \eqref{eq:Qmat}--\eqref{eq:Q_Piece}; for example, in the first non-zero component $Q^{(2,5)}(\bm U)$, $A_1 = - 1$, $A_{\ell + 1 } = -\phi_\ell(1)$ for $\ell = 1, \dots, M$, $\chi_1(y) = ha_m$, $\chi_{\ell+1}(y) = h\gamma_\ell$ for $\ell = 1, \dots, M$, and $\psi_k(y) = hb_m$ for all $k$.

Therefore, to evaluate the integrals in \eqref{eq:noncons} exactly, we only need to evaluate the remaining integrals in \eqref{eq:noncons_gen}. 
After removing the subscript $k$'s for the sake of brevity and substituting in the global interpolant from \eqref{eq:gen_rec}, we obtain the generic integral evaluations needed to compute the nonconservative product terms in \eqref{eq:noncons}:
\begin{equation}\label{eq:intQ}
\begin{aligned}
	\int\limits_{y_\jmh}^{y_\jph}
	\frac{\widetilde{\chi}(y)}{\widetilde{h}(y)}\widetilde{\psi}(y)_y\, dy  
	&=\left\{\begin{aligned}
	&\frac{\xbar{\chi}_j}{\xbar{h}_j}\p{\psi_y}_j\dy &&\mbox{if}~(h_y)_{j}=0,\\
	&\frac{\p{\psi_y}_j}{\p{h_y}_j}\br{ \p{\xbar{\chi}_j - \frac{\p{\chi_y}_j}{\p{h_y}_j}\xbar{h}_j }
	\ln\abs{\frac{h_j^{\rm N}}{h_j^{\rm S}}} + \p{\chi_y}_j \dy }
	&&\mbox{otherwise},
	\end{aligned}\right.\\[4pt]
	\int\limits_0^1
	\frac{\chi\p{\bm\Upsilon_{\jph}(s)}}
 	{h\p{\bm\Upsilon_{\jph}(s)}}
	[\psi]_\jph\, ds 
	& = \left\{\begin{aligned}
	& \frac{1}{2}\frac{1}{h_j^{\rm N}}\p{\chi_{j+1}^{\rm S} + \chi_{j}^{\rm N}} [\psi]_{\jph}
	&& \mbox{if}~[h]_{\jph}=0,\\
	& \frac{[\psi]_\jph}{[h]_\jph}\br{\frac{\p{\chi_j^{\rm N} [h]_\jph 
	- h_j^{\rm N} [\chi]_\jph }}{[h]_\jph}\ln\abs{\frac{h_{j+1}^{\rm S}}{h_j^{\rm N}}} 
	+ [\chi]_\jph} && \mbox{otherwise}.
	\end{aligned}\right.
\end{aligned}
\end{equation}

\subsection{Scheme for MRSW reference system}\label{sec5.2}
The 2-D MRSW reference system stems from that in \eqref{eq:newh}--\eqref{eq:newhb}, but is now reduced to only have $y$- and $\zeta$-spatial dependence.
This reduction is applied so that it represents the expected solution of the 1-D MRSWME (which is spatially dependent on $y$ but reduces the $\zeta$-dependence with a polynomial expansion).

For computational purposes, we additionally make two adjustments to the $y\zeta$-reduced form of \eqref{eq:newh}--\eqref{eq:newhb}: (i) we replace $v_m$ with $\omega$ in the height balance \eqref{eq:newh} using \eqref{eq:omega_rel}; and (ii) we include the Godunov-Powell source term for purposes of numerically treating the divergence-free condition of the magnetic field.
Directly using the reference system divergence condition in \eqref{eq:newdivb}, the Godunov-Powell source would be associated to $(hb_m)_y$, but this does not account for the the $\zeta$-dependence.
Alternatively, we also replace $(hb_m)_y$ using condition in \eqref{eq:mcC_rel}; that is, we substitute in $(hb_m)_y = (h\tilde b)_y + (h\mc C)_\zeta = 0$ for the divergence-free condition. 
Applying these modifications to the reference system \eqref{eq:newh}--\eqref{eq:newhb}, omitting tildes for the scaled variables, and rewriting in vector form, the 2-D MRSW reference system used for numerical experiments reads
\begin{equation}\label{eq:vecform2D}
	\bm U_t + \bm G(\bm U)_y + \bm H(\bm U)_\zeta 
	= \bm Q(\bm U)\br{ (hb)_y + (h\mc C)_\zeta} + \bm S (\bm U),
\end{equation}
where $\bm U = (h, hu, hv, ha, hb)^\top$,
\begin{equation}\label{eq:flux2D}
	\bm G(\bm U) = 
	\begin{pmatrix}
		hv \\
		huv-hab \\
		\textstyle hu^2 + \frac{g}{2}h^2 - hb^2 \\
		hav -hbu \\
		0
	\end{pmatrix}, \qquad 
	\bm H(\bm U) = 
	\begin{pmatrix}
		h\omega \\
		hu\omega-ha\mc C \\
		hv \omega - hb \mc C \\
		ha\omega -hu\mc C \\
		hb\omega - hv\mc C
	\end{pmatrix},
\end{equation}
are the fluxes in the $y$- and $\zeta$-directions, respectively, $\bm Q (\bm U) = -(a,b,u,v)^\top$ is the only non-zero vector of the full nonconservative matrix $Q(\bm U)$, and $\bm S(\bm U) = (0, fhv, -fhu-ghZ_y, 0, 0)^\top$ is the source term.
The system \eqref{eq:vecform2D}--\eqref{eq:flux2D} is paired with the updated $y\zeta$ divergence-free condition $(hb)_y + (h\mc C)_\zeta = 0$, and is closed by the definitions of $\omega$ and $\mc C$ in \eqref{eq:omega_new}--\eqref{eq:mcC_new}.

\subsubsection{Path-conservative central upwind scheme}\label{sec5.2.1}
Again following the scheme in \cite{Chertock2024Locally}, we assume at time $t$, the numerical solution in terms of its cell averages 
\begin{equation*}\label{eq:cellavg2D}
	\xbar{\bm U}_{j,k}(t) \approx \frac{1}{\dy\dz}\iint\limits_{C_{j,k}} \bm U(y,t) \ dy,
\end{equation*}
is available. 
Here, we introduce a uniform mesh with finite volume cells $C_{j,k} = [y_\jmh, y_\jph]\times[\zeta_\kmh,\zeta_\kph]$, where $y_{\jph} - y_{\jmh} = \dy$ for $j = 1,\dots,N_y$ and $\zeta_{\kph} - \zeta_{\kmh} = \dz$, $j = 1,\dots,N_\zeta$.
Note we will again omit the time dependence in the semi-discretization. 

To advance the cell averages $\xbar{\bm U}_{j,k}$ in time, we again apply the PCCU scheme from \cite{Castro2019PCCU} to the system \eqref{eq:vecform2D}--\eqref{eq:flux2D}.
Therefore, the semi-discretization for the 2-D MRSW reference system is 
\begin{equation}\label{eq:PCCU2D}
\begin{aligned}
	\frac{{d}}{{d}t}\,\xbar{\bm U}_{j,k}=
	&-\frac{1}{\dy}\br{\bm{{\cal G}}_{\jph,k}-\bm{{\cal G}}_{\jmh,k}-\bm Q_{j,k}^{(y)}
	-\frac{s_{\jmh,k}^+}{s_{\jmh,k}^+-s_{\jmh,k}^-}\,\bm Q_{\Upsilon,\jmh,k}^{(y)}
	+\frac{s_{\jph,k}^-}{s_{\jph,k}^+-s_{\jph,k}^-}\,\bm Q_{\Upsilon,\jph,k}^{(y)}}\\
	&-\frac{1}{\dz}\br{\bm{{\cal H}}_{j,\kph}-\bm{{\cal H}}_{j,\kmh}-\bm Q_{j,k}^{(\zeta)}
	-\frac{s_{j,\kmh}^+}{s_{j,\kmh}^+-s_{j,\kmh}^-}\,\bm Q_{\Upsilon,j,\kmh}^{(\zeta)}
	+\frac{s_{j,\kph}^-}{s_{j,\kph}^+-s_{j,\kph}^-}\,\bm Q_{\Upsilon,j,\kph}^{(\zeta)}},
\end{aligned}
\end{equation}
where $\bm{\mc G}_{\jph,k}$, defined in \eqref{eq:CUflux}, and 
\begin{equation}\label{eq:CUfluxz}
	\bm{{\cal H}}_{j,\kph}
	=\frac{s_{j,\kph}^+\bm H(\bm U^{\rm U}_{j,k})
	-s_{j,\kph}^-\bm H(\bm U^{\rm D}_{j,k+1})}
	{s_{j,\kph}^+-s_{j,\kph}^-}+\frac{s_{j,\kph}^+s_{j,\kph}^-}{s_{j,\kph}^+
	-s_{j,\kph}^-}\left(\bm U^{\rm U}_{j,k+1}-\bm U^{\rm D}_{j,k}\right)
\end{equation}
are the CU numerical fluxes introduced in \cite{Kurganov2002Sol}, $\bm U_{j,k}^{\rm N,S,U,D}$ are the right, left, top, and bottom cell interface point values of $C_{j,k}$, respectively (see \S\ref{sec5.2.2}), $s^\pm_{\jph,k}$ and $s^\pm_{j,\kph}$ are the one-sided local speeds of propagation in the $y$- and $\zeta$-directions, respectively, approximated by the maximum and minimum eigenvalues of the matrices $\pd{\bm G}{\bm U}(\bm U)-Q(\bm U)$ and $\pd{\bm H}{\bm U}(\bm U)$:
\begin{equation*}\label{eq:wavespeeds2D}
\begin{aligned}
	s_{\jph,k}^+&=\max\set{v_{j,k}^{\rm N}+\sqrt{(b_{j,k}^{\rm N})^2+ gh_{j,k}^{\rm N}}\ ,
	v_{j+1,k}^{\rm S}+\sqrt{(b_{j+1,k}^{\rm S})^2 + gh_{j+1,k}^{\rm S}},\ 0},\\
	s_{\jph,k}^-&=\min\set{v_{j,k}^{\rm N}-\sqrt{(b_{j,k}^{\rm N})^2+ gh_{j,k}^{\rm N}}\ ,
	v_{j+1,k}^{\rm S}-\sqrt{(b_{j+1,k}^{\rm S})^2 + gh_{j+1,k}^{\rm S}},\ 0},\\
	s_{j,\kph}^+&= \max\set{\omega_{j,k}^{\rm U}+|\mc C_{j,k}^{\rm U}|,\ 
	\omega_{j,k+1}^{\rm D}+|\mc C_{j,k+1}^{\rm D}|,\ 0}, \\
	s_{j,\kph}^-&= \min\set{\omega_{j,k}^{\rm U}-|\mc C_{j,k}^{\rm U}|,\ 
	\omega_{j,k+1}^{\rm D}-|\mc C_{j,k+1}^{\rm D}|,\ 0}. \\
\end{aligned}
\end{equation*}
In addition, $\bm Q_{j,k}^{(y)}$, $\bm Q_{j,k}^{(\zeta)}$, $\bm Q_{\Upsilon,\jph,k}^{(y)}$, and $\bm Q_{\Upsilon,j,\kph}^{(\zeta)}$ denote the discretizations of the nonconservative products on the right hand side of \eqref{eq:vecform2D} (see \S\ref{sec5.2.3}),
and $\omega_{j,k},\ \omega_{j,k}^{\rm U,D}$ and $\mc C_{j,k},\ \mc C_{j,k}^{\rm U,D}$ are the cell center, upper cell interface, and lower cell interface approximations to the vertical coupling terms $\omega$ and $\mc C$ from \eqref{eq:omega_new}, \eqref{eq:mcC_new}, respectively (see \S\ref{sec5.2.4}).

To advance the system of ODEs in \eqref{eq:PCCU2D} in time, we numerically integrate and apply an appropriate ODE solver. 
Similarly to that of the MRSWME, the numerical results presented in \S\ref{sec6} use the explicit third-order SSP RK method; see, e.g., 
\cite{Gottlieb2001SSP,
Gottlieb2011SSP};
in which we ensure stability by meeting the CFL condition
\begin{equation}\label{eq:CFL2D}
	\Delta t \leq \min{\p{
	\frac{\nu \dy}{\max\limits_{j,k}\set{s_{\jph,k}^+,-s_{\jph,k}^-}},
	\frac{\nu \dz}{\max\limits_{j,k}\set{s_{j,\kph}^+,-s_{j,\kph}^-}}}},
\end{equation}
for CFL number $\nu \leq 0.5$.

\subsubsection{Reconstruction}\label{sec5.2.2}
The reconstruction for most variables in \eqref{eq:vecform2D}--\eqref{eq:flux2D} is a 2-D extension of that presented in \S\ref{sec5.1.2}. 
To compute the cell interface point values we use the conservative piecewise linear reconstruction 
\begin{equation}\label{eq:gen_rec2D}
	\widetilde{\bm U}(y,\zeta) = \xbar{\bm U}_{j,k} + \p{\bm U_y}_{j,k}(y-y_j)
	+ \p{\bm U_\zeta}_{j,k}(\zeta-\zeta_k),
	\qquad y,\zeta \in C_{j,k}.
\end{equation}
This results in the following cell interface values:
\begin{equation*}\label{eq:rec2D}
\begin{aligned}
	&\bm U^{\rm N}_{j,k}=\xbar{\bm U}_{j,k}+\frac{\dy}{2}(\bm U_y)_{j,k},\quad
	\bm U^{\rm S}_{j,k}=\xbar{\bm U}_{j,k}-\frac{\dy}{2}(\bm U_y)_{j,k},\\
	&\bm U^{\rm U}_{j,k}=\xbar{\bm U}_{j,k}+\frac{\dz}{2}(\bm U_\zeta)_{j,k},\quad
	\bm U^{\rm D}_{j,k}=\xbar{\bm U}_{j,k}-\frac{\dz}{2}(\bm U_\zeta)_{j,k},
\end{aligned}
\end{equation*}
where the slopes $(\bm U_y)_{j,k}$ are computed using the generalized minmod limiter in \eqref{eq:limiter}, and the slopes $(\bm U_\zeta)_{j,k}$ are computed analogously in the $\zeta$-direction:
\begin{equation}\label{eq:limiter2D}
    \p{\bm U_\zeta}_{j,k} = {\textrm{minmod}}\p{
    \theta \frac{\xbar{\bm U}_{j,k+1}-\xbar{\bm U}_{j,k}}{\dz},\ 
    \frac{\xbar{\bm U}_{j,k+1}-\xbar{\bm U}_{j,k-1}}{2\dz},\ 
    \theta \frac{\xbar{\bm U}_{j,k}-\xbar{\bm U}_{j,k-1}}{\dz} },
\end{equation}
where the minmod function is defined in \eqref{eq:minmod}.

Unlike the 1-D MRSWME, the 2-D reference system now requires additional treatment to implement the divergence-free condition of the magnetic field $(hb)_y+(h\mc C)_\zeta = 0$. 
To ensure a locally-divergence-free scheme, we use the approach of \cite{Chertock2024Locally}; 
that is, we (i) write an evolution equation for $B = (hb)_y$ by differentiating the $hb$ equation in \eqref{eq:vecform2D}--\eqref{eq:flux2D}, which reads
\begin{equation*}\label{eq:Bevol}
	B_t + \p{vB-h\mc C v_\zeta}_y + \p{\omega B + hb\omega_y}_\zeta = 0;
\end{equation*}
(ii) evolve $B$ simultaneously with \eqref{eq:vecform2D}--\eqref{eq:flux2D} using the same CU fluxes presented in \eqref{eq:CUflux} and \eqref{eq:CUfluxz}, where we compute $v_\zeta$ and $\omega_y$ using the generalized minmod limiter slopes $(v_\zeta)_{j,k}$ from \eqref{eq:limiter2D} and $(\omega_y)_{j,k}$ from \eqref{eq:limiter}; 
and (iii) use the evolved values $B_{j,k}$ for the $y$-direction slopes within the reconstruction of $hb$. 
However, we also cannot reconstruct $hb$ using $B_{j,k}$ without including some limiter to avoid oscillations.
Therefore, to reconstruct $hb$ in the $y$-direction, we use 
\begin{equation}\label{eq:recB}
	(hb)^{\rm N}_{j,k}=(\xbar{hb})_{j,k}+\frac{\dy}{2}B_{j,k}\sigma_{j,k},\qquad
	(hb)^{\rm S}_{j,k}=(\xbar{hb})_{j,k}-\frac{\dy}{2}B_{j,k}\sigma_{j,k},
\end{equation}
with
\begin{equation*}\label{eq:msig}
	\sigma_{j,k}:=\left\{\begin{aligned}
		&\min\bigg\{1,\frac{((hb)_y)_{j,k}}{\xbar B_{j,k}}\bigg\}
		&&\mbox{if}~((hb)_y)_{j,k}\,\xbar B_{j,k}>0,\\
		&0&&\mbox{otherwise},
	\end{aligned}\right.
\end{equation*}
where $((hb)_y)_{j,k}$ is computed using the generalized minmod limiter in \eqref{eq:limiter}.

Since we only have the integral equation in \eqref{eq:mcC_new} for $h\mc C$, as opposed to an evolution equation, we cannot exactly follow discretization in \cite{Chertock2024Locally} for this term. 
We share the details on how $(h\mc C)_{j,k}^{\rm U,D}$ are computed in \S\ref{sec5.2.4}.

{\rmk\label{rmk5.1}
Notice that, to ensure a second-order method for $hb$, the evolution of $B$ only needs to be first-order accurate. This is because $B$ is utilized in the reconstruction \eqref{eq:recB} and always multiplied by $\dy$. 
}

\subsubsection{Discretization of nonconservative products}\label{sec5.2.3}
The discretization of the nonconservative terms $\bm Q(\bm U)\br{ (hb)_y + (h\mc C)_\zeta}$ in \eqref{eq:vecform2D} follow very closely to that presented in \S\ref{sec5.1.3}. 
Following \cite{Castro2019PCCU,Chu2023Fifth}, we compute the terms $\bm Q_{j,k}^{(y)}$, $\bm Q_{j,k}^{(\zeta)}$, $\bm Q_{\Upsilon,\jph,k}^{(y)}$, and $\bm Q_{\Upsilon,j,\kph}^{(\zeta)}$ in \eqref{eq:PCCU2D} by evaluating the following integrals exactly:
\begin{equation}\label{eq:noncons2D}
\begin{aligned}
	\bm Q_{j,k}^{(y)}=\int\limits_{y_\jmh}^{y_\jph}
	\bm Q\p{\widetilde{\bm U}(y,\zeta_k)} (\widetilde{hb})(y,\zeta_k)_y\, dy, \qquad
	\bm Q_{\Upsilon,\jph,k}^{(y)}=\int\limits_0^1
	\bm Q\p{\bm U\p{\Upsilon_{\jph,k}(s)}}\Upsilon_{\jph,k}'(s)\, ds,\\
	\bm Q_{j,k}^{(\zeta)}=\int\limits_{\zeta_\kmh}^{\zeta_\kph}
	\bm Q\p{\widetilde{\bm U}(y_j,\zeta)}
	(\widetilde{h\mc C})(y_j,\zeta)_\zeta\, d\zeta, \qquad
	\bm Q_{\Upsilon,j,\kph}^{(\zeta)}=\int\limits_0^1
	\bm Q\p{\bm U\p{\Upsilon_{j,\kph}(s)}}\Upsilon_{j,\kph}'(s)\, ds,
\end{aligned}
\end{equation}
where $\bm Q(\bm U) = -(a,b,u,v)^\top$, $\widetilde{\bm U}(y,\zeta)$ is the global interpolant defined in \eqref{eq:gen_rec2D}, and 
\begin{equation*}\label{eq:Upsilon2D}
\begin{aligned}
	\Upsilon_{\jph,k}(s)&=(hb)^{\rm N}_{j,k}+s[hb]_{\jph,k}, \qquad 
	&&[hb]_{\jph,k} = (hb)_{j+1,k}^{\rm S} - (hb)_{j,k}^{\rm N},\\
	\Upsilon_{j,\kph}(s)&=(h\mc C)^{\rm U}_{j,k}+s[h\mc C]_{j,\kph}, \qquad 
	&&[h\mc C]_{j,\kph} = (h\mc C)_{j,k+1}^{\rm D} - (h\mc C)_{j,k}^{\rm U},
\end{aligned}
\end{equation*}
are linear paths connecting states $(hb)_{j,k}^{\rm N}$ to $(hb)_{j+1,k}^{\rm S}$ and $(h\mc C)_{j,k}^{\rm U}$ to $(h\mc C)_{j,k+1}^{\rm D}$, respectively.

Notice that, due to the definition of $\bm Q(\bm U)$, all four integrals in \eqref{eq:noncons2D} also fall under the general form in \eqref{eq:intQ}.
We thus omit the final evaluations of \eqref{eq:noncons2D} for the sake of brevity.

\subsubsection{Discretization of vertical coupling terms}\label{sec5.2.4}
To discretize the vertical coupling term corresponding to velocity $\omega$ at each time step, we follow that in \cite{Kowalski2019Moment}; that is, we compute vertical cell interface values such that the $\omega$ discretization is continuous across cells. 
Starting from the definition in \eqref{eq:omega_new}, we obtain the discretization of $\omega_{j,k}^{\rm U} = \omega_{j,k+1}^{\rm D} \equiv \omega_{j,\kph}$ by 
(i) evaluating \eqref{eq:omega_new} at $\zeta_\kph$, noticing that the lower bound is $\zeta_{\frac{1}{2}} = 0$;
(ii) apply midpoint rule to \eqref{eq:omega_new}; and 
(iii) use the $y$-direction numerical fluxes \eqref{eq:CUflux} of the $h$ equation \eqref{eq:newh} to compute $((hv)_y)_{j,k}$ and the corresponding vertical average $((hv_m)_y)_{j}$. 
Doing so gives the following discretization of the velocity vertical coupling term:
\begin{equation}\label{eq:omega_rec}
	\omega_{j,k}^{\rm U} = \omega_{j,k+1}^{\rm D} = \frac{2}{h_{j,k+1}^{\rm D} +h_{j,k}^{\rm U}}
	\sum_{\ell = 1}^k \dz \p{\p{(hv_m)_y}_{j} 
	- \frac{1}{\dy}\br{\mc G_{\jph,\ell}^{(h)}-\mc G_{\jmh,\ell}^{(h)}}},
\end{equation}
where 
\begin{equation}\label{eq:hvm}
	\p{(hv_m)_y}_{j} = \sum_{k = 1}^{N_\zeta}\frac{\dz}{\dy}\br{\mc G_{\jph,k}^{(h)}-\mc G_{\jmh,k}^{(h)}}.
\end{equation}

The vertical coupling for the magnetic field must be computed differently since it contributes to the 2-D divergence-free condition $(hb)_y + (h\mc C)_\zeta = 0$. 
To ensure the preservation of the local divergence-free condition, the reconstruction of $h\mc C$ in the $\zeta$-direction must be consistent with that of $hb$ in the $y$-direction.
Therefore, to get the magnetic field vertical coupling approximation, we 
(i) take an average of both sides of \eqref{eq:mcC_new} through integration over cell $C_{j,k}$;
(ii) apply trapezoid rule to the averaging integral;
(iii) apply midpoint rule to the remaining integrals originating from \eqref{eq:mcC_new}; and
(iv) divide by $h_{j,k}$ and use the reconstruction slope from \eqref{eq:recB} for the computation of $[(hb)_y]_{j,k}$.
Doing so gives the following cell-centered point value of the magnetic field vertical coupling:
\begin{equation}\label{eq:mcC_cen}
	(h\mc C)_{j,k} = -\p{\frac{1}{2}\sigma_{j,k}B_{j,k} 
	+ \sum_{\ell = 1}^k \sigma_{j,\ell}B_{j,\ell}}.
\end{equation}
Then, to ensure the reconstruction of $(h\mc C)_{j,k}$ within the divergence-free condition is consistent, we simply use $-[(hb)_y]_{j,k}$ as the slope, in turn giving the cell interface point values of the magnetic field coupling:
\begin{equation}\label{eq:mcC_rec}
	(h\mc C)_{j,k}^{\rm U} = (h\mc C)_{j,k} + \frac{\dz}{2}\p{-\sigma_{j,k}B_{j,k} },\qquad
	(h\mc C)_{j,k}^{\rm D} = (h\mc C)_{j,k} - \frac{\dz}{2}\p{-\sigma_{j,k}B_{j,k} }.
\end{equation}
Note that the discretizations of $\omega$ in \eqref{eq:omega_rec}--\eqref{eq:hvm} and $\mc C$ in \eqref{eq:mcC_cen}--\eqref{eq:mcC_rec} are updated in each stage of the time update.

\section{Numerical Simulations}\label{sec6}

In this section, we demonstrate the viability of the 1-D MRSWME equations in comparison to the corresponding 2-D reference solution through four numerical tests. 
The implementation of these numerical experiments has been made publicly available \cite{Redle2025}. 
All simulations use a CFL number of $\nu = 0.45$ in \eqref{eq:CFL} and \eqref{eq:CFL2D}, gravitational constant $g = 1$, and, unless otherwise stated, a generalized minmod limiting parameter $\theta = 1.3$.

As discussed in \S\ref{sec4.4.2}, Examples 2-4 of the moment model experiments do indeed produce complex eigenvalues; however, the complex portions are several orders of magnitude smaller than those of their real parts. Therefore, in these simulations, we project the eigenvalues to the real plane to compute the wave speeds in \eqref{eq:wavespeeds}. 
We plan to address the difficulties regarding hyperbolicity regularization in future work.

\subsubsection*{Example 1---Small hump without magnetic field}\label{ex1}
While fairly simple to see within \eqref{eq:gen_form_1D}--\eqref{eq:Q_Piece}, we first confirm numerically that the MRSWME reduce to the SWME when a magnetic field is not present. 
Furthermore, this example also tests the reference solution and the viability of the proposed methods on an existing test case. 
To do so, we consider the smooth bump initial conditions suggested in \cite{Kowalski2019Moment}:
\begin{equation}\label{eq:ex1_ICs}
	h(y,0) = 1+\exp\p{3\cos\p{\pi (y + 0.5)}-4},\qquad
	v(y,0,\zeta) = \begin{cases}
		\frac{1}{4} & {\rm constant\ case}, \\
		\frac{1}{2}\zeta & {\rm linear\ case}, \\
		\frac{3}{2}\zeta(1-\zeta) & {\rm quadratic\ case}, \\
		3\zeta(1-\frac{5}{2}\zeta -\frac{5}{6}\zeta^2) & {\rm cubic\ case}, 
	\end{cases}
\end{equation}
with $u = a = b \equiv 0$, on the domain $y \in \br{-1,1}$ with periodic boundary conditions, and zero bottom topography or rotational forces.
Notice that the linear case is equivalent to taking $\beta_1 = -\frac{1}{4}$; likewise the quadratic case is equivalent to taking $\beta_1 = 0$ and $\beta_2 = -\frac{1}{4}$, and the cubic case is equivalent to using $\beta_1 = \beta_2 = 0$ and $\beta_3 = -\frac{1}{4}$.
For all four cases, $v_m = \frac{1}{4}$, implying the depth-averaged standard shallow water system results in the same solution for the four different vertical profiles. 
Note that for this example, we take the generalized minmod parameter to instead be $\theta = 1$ to help avoid some small oscillations that arise. 

\begin{figure}[ht!]
\centerline{
\includegraphics[trim=0.9cm 0.1cm 1.7cm 0.2cm, clip, width=8cm]{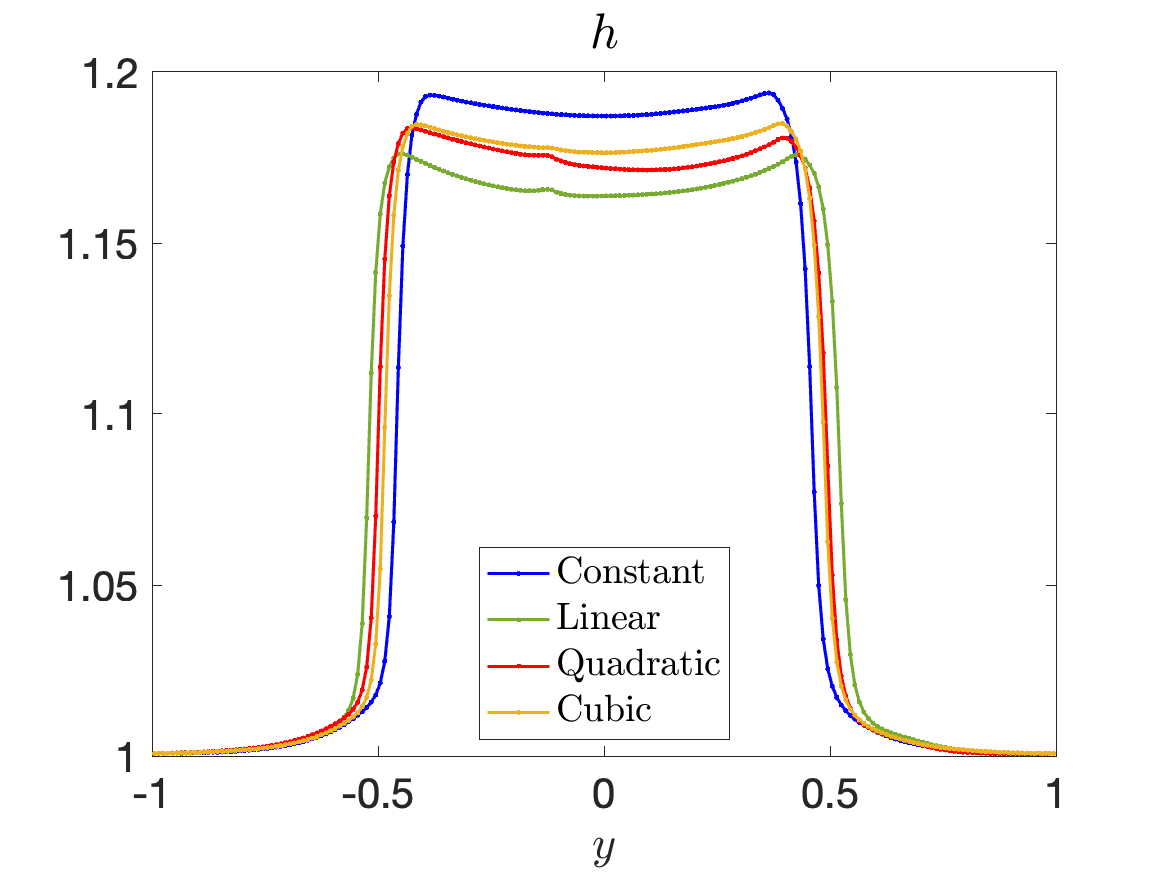}
\hspace{10pt}
\includegraphics[trim=0.9cm 0.1cm 1.7cm 0.2cm, clip, width=8cm]{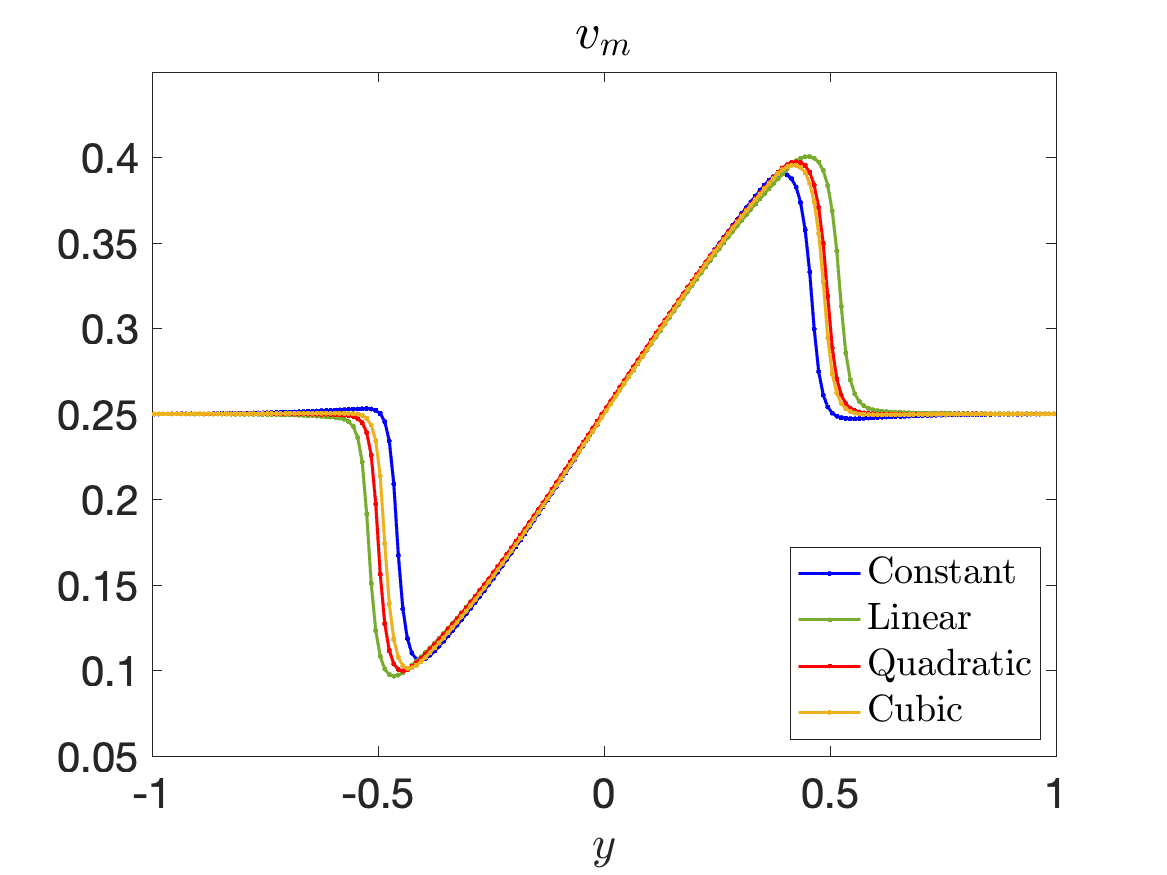}}
\caption{\sf Example 1: Reference system solutions $h$ (left) and $v_m$ (right) for the four different polynomial vertical profiles.}
\label{fig:bump_ref}
\end{figure}

We run the smooth bump problem with various vertical profiles to a final time $t = 2$ on a mesh of 200 cells in the $y$-direction, and for the reference solution additionally discretize in the $\zeta$-direction with 100 cells. 
First, we present the reference solutions of $h$ and the depth-averaged velocity $v_m$ for the four different vertical profiles in Figure \ref{fig:bump_ref}. 
The initial bump splits into a left and right wave that travel through each other due to the periodic boundary conditions. At $t = 2$, the waves traveled through one another twice, resulting in the appearance of two shock waves. 
Here we see strong agreement in our solutions with those presented in \cite{Kowalski2019Moment}; that is, we see that the inclusion of a vertical profile returns different wave speeds and amplitudes of $h$ and $v_m$ in comparison to that of the standard shallow water model. 

\begin{figure}[ht!]
\centerline{
\includegraphics[trim=0.9cm 0.1cm 1.7cm 0.2cm, clip, width=4.5cm]{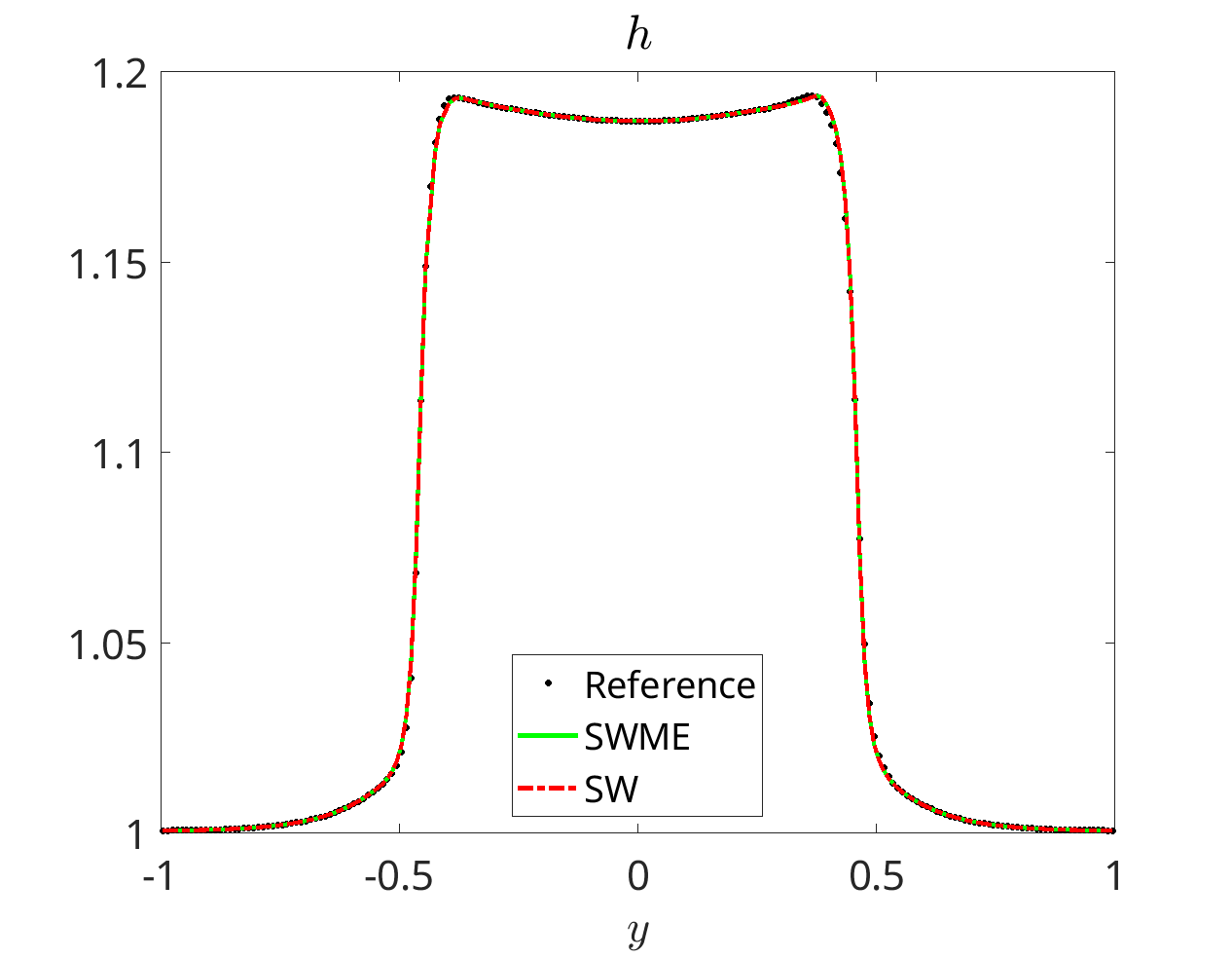}
\hspace{2pt}
\includegraphics[trim=0.9cm 0.1cm 1.7cm 0.2cm, clip, width=4.5cm]{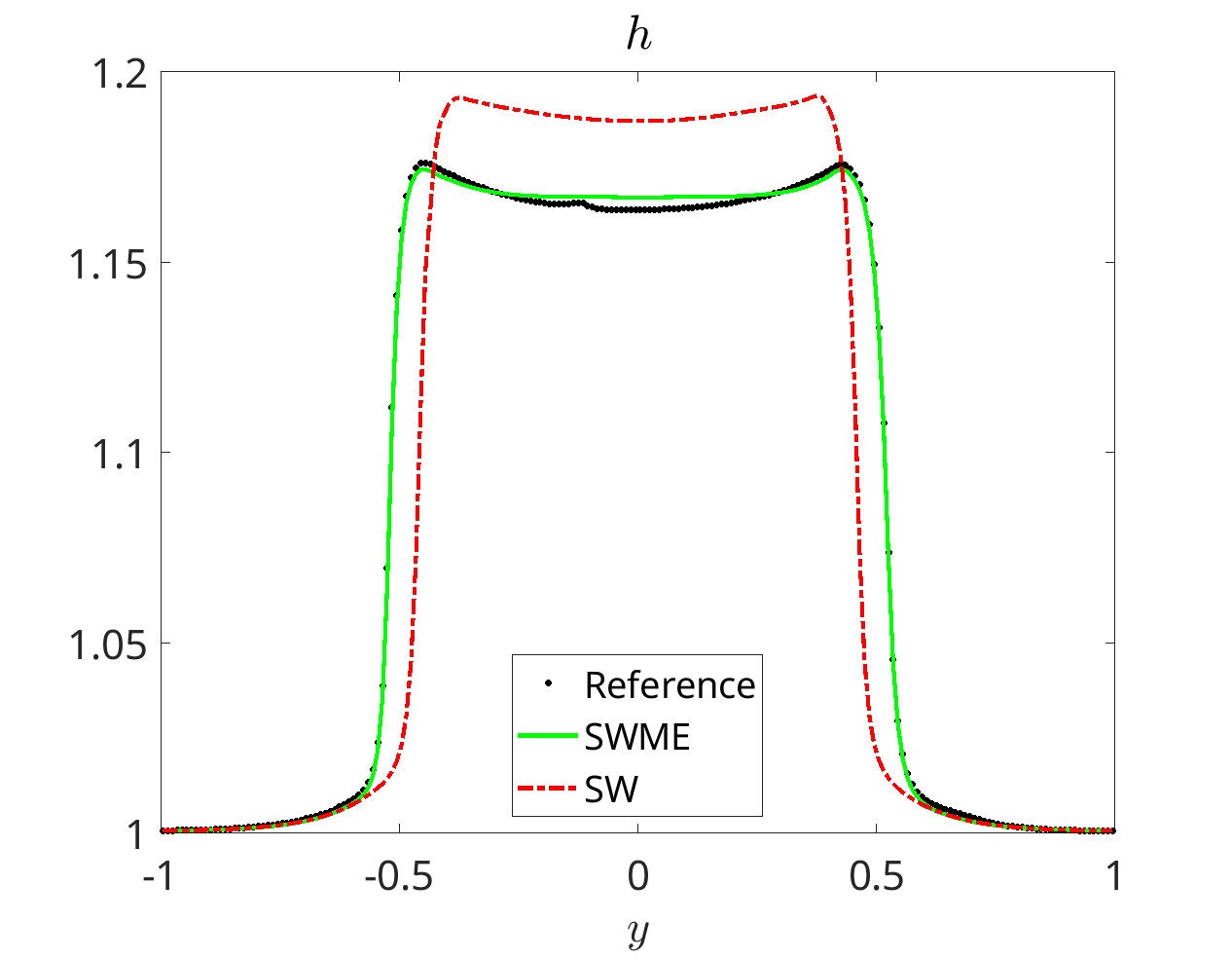}
\hspace{2pt}
\includegraphics[trim=0.9cm 0.1cm 1.7cm 0.2cm, clip, width=4.5cm]{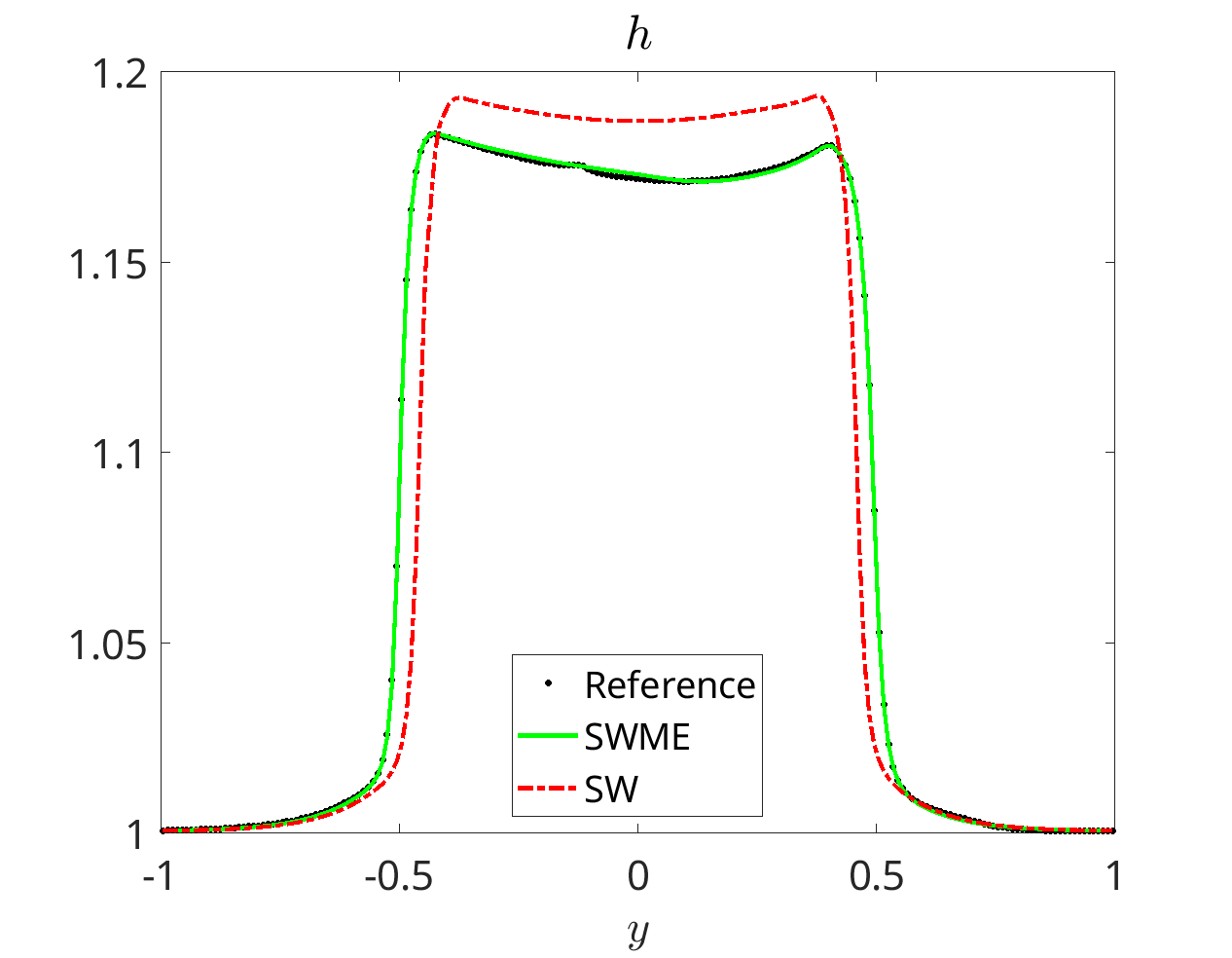}
\hspace{2pt}
\includegraphics[trim=0.9cm 0.1cm 1.7cm 0.2cm, clip, width=4.5cm]{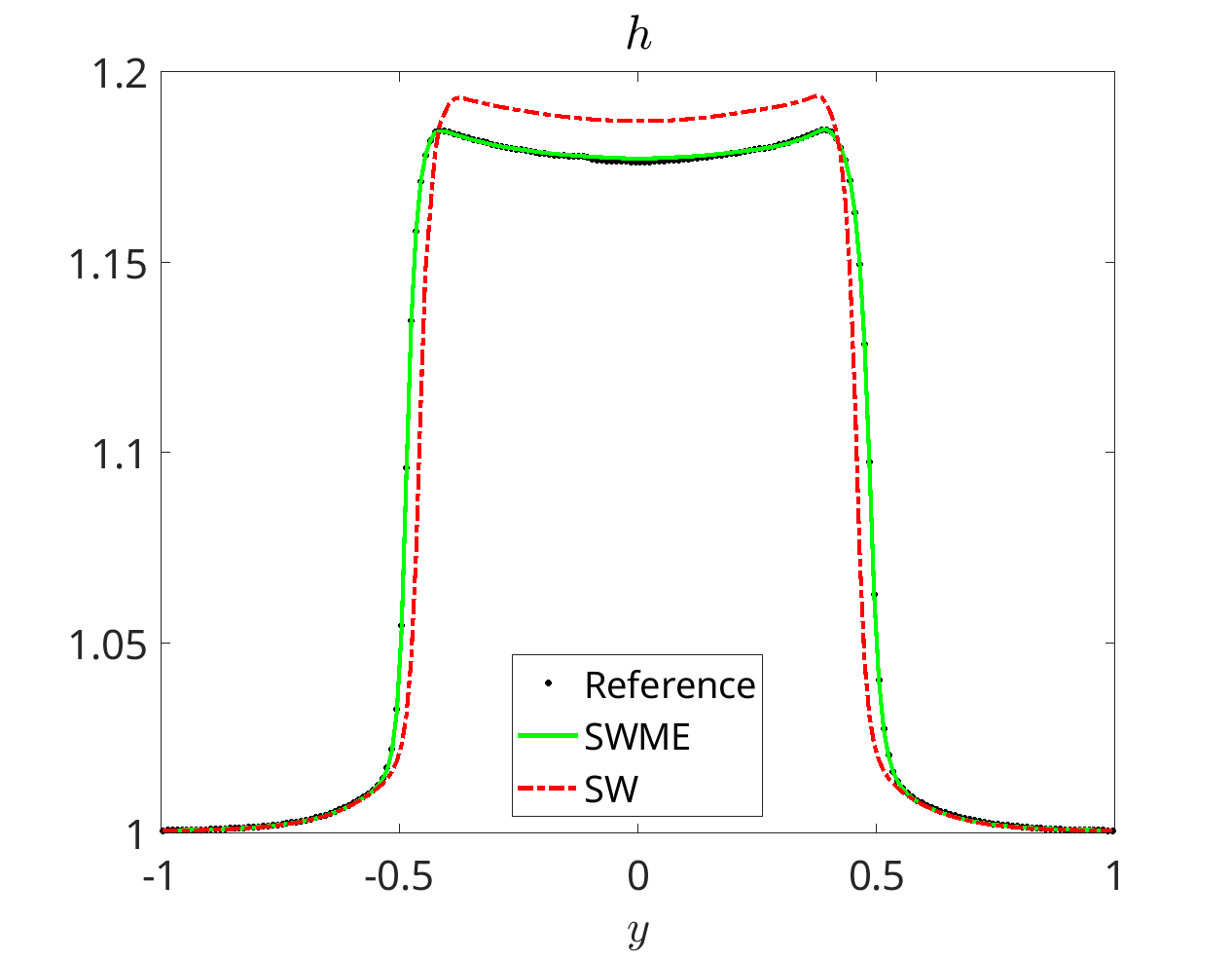}}
\vspace{7pt}
\centerline{
\includegraphics[trim=0.9cm 0.1cm 1.7cm 0.2cm, clip, width=4.5cm]{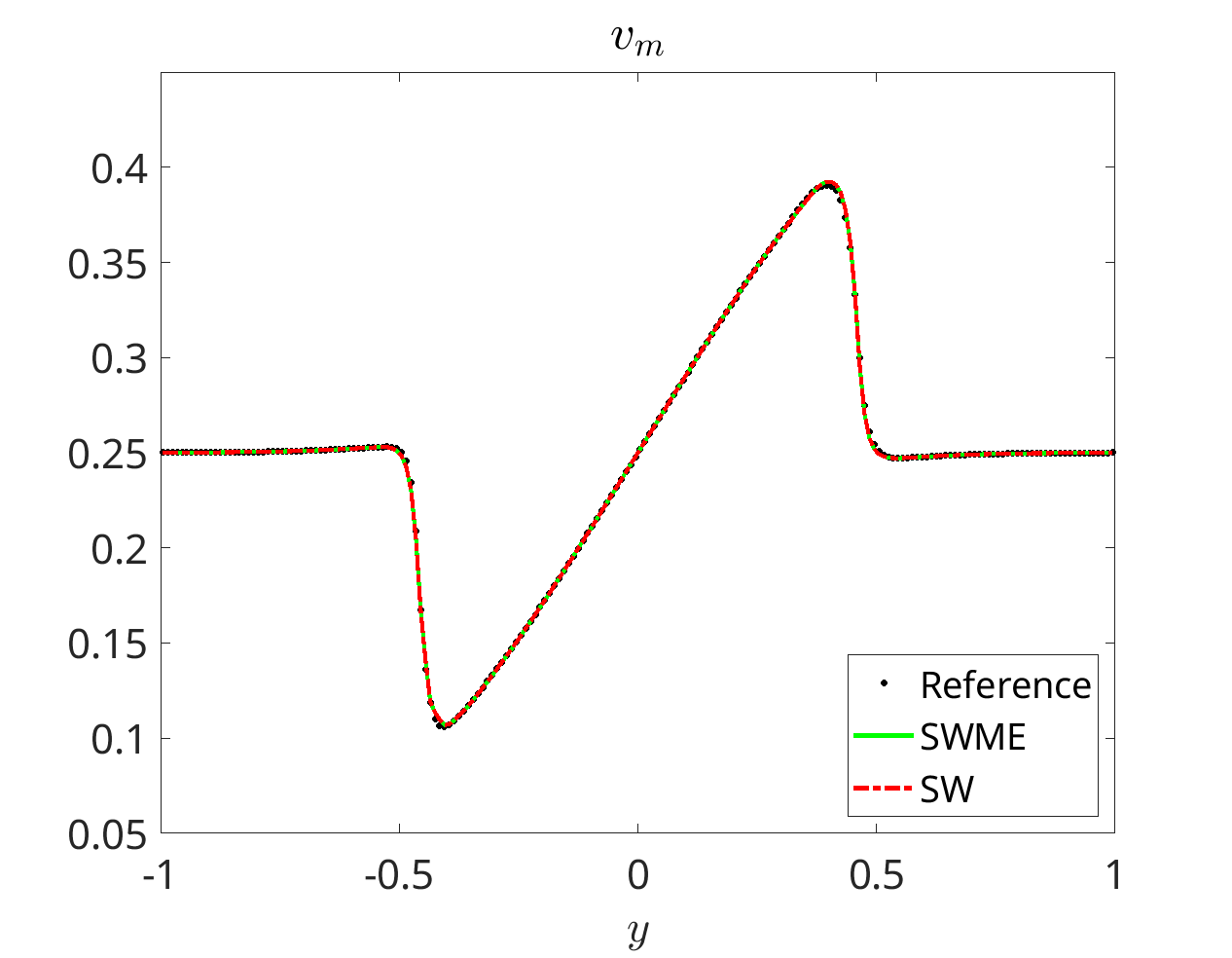}
\hspace{2pt}
\includegraphics[trim=0.9cm 0.1cm 1.7cm 0.2cm, clip, width=4.5cm]{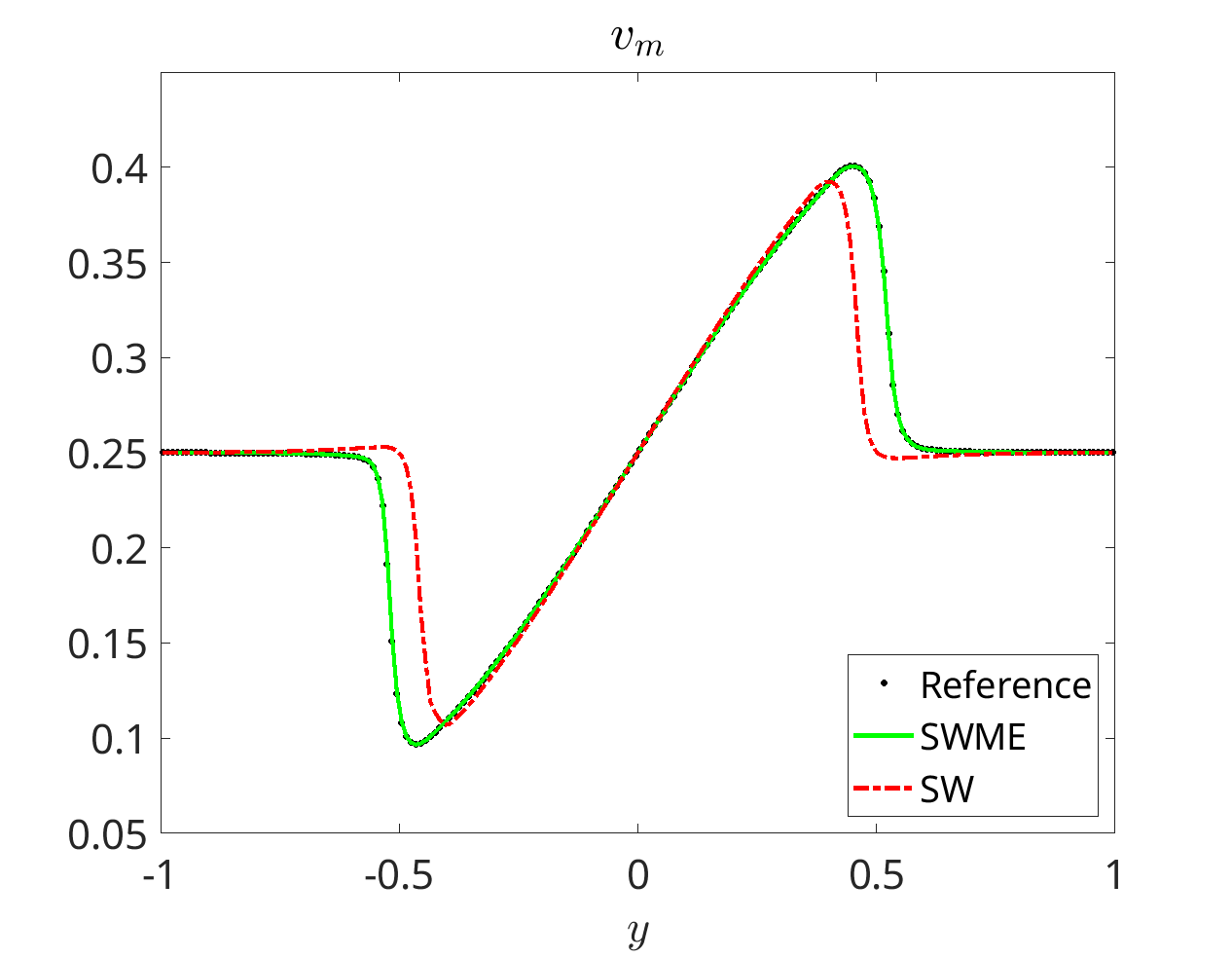}
\hspace{2pt}
\includegraphics[trim=0.9cm 0.1cm 1.7cm 0.2cm, clip, width=4.5cm]{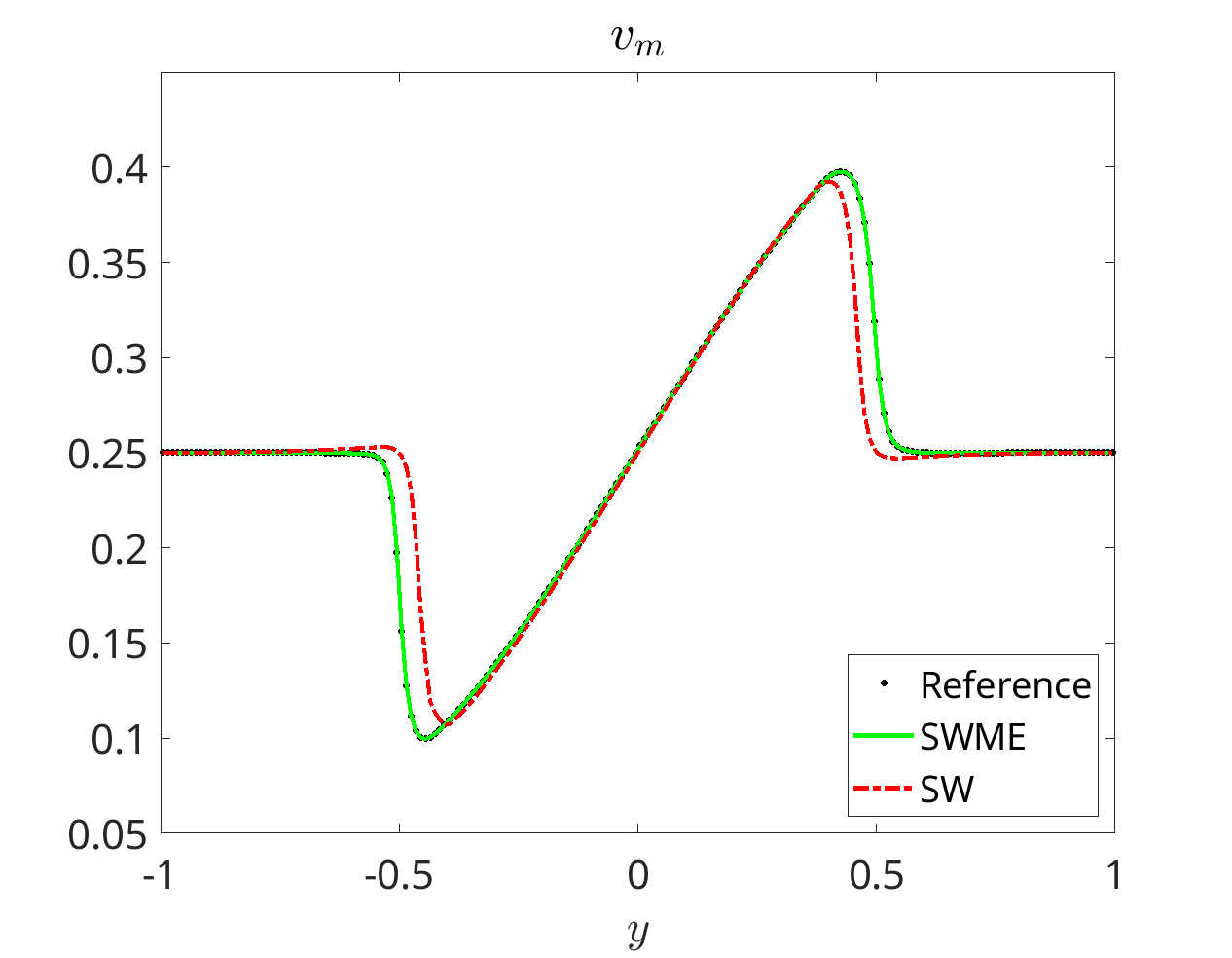}
\hspace{2pt}
\includegraphics[trim=0.9cm 0.1cm 1.7cm 0.2cm, clip, width=4.5cm]{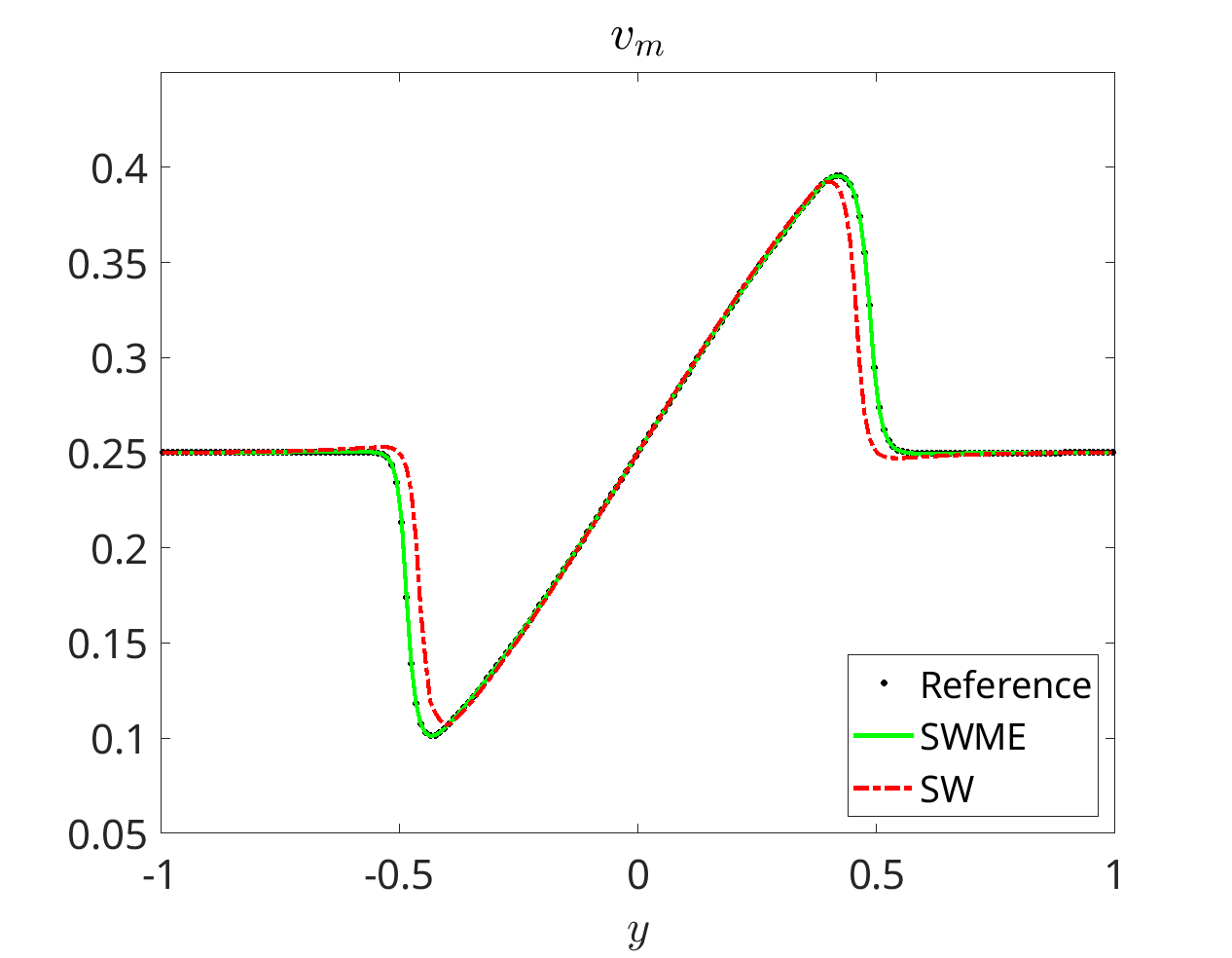}}
\caption{\sf Example 1: Solutions for $h$ (top row) and $v_m$ (bottom row) computed using the 2-D reference system, the 1-D standard shallow water (SW), and the SWME for the different vertical profile cases of $M = 0$ (first column), $M = 1$ (second column), $M = 2$ (third column), and $M = 3$ (last column).}
\label{fig:bump_comp}
\end{figure}

To then test the (MR)SWME in \eqref{eq:gen_form_1D}--\eqref{eq:Q_Piece}, we present a results comparison of the 2-D reference solution, the 1-D standard shallow water equations (since magnetic field is zero for all time), and the SWME in Figure \ref{fig:bump_comp}. 
In agreement with \cite{Kowalski2019Moment}, we also conclude that in the case of zero magnetic field, the SWME solutions better approximate the reference solution when a vertical profile is present, as the moment models better capture the appropriate wave speeds and amplitudes of the corresponding reference solutions.

\subsubsection*{Example 2---Small hump with magnetic field}\label{ex2}
In this simulation, we take the same problem as that proposed in Example 1, except now impose a comparatively strong magnetic field with a vertical profile of its own. 
The initial conditions are the same as that in Example 1 \eqref{eq:ex1_ICs}, but now additionally with the initial magnetic field being 
$$
	(hb)(y, 0, \zeta) = \begin{cases}
		1.1 & {\rm constant\ case}, \\
		1.1-\frac{1}{4}\phi_1(\zeta) & {\rm linear\ case}, \\
		1.1-\frac{1}{4}\phi_2(\zeta) & {\rm quadratic\ case}, \\
		1.1-\frac{1}{4}\phi_3(\zeta) & {\rm cubic\ case},
	\end{cases}
$$
where $\phi_\ell(\zeta),\ \ell = 1,\dots,3$, are the Legendre basis polynomials presented in \eqref{eq:basis}.
The above initial conditions are equivalent to initially taking $h\eta_1 = -\frac{1}{4}$ in the linear case, $h\eta_1 = 0$ with $h\eta_2 = -\frac{1}{4}$ in the quadratic case, and $h\eta_1 = h\eta_2 = 0$ with $h\eta_3 = -\frac{1}{4}$ in the cubic case.
All cases return a depth-averaged $y$-direction magnetic field of $hb_m = 1.1$.
Notice that these initial conditions are listed for $hb$ instead of $b$, as the initial conditions must satisfy the divergence-free constraint $(hb_m)_y = 0$. 

\begin{figure}[ht!]
\centerline{
\includegraphics[trim=0.9cm 0.1cm 1.7cm 0.2cm, clip, width=8cm]{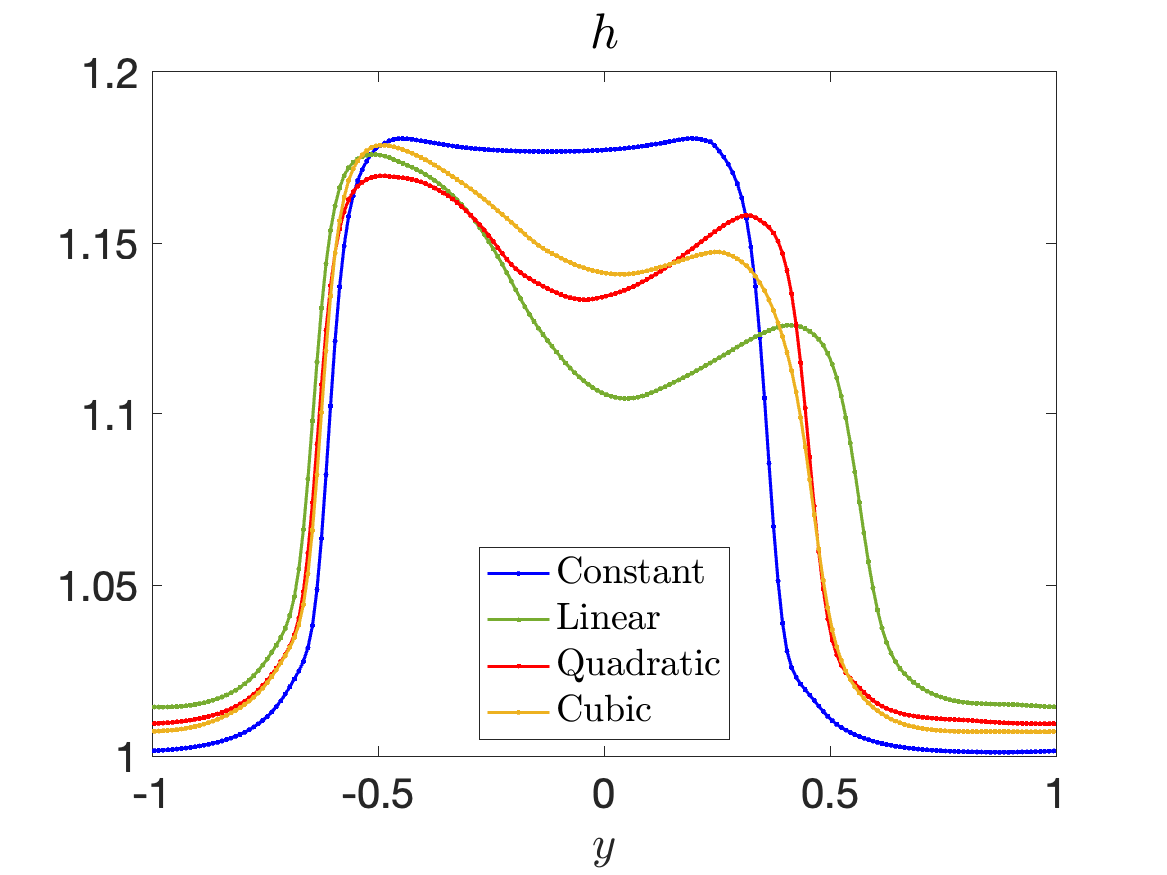}
\hspace{10pt}
\includegraphics[trim=0.9cm 0.1cm 1.7cm 0.2cm, clip, width=8cm]{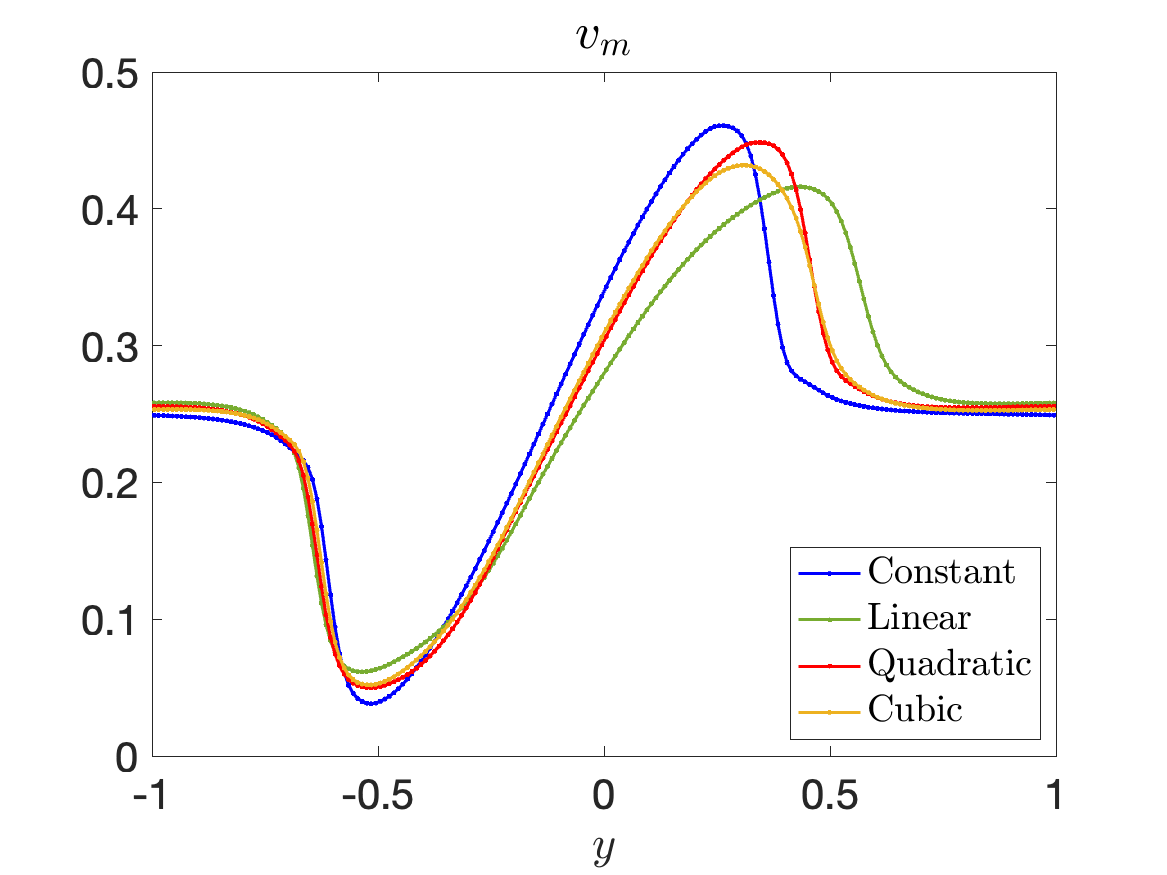}}
\caption{\sf Example 2: Reference system solutions $h$ (left) and $v_m$ (right) for the four different polynomial vertical profiles.}
\label{fig:bumpB_ref}
\end{figure}

To get a similar structure in solution to that of Example 1, we run this example to a final time of $t = 1.5$, and do so on the same spatial discretization.
The final time is selected earlier than that of Example 1 since the magnetic field causes the waves to travel faster; see, e.g., the wave speeds for the MRSW equations in \eqref{eq:MRSWspeeds}; and doing so will result in this same `plateau-type' structure in $h$. 
We present the reference solutions of $h$ and the depth-averaged velocity $v_m$ for the four different initial vertical profiles in Figure \ref{fig:bumpB_ref}, noting that we do not show the depth-averaged magnetic field $b_m$ since $hb_m$ is constant for all time (even if $hb$ is changing in the $\zeta$-direction).
Much like that of the non-magnetic field case in Example 1, we see the left and right traveling waves meet for the second time near the center of the periodic domain. 
However, unlike the previous example where the vertical profiles only shift the wave speeds and amplitudes slightly, the introduction of a magnetic field and its corresponding vertical profile significantly change the structure of the 2-D reference solution---especially that of the right-moving wave.

\begin{figure}[ht!]
\centerline{
\includegraphics[trim=0.9cm 0.1cm 1.7cm 0.2cm, clip, width=8cm]{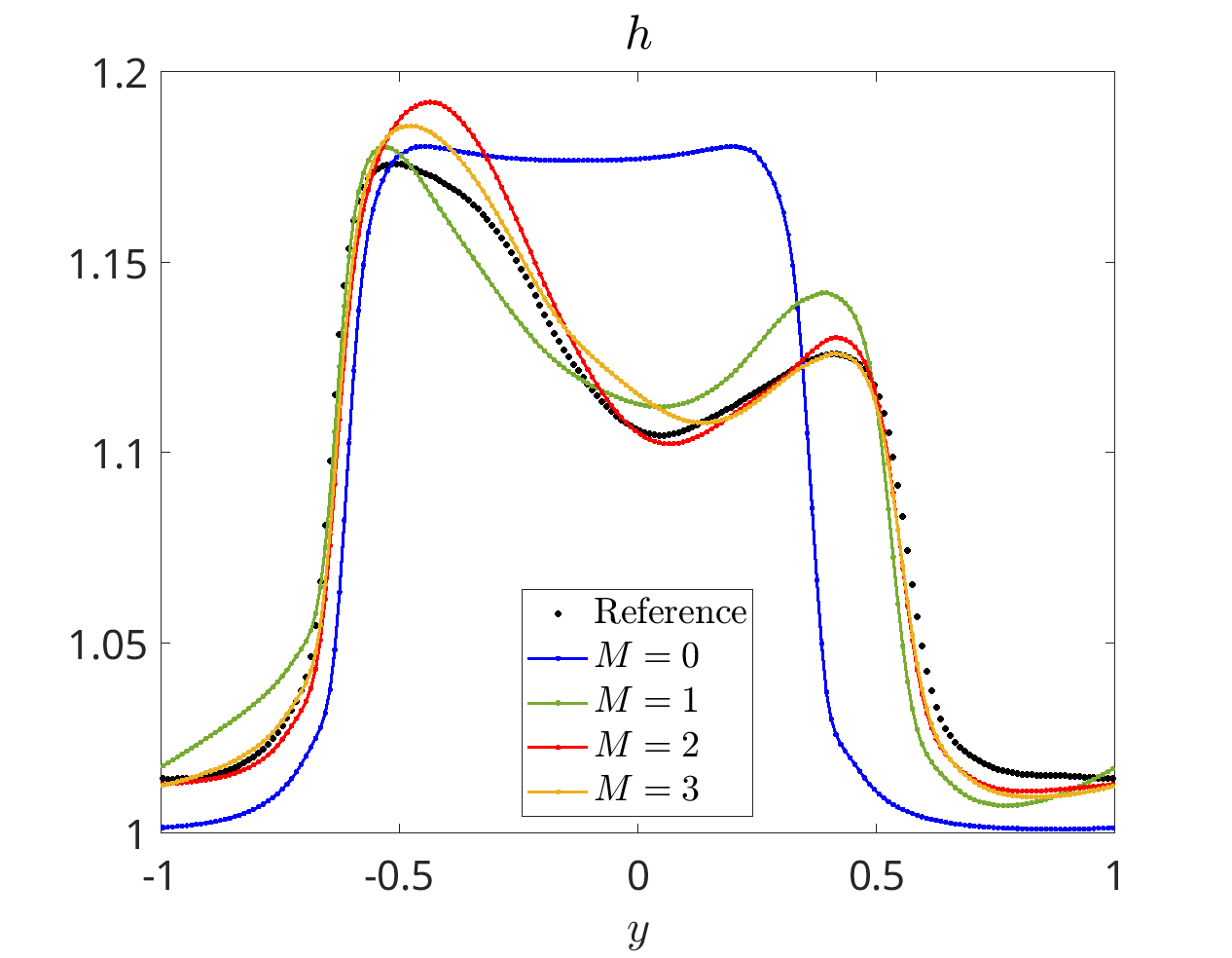}
\hspace{10pt}
\includegraphics[trim=0.9cm 0.1cm 1.7cm 0.2cm, clip, width=8cm]{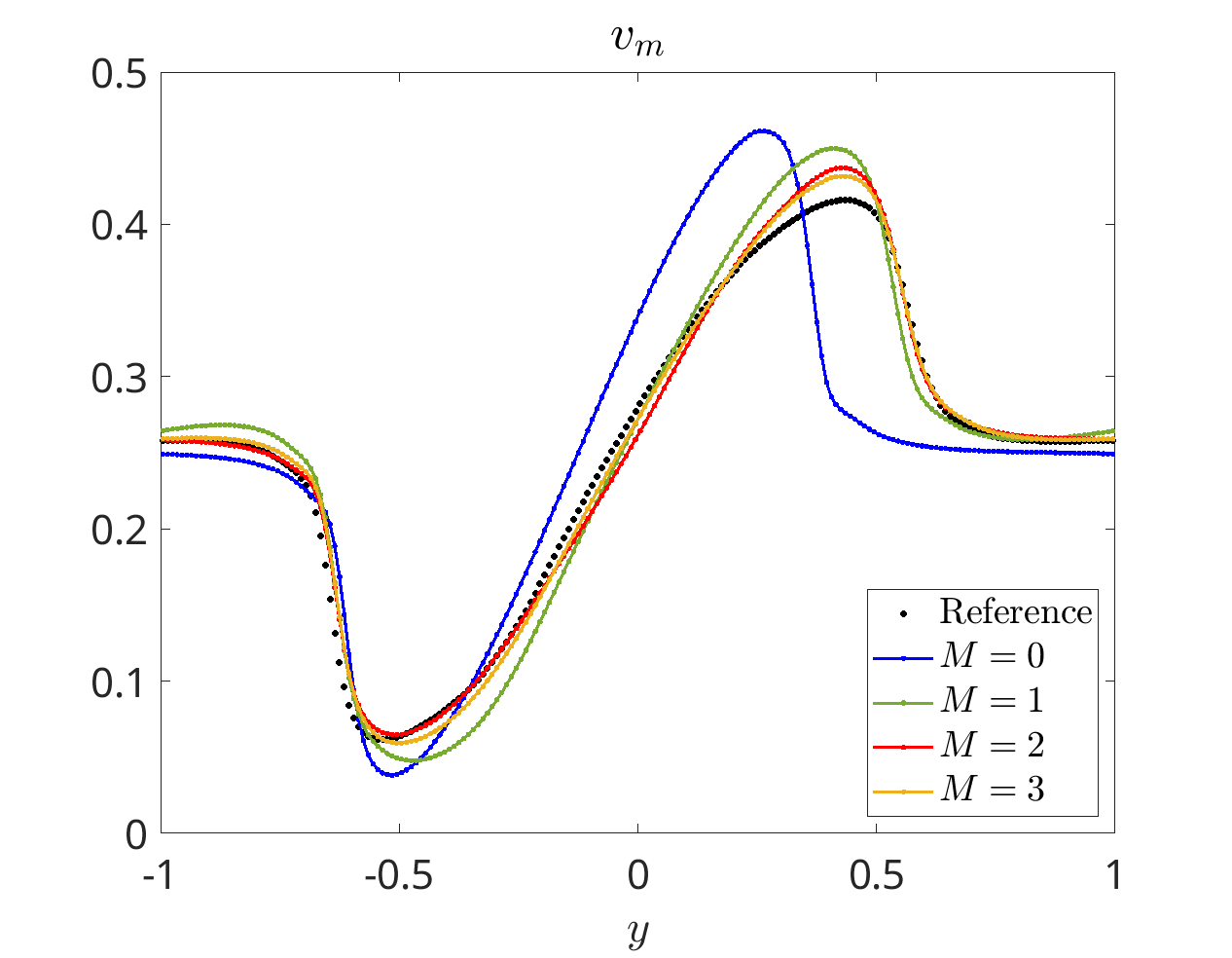}}
\caption{\sf Example 2 (Linear case): The solutions for $h$ (left) and $v_m$ (right) of the $Mth$-order MRSWME, $M = 0,\dots,3,$ against the reference solution.}
\label{fig:bumpB_comp}
\end{figure}

\begin{figure}[ht!]
\centerline{
\includegraphics[trim=0.6cm 0.1cm 1.5cm 0.2cm, clip, width=8cm]{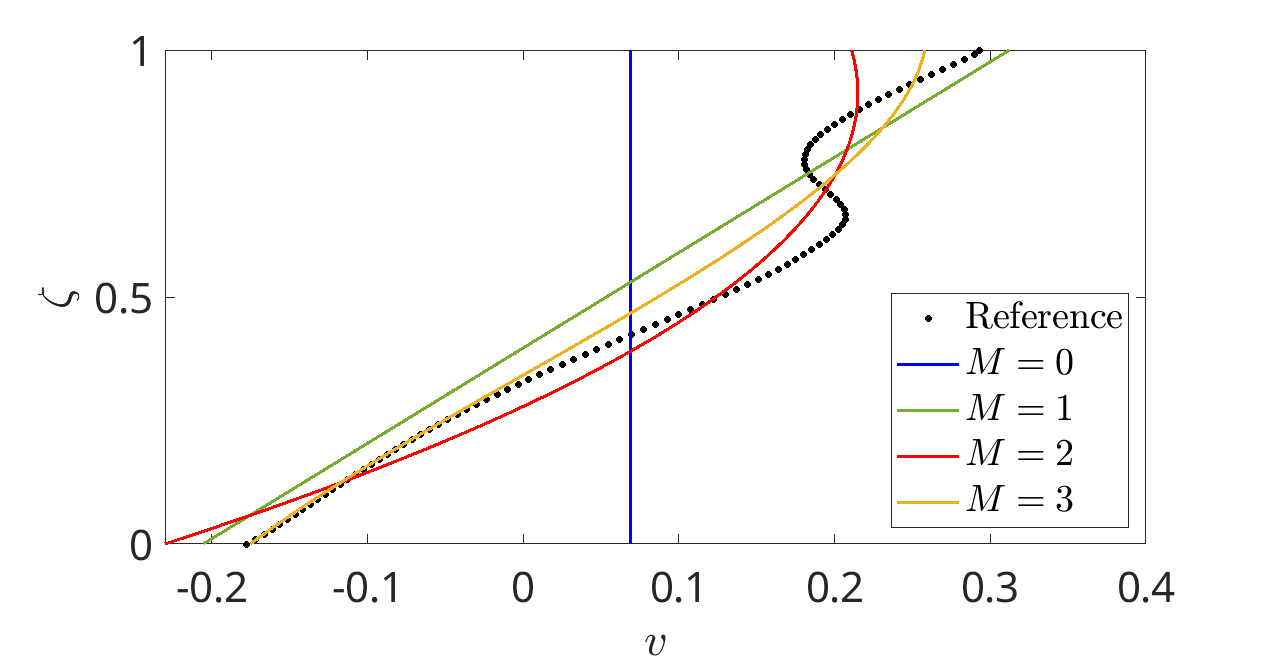}
\hspace{10pt}
\includegraphics[trim=0.6cm 0.1cm 1.5cm 0.2cm, clip, width=8cm]{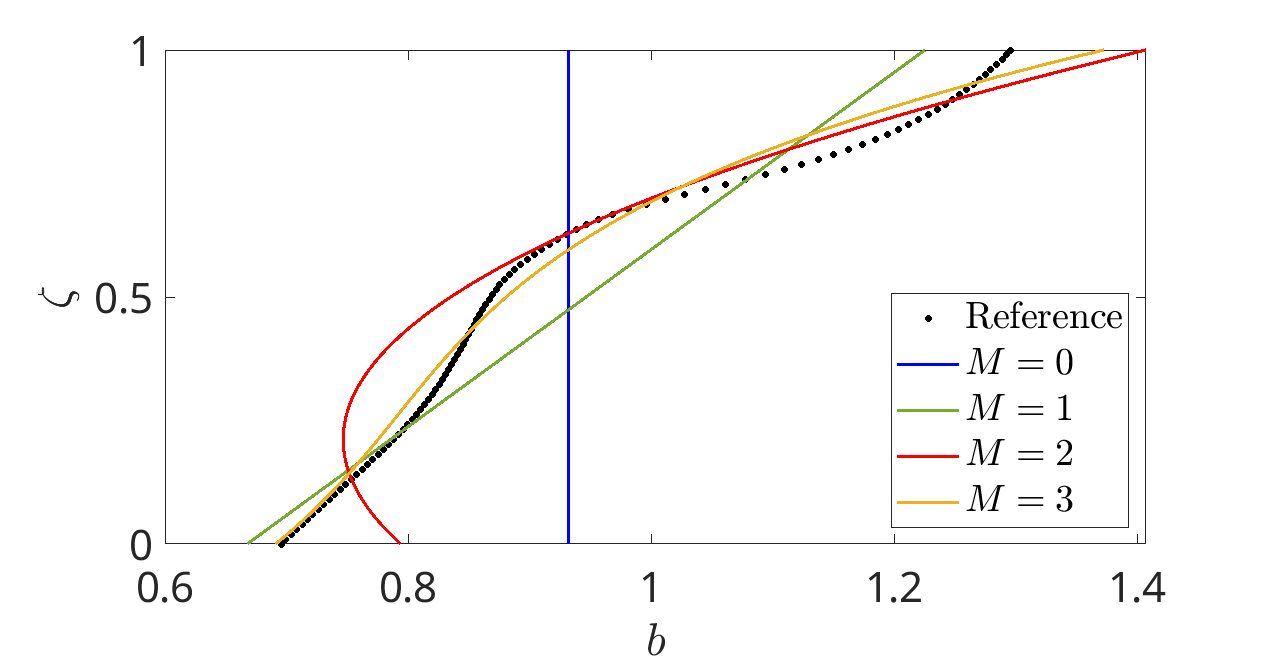}}
\caption{\sf Example 2 (Linear case): Vertical profiles of $v$ (left) and $b$ (right) at $y = -2/5$ of the reference solution and the $M$th-order MRSWME, $M = 0,\dots,3$.}
\label{fig:bumpB_prof}
\end{figure}

To assess the capability of how the moment models approximate the reference system, we first present the 1-D first-, second-, and third-order MRSWME against the reference solution for the linear case.
For the initial conditions of all moment models of order $M$, we take $\beta_1 = \eta_1 = -\frac{1}{4}$ and all other moments equal to zero, to match the linear profile initial conditions of the reference solution.
The comparisons against the reference solution for $h$ and mean velocity $v_m$ are shown in Figure \ref{fig:bumpB_comp}; likewise, we show the $v$ and $b$ vertical profiles at $y = -2/5$ of the first-, second-, and third-order MRSWME against that of the reference solution in Figure \ref{fig:bumpB_prof}.

\begin{figure}[ht!]
\centerline{
\includegraphics[trim=0.0cm 0.1cm 1.5cm 0.2cm, clip, width=8cm]{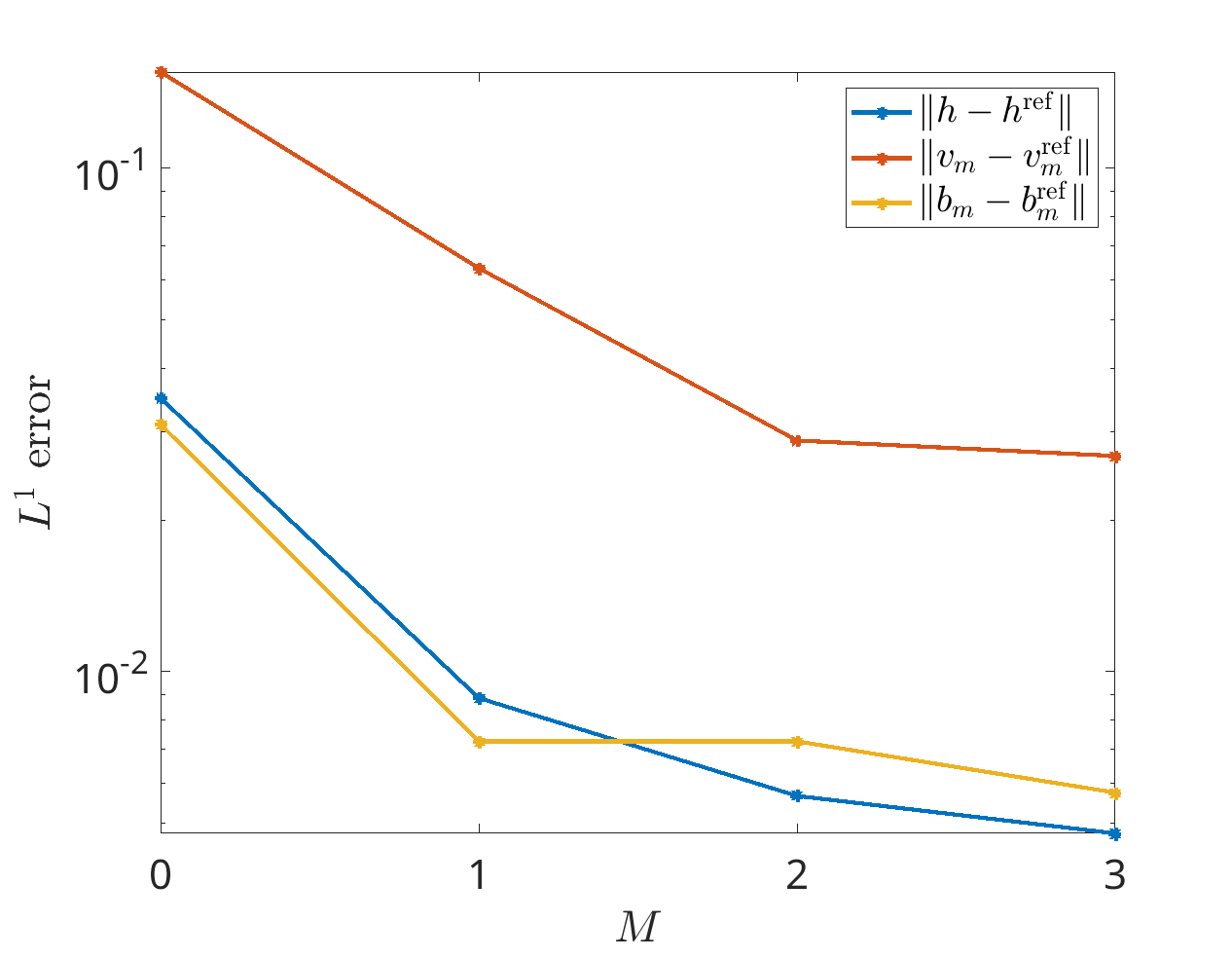}}
\caption{\sf Example 2 (Linear case): $L^1$ errors of $h$, $v_m$, and $b_m$ between moment models of order $M$ and the reference solution.}
\label{fig:bumpB_lin_err}
\end{figure}

In Figure \ref{fig:bumpB_comp}, we see that as $M$ increases, the moment model solution appears to approach the mean values of the reference solution, most clearly seen on the right traveling wave caused by the introduction of the linear vertical profile.
This is further confirmed by computing the $L^1$ error of $h$, $v_m$, and $b_m$ against the reference solution, presented in Figure \ref{fig:bumpB_lin_err} for the different moment models.
Further `visual convergence' can be seen in the vertical profile comparison in Figure \ref{fig:bumpB_prof}, especially in that for the magnetic field $b$. 

 \medskip



\begin{figure}[ht!]
\centerline{
\includegraphics[trim=0.9cm 0.1cm 1.7cm 0.2cm, clip, width=8cm]{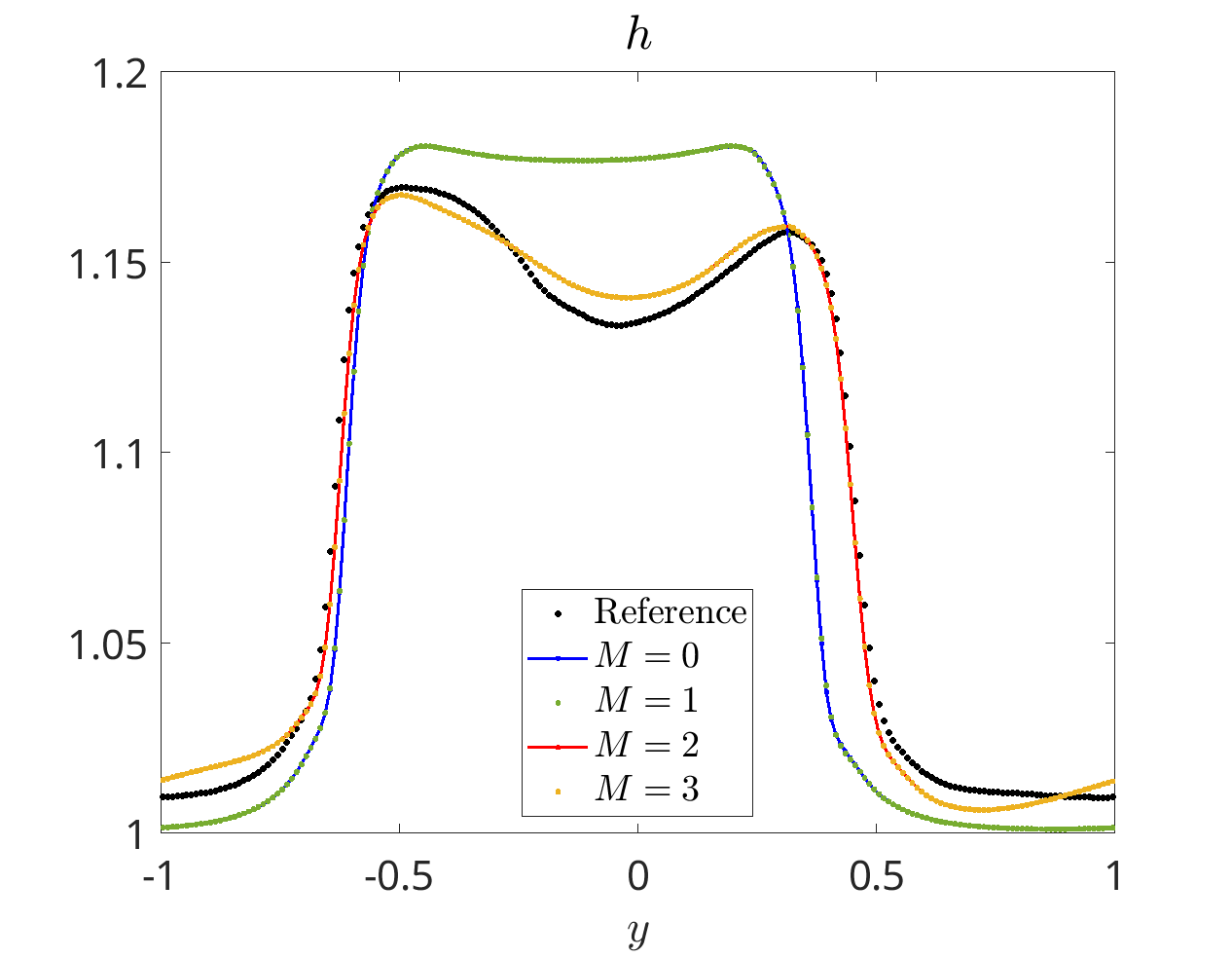}
\hspace{10pt}
\includegraphics[trim=0.9cm 0.1cm 1.7cm 0.2cm, clip, width=8cm]{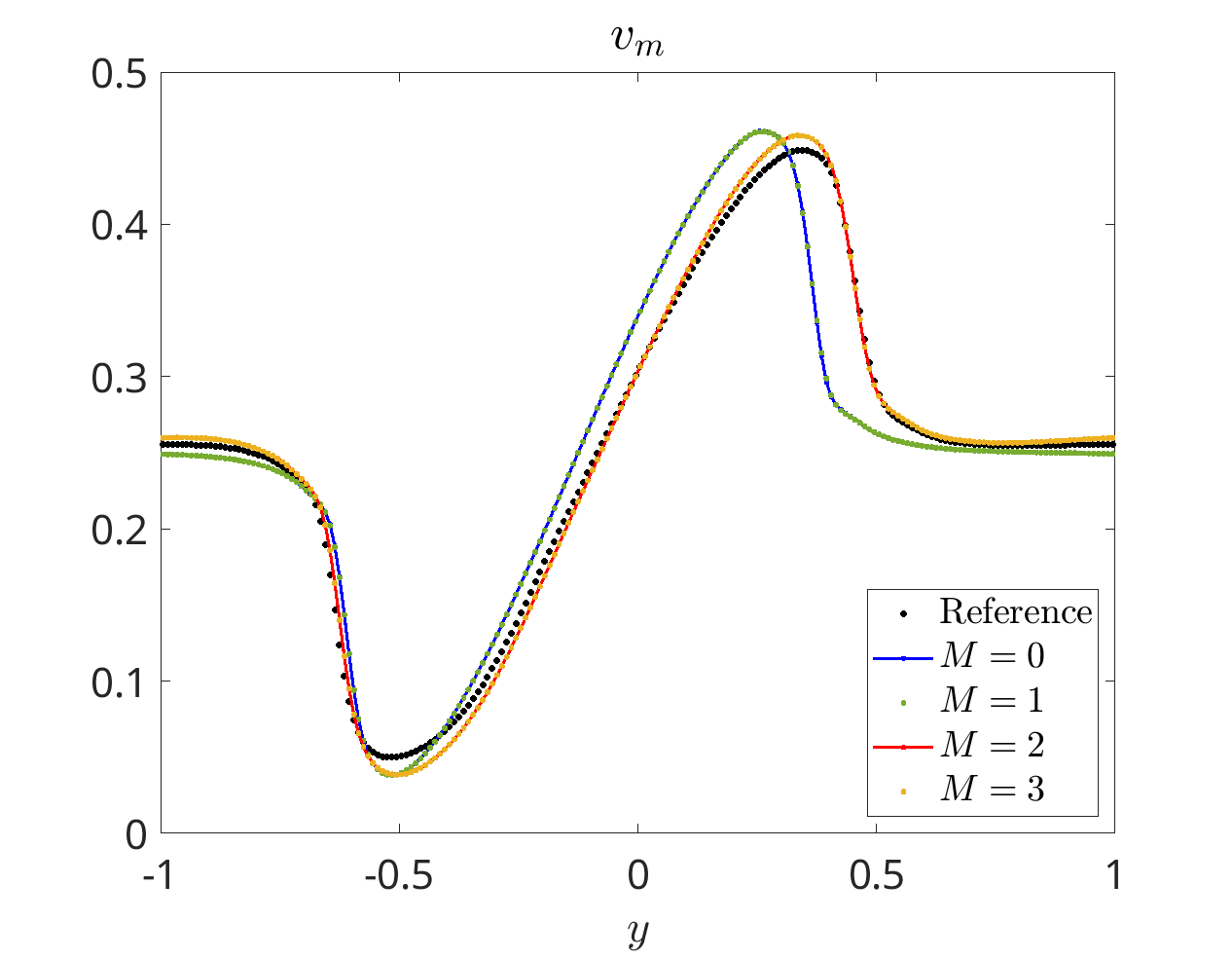}}
\caption{\sf Example 2 (Quadratic case): The solutions for $h$ (left) and $v_m$ (right) of the $Mth$-order MRSWME, $M = 0,\dots,3,$ against the reference solution.}
\label{fig:bumpB_comp_quad}
\end{figure}

\begin{figure}[ht!]
\centerline{
\includegraphics[trim=0.6cm 0.1cm 1.5cm 0.2cm, clip, width=8cm]{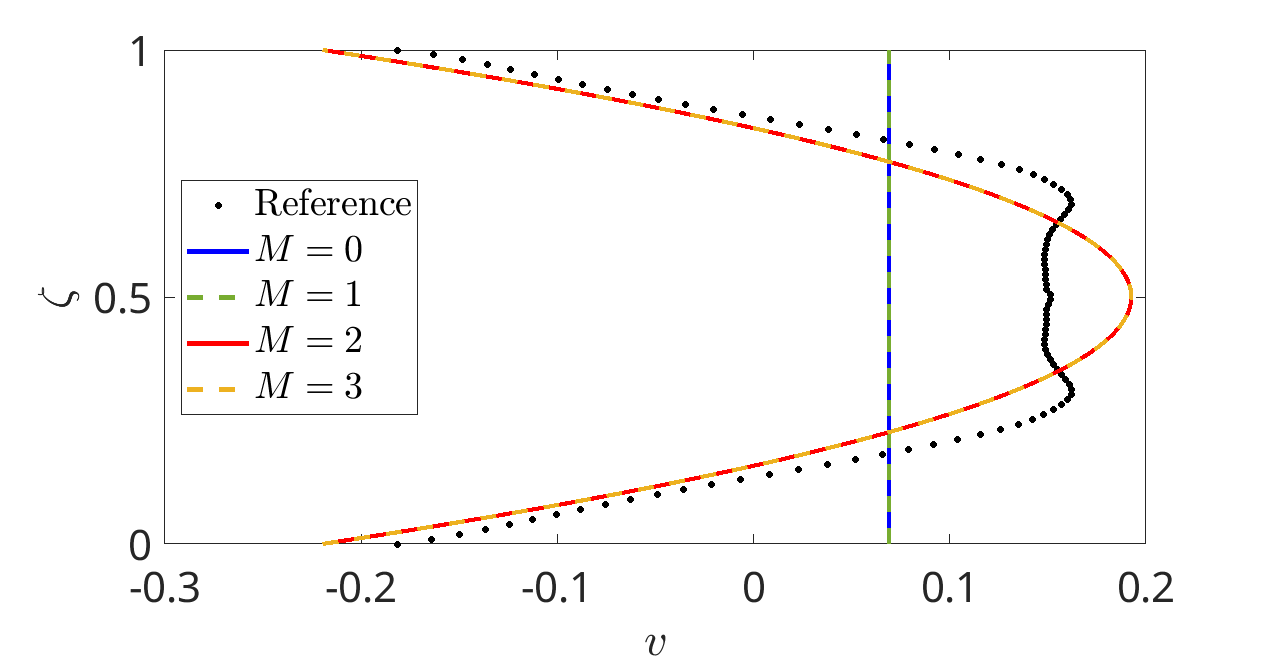}
\hspace{10pt}
\includegraphics[trim=0.6cm 0.1cm 1.5cm 0.2cm, clip, width=8cm]{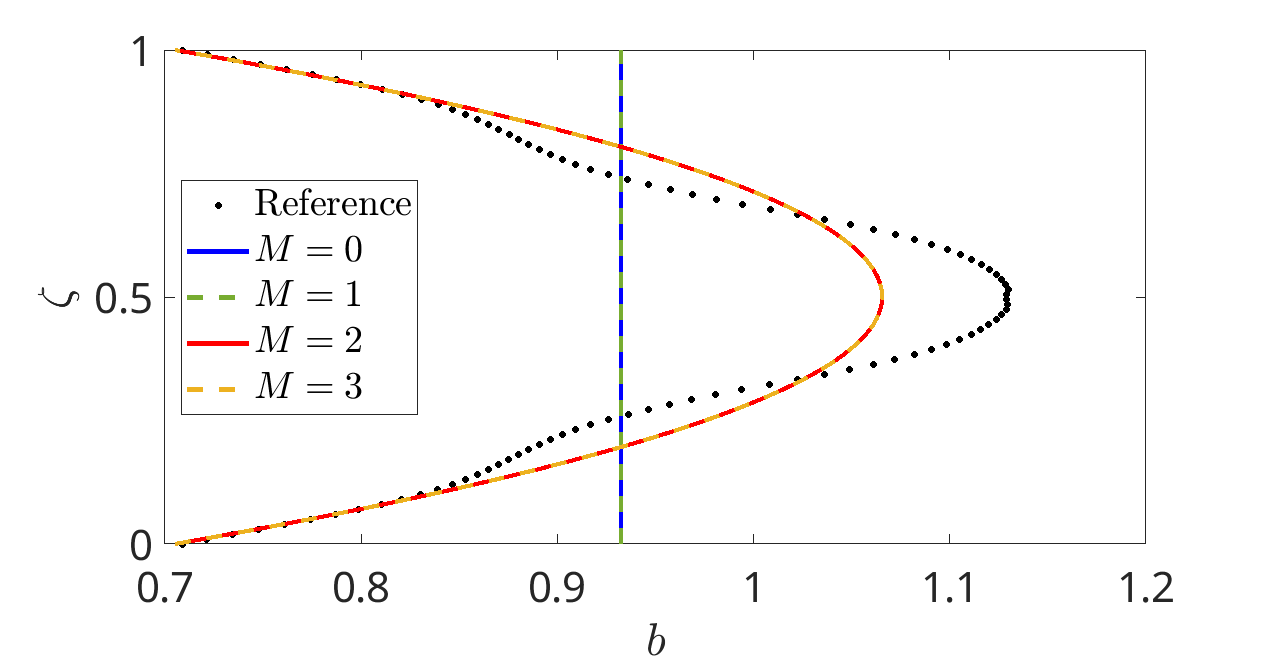}}
\caption{\sf Example 2 (Quadratic case): Vertical profiles of $v$ (left) and $b$ (right) at $y = -2/5$ of the reference solution and the $M$th-order MRSWME, $M = 0,\dots,3$.}
\label{fig:bumpB_prof_quad}
\end{figure}

Next, we compare the results of the MRSWME against the reference solution for the case in which the initial $v$ and $b$ vertical profiles are both quadratics.
Here, the initial conditions for the moment models take $\beta_2 = \eta_2 = -\frac{1}{4}$ and all other moments zero; thus, the $M = 1$ case reduces to the $M = 0$ result, as the linear coefficients are zero. 
The comparison of mean values is presented in Figure \ref{fig:bumpB_comp_quad}, and the vertical profile comparison is shown in Figure \ref{fig:bumpB_prof_quad}.

\begin{figure}[ht!]
\centerline{
\includegraphics[trim=0.0cm 0.1cm 1.5cm 0.2cm, clip, width=8cm]{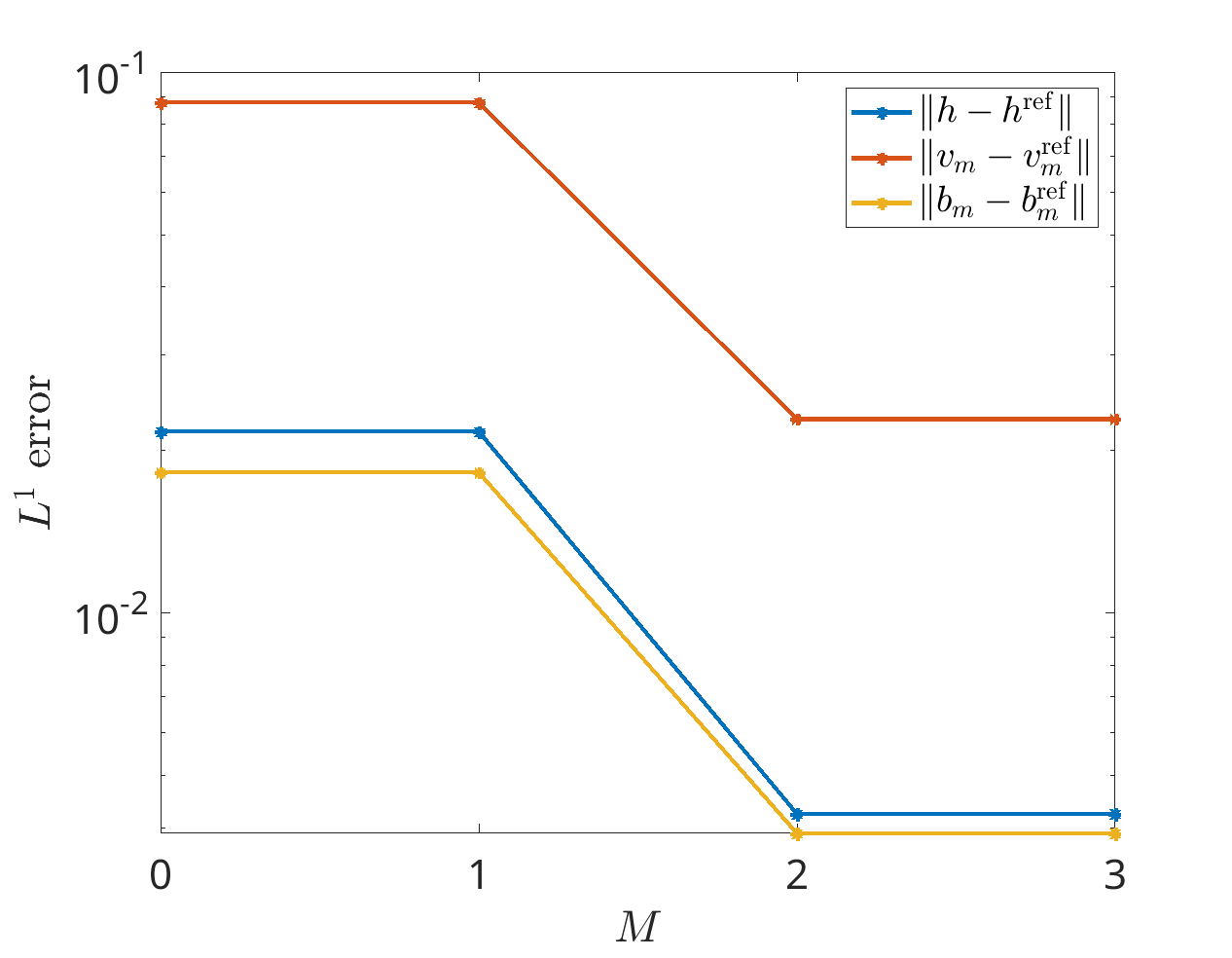}}
\caption{\sf Example 2 (Quadratic case): $L^1$ errors of $h$, $v_m$, and $b_m$ between moment models of order $M$ and the reference solution.}
\label{fig:bumpB_quad_err}
\end{figure}

Unlike previous work on the SWME, we do not consider friction at the bottom of the domain.
As a result, the vertical profiles remain symmetric about $\zeta = \frac{1}{2}$, as seen in Figure \ref{fig:bumpB_prof_quad}.
Furthermore, in a symmetric vertical profile case such as this one, any moment models of order $2M+1$ reduce to that of order $2M$, hence why we only see overlapping curves for the MRSWME in Figures \ref{fig:bumpB_comp_quad} and \ref{fig:bumpB_prof_quad}.
This is further confirmed when computing the $L^1$ errors of $h$, $v_m$, and $b_m$ in Figure \ref{fig:bumpB_quad_err}, where we see the error only reduces when increasing to an even-order moment model. 
However, the error reduction at the even-order model is significant, as all variables see an error reduction by a factor of $\gtrsim 4$.
While not presented in this work, we expect this step-like error reduction seen in Figure \ref{fig:bumpB_quad_err} to continue as the order of the MRSWME is increased. 


\medskip

\begin{figure}[ht!]
\centerline{
\includegraphics[trim=0.9cm 0.1cm 1.7cm 0.2cm, clip, width=8cm]{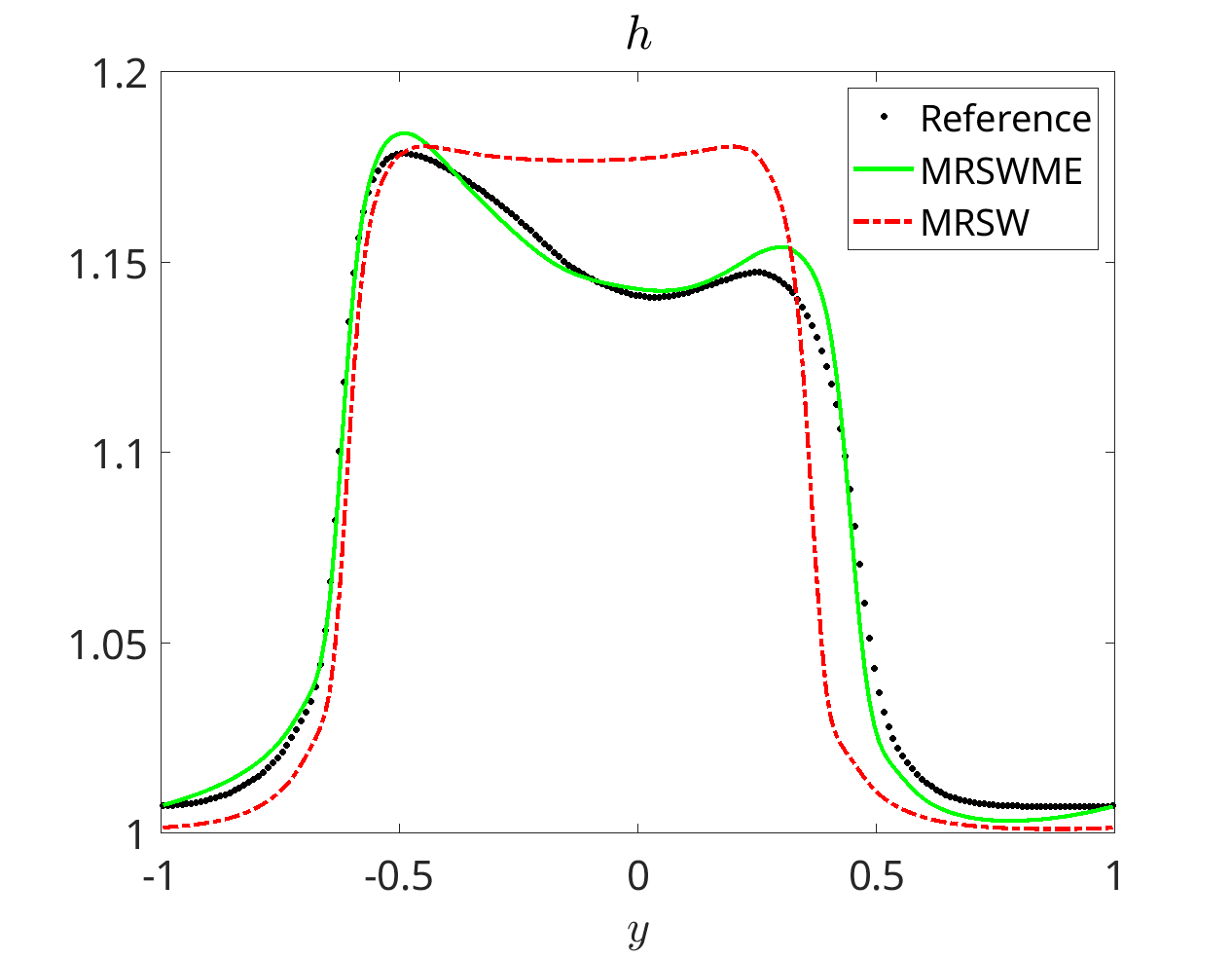}
\hspace{10pt}
\includegraphics[trim=0.9cm 0.1cm 1.7cm 0.2cm, clip, width=8cm]{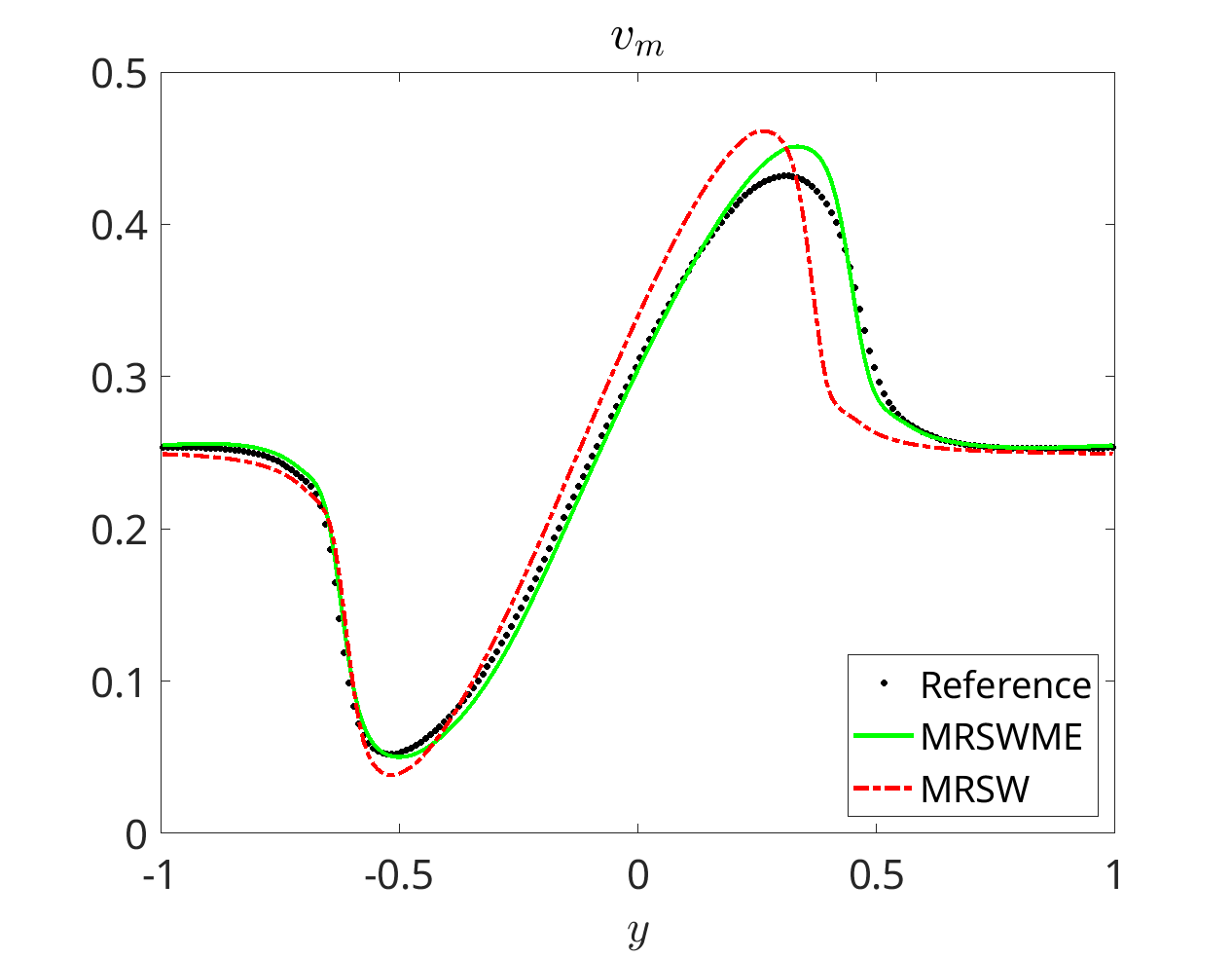}}
\caption{\sf Example 2 (Cubic case): The solutions for $h$ (left) and $v_m$ (right) of the MRSW system and the third-order MRSWME against the reference solution.}
\label{fig:bumpB_comp_cub}
\end{figure}

\begin{figure}[ht!]
\centerline{
\includegraphics[trim=0.6cm 0.1cm 1.5cm 0.2cm, clip, width=8cm]{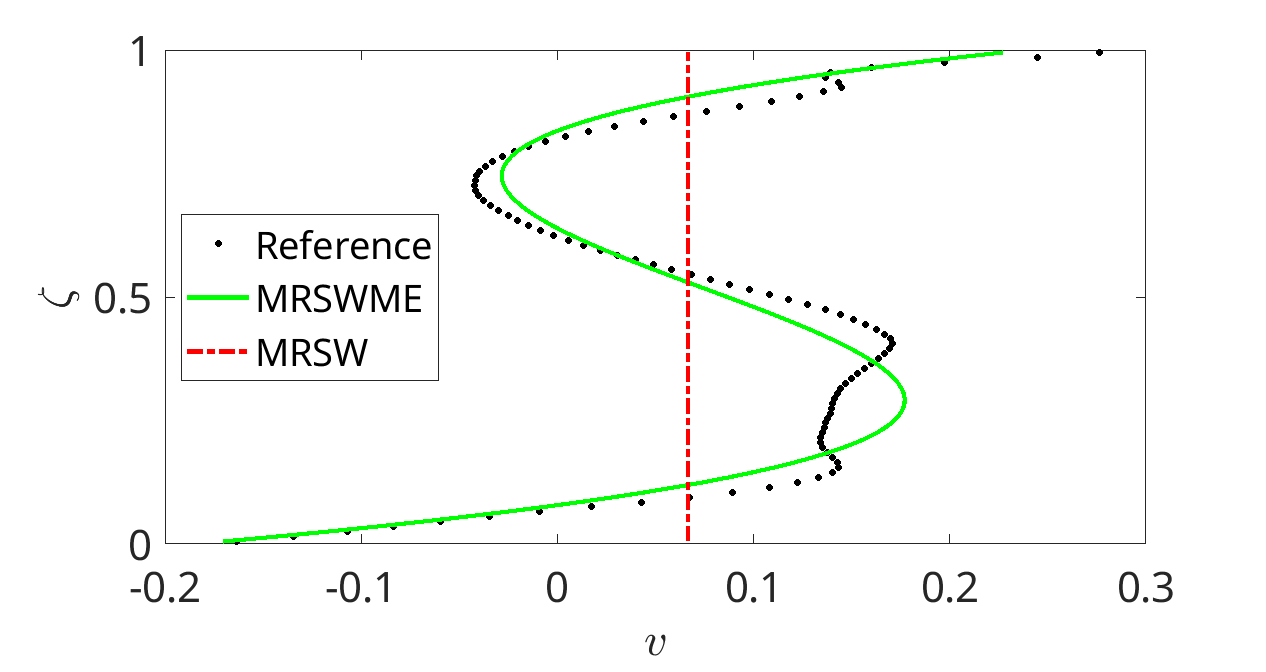}
\hspace{10pt}
\includegraphics[trim=0.6cm 0.1cm 1.5cm 0.2cm, clip, width=8cm]{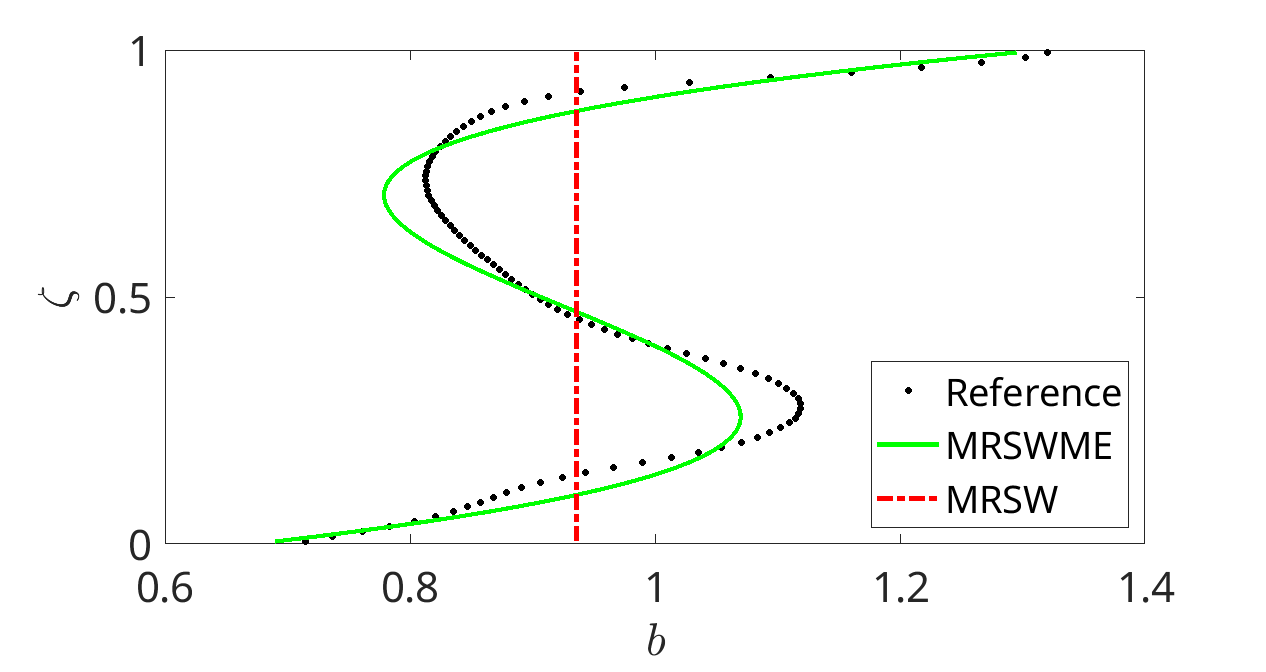}}
\caption{\sf Example 2 (Cubic case): Vertical profiles of $v$ (left) and $b$ (right) at $y = -2/5$ of the reference solution, MRSW system, and the third-order MRSWME.}
\label{fig:bumpB_prof_cub}
\end{figure}

Lastly for this example, we present the third-order moment model against the reference solution with initially cubic vertical profiles for $v$ and $b$; that is, we take the $\beta_3 = \eta_3 = -\frac{1}{4}$ initially.
The $h$ and $v_m$ results are presented in Figure \ref{fig:bumpB_comp_cub}, with the profiles of $v$ and $b$ at $y = -2/5$ presented in Figure \ref{fig:bumpB_prof_cub}.
Much like the linear and quadratic profile cases, we see that the third-order MRSWME better approximates the reference solution in comparison to the MRSW---capturing the wave speeds, wave amplitudes, and vertical profiles much more accurately than the standard MRSW system.

\subsection{Geostrophic adjustment}\label{sec6.1}
For the remaining examples, we study the dynamical regimes of magneto-geostrophic adjustment problems, and the influence of vertical velocity and magnetic field profiles on these problems. 
The (magneto-)geostrophic adjustment process describes the evolution toward the corresponding (magneto-)geostrophic equilibrium---which is of particular importance in astrophysical and geophysical applications; see, e.g., \cite{REZNIK2001geostr,Tobias2007MHDTurb,Zeitlin2015Geostrophic}. 
The steady-state magneto-geostrophic equilibria for the MRSW equations read
\begin{equation*}\label{eq:mgeo_eq}
	fu = -g(h+Z)_y+ bb_y ,\qquad v = 0.
\end{equation*}
In the absence of a magnetic field, it has been shown that any non-trivial initial condition evolves toward the state of corresponding geostrophic equilibrium through emission of inertia-gravity waves; see \cite{zeitlin2018geophysical} and references therein.
However, the presence of the magnetic field makes the magneto-geostrophic adjustment process more complicated, as the additional equation for $a_m$ adds a restriction on the equilibrium solution of $u_m$; see, e.g., \cite{CHERTOCK2024WB}.
More so, the addition of the vertical profiles complicates the (magneto-)geostrophic equilibria further, as expansion terms remain present in the equilibrium solution---even when taking $v = v_m = 0$ everywhere along the vertical profile.

Recall that the regimes near (magneto-)geostrophic equilibria are characterized by small Rossby and magnetic Rossby numbers
$$
	{\rm Ro} = \frac{\mathfrak U}{f \mathfrak L}, \qquad 
	{\rm Ro}_m = \frac{\mathfrak B}{f \mathfrak L},
$$
where $\mathfrak U$ and $\mathfrak B$ are the scales of velocity and magnetic field, respectively, and $\mathfrak L$ is the typical length scale of the motions under study. 
In the following examples, we look at two cases: Example 3 considers a low Rossby number (${\rm Ro}<1$ and ${\rm Ro}_m < 1$), and Example 4 considers a high Rossby number (${\rm Ro}>1$ and ${\rm Ro}_m > 1$).
In addition, both cases will have an initial sinusoidal vertical profile in the $y$-velocity, allowing us to analyze the influence of this vertical profile on (i) the magneto-geostrophic solution; and (ii) how it introduces a vertical profile in the magnetic field.

\subsubsection*{Example 3---Magneto-geostrophic adjustment at low Rossby numbers}\label{ex3}
The first example considering magneto-geostrophic adjustment is that with low Rossby numbers; that is, with ${\rm Ro}={\rm Ro}_m= 0.1$; in which we analyze the influence of an initial vertical perturbation in $v$.
This problem remains in the dynamical regime that is near magneto-geostrophic equilibrium and the outward-moving waves are expected to be smooth; this remains true with our choice of vertical profile. 
For this problem, we take the initial conditions in which the vertical profile of $v$ is described by a sine wave; that is, we take
$$
    h(y,0) = 1, \qquad 
    u(y,\zeta, 0) = 0.1e^{-y^2}, \qquad
    v(y,\zeta,0) = \frac{1}{4}\sin(2\pi\zeta),
    \qquad 
    b(y,\zeta,0) = 0.1,
$$
with $a(y,\zeta,0) = 0$, a constant Coriolis parameter $f(y) = 1$, and zero bottom topography $Z(y) = 0$ on the domain $[-20,20]$ subject to outflow boundary conditions. 
Notice that the vertical profile of $v$ has a mean value of zero, agreeing with the MRSW initial state of \cite{CHERTOCK2024WB}.
In addition, since the Coriolis term is non-zero, we also expect to see changes to the $x$-direction velocity and magnetic field as the solution is advanced in time.

To appropriately choose the profile coefficients needed for the initial state of the MRSWME, 
we use \eqref{eq:moment} to compute these coefficients exactly.
Doing so gives the following initial conditions for the MRSWME:
$$
    h(y,0) = 1, 
    \qquad u(y,\zeta,0) = 0.1e^{-y^2},
    \qquad v(y,\zeta,0) = \sum_{\ell = 1}^M\beta_\ell(y,0) \phi_\ell(\zeta), \qquad 
    b(y,\zeta,0) = 0.1,
$$
%
with $a = 0$, 
and after using \eqref{eq:moment} for $M = 1,\dots, 3$, the initial coefficients for the vertical velocity profile read
\begin{equation} \label{eq:beta_ic}
    \beta_1(y,0) = \frac{3}{4\pi}, \qquad 
    \beta_2(y,0) = 0, \qquad 
    \beta_3(y,0) = \frac{7}{4\pi^3}(\pi^2 -15).
\end{equation}

\begin{figure}[ht!]
\centerline{
\includegraphics[trim=0.9cm 0.1cm 1.6cm 0.2cm, clip, width=4.5cm]{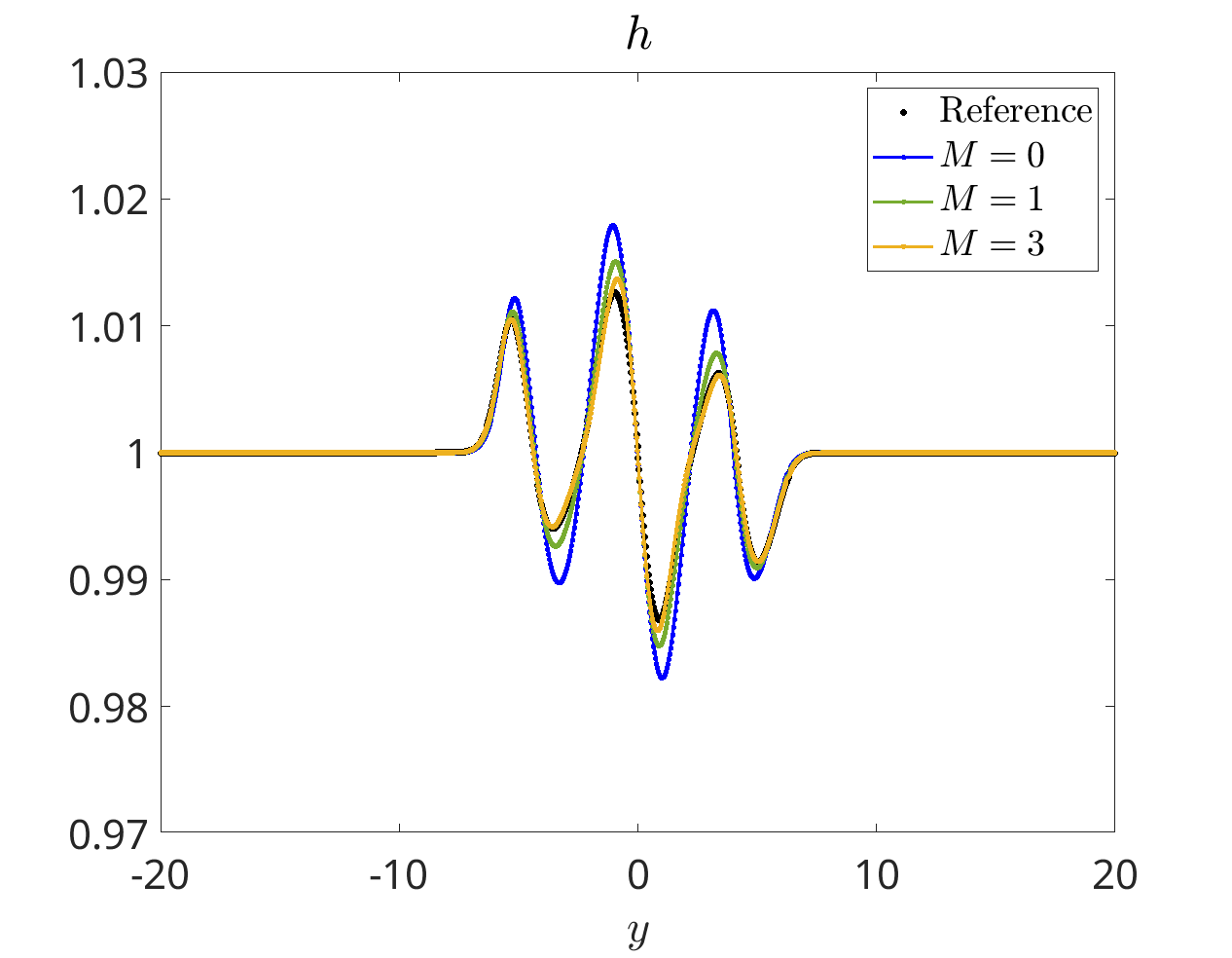}
\includegraphics[trim=0.9cm 0.1cm 1.6cm 0.2cm, clip, width=4.5cm]{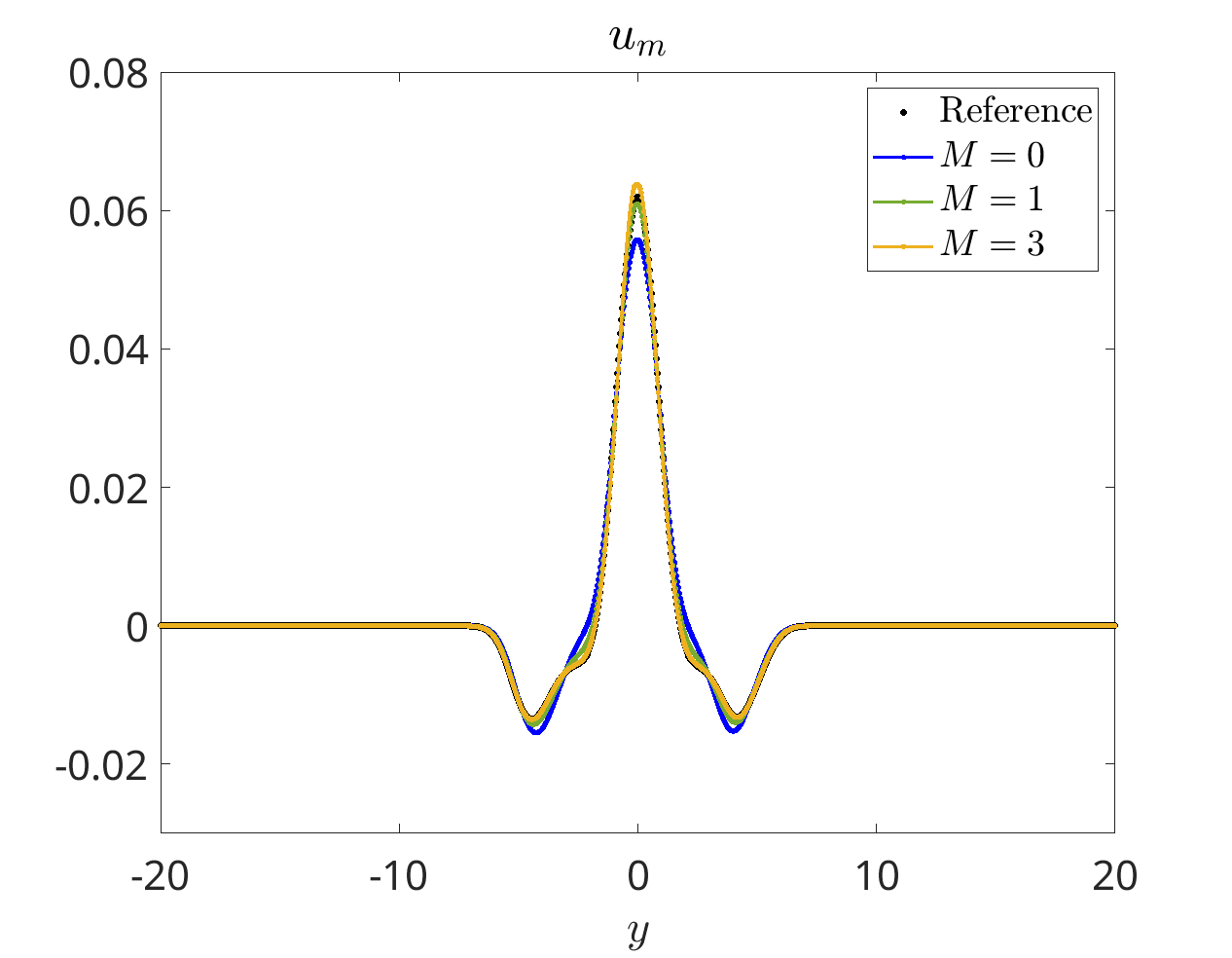}
\includegraphics[trim=0.9cm 0.1cm 1.6cm 0.2cm, clip, width=4.5cm]{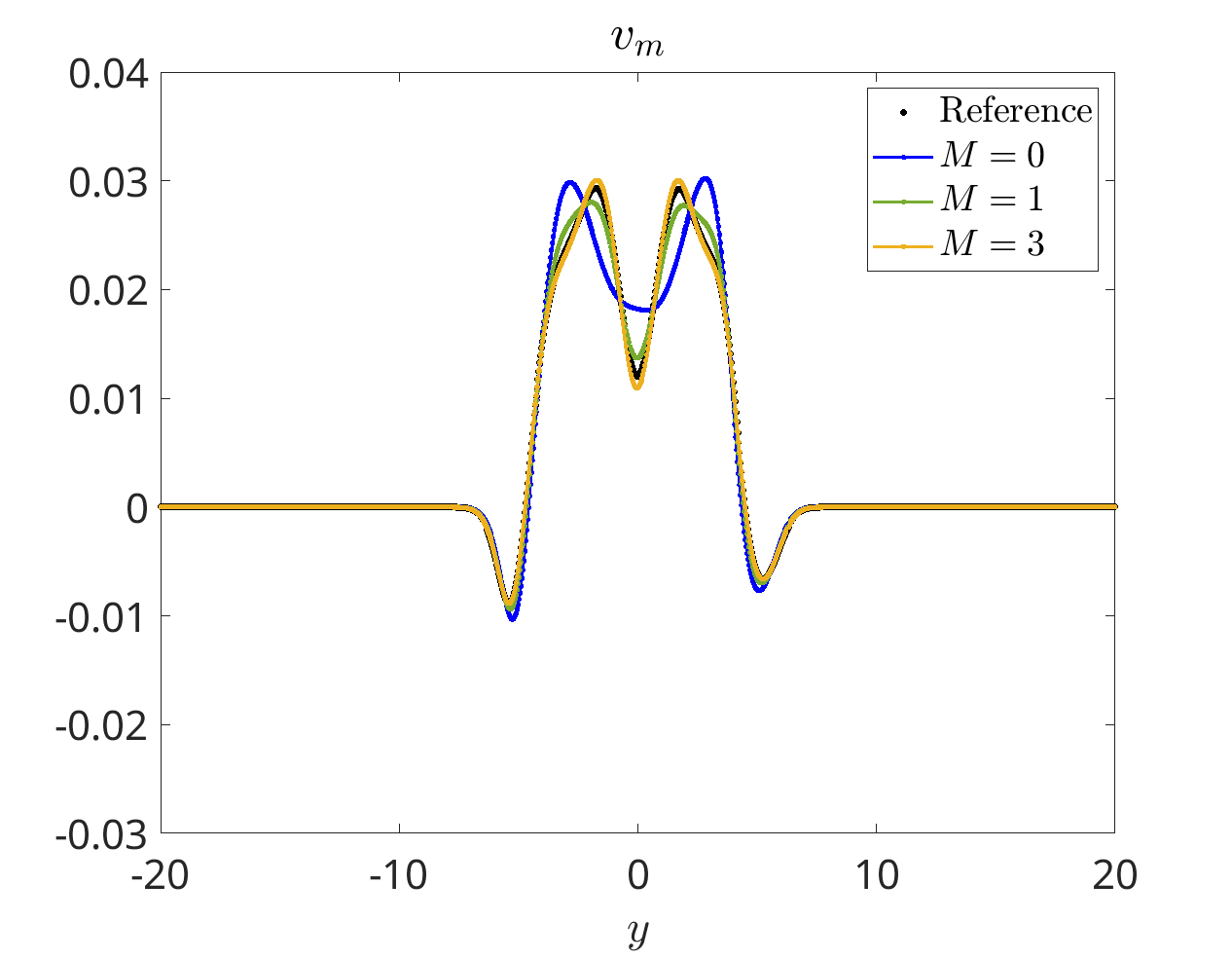}
\includegraphics[trim=0.9cm 0.1cm 1.6cm 0.2cm, clip, width=4.5cm]{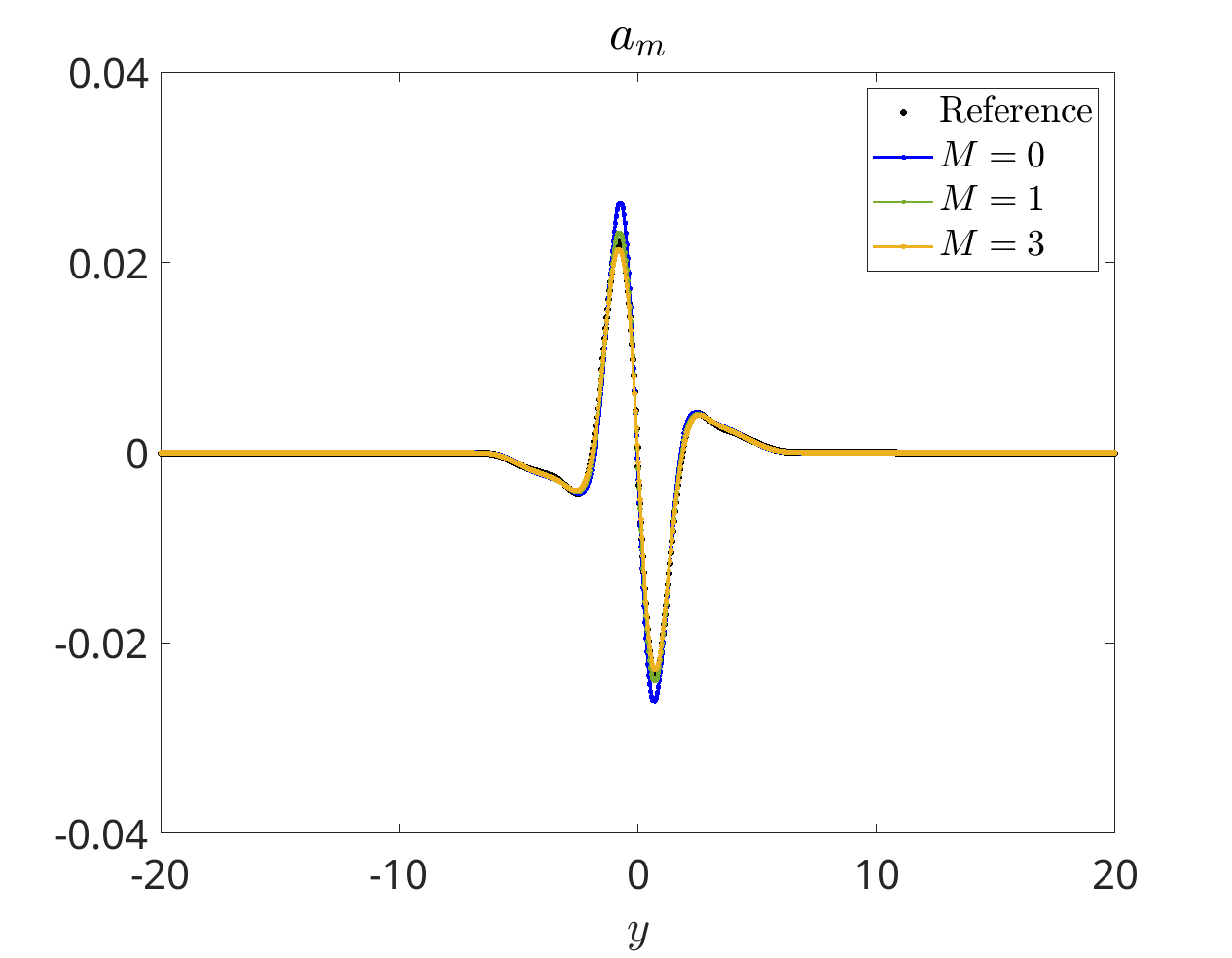}}
\vspace{7pt}
\centerline{
\includegraphics[trim=0.9cm 0.1cm 1.6cm 0.2cm, clip, width=4.5cm]{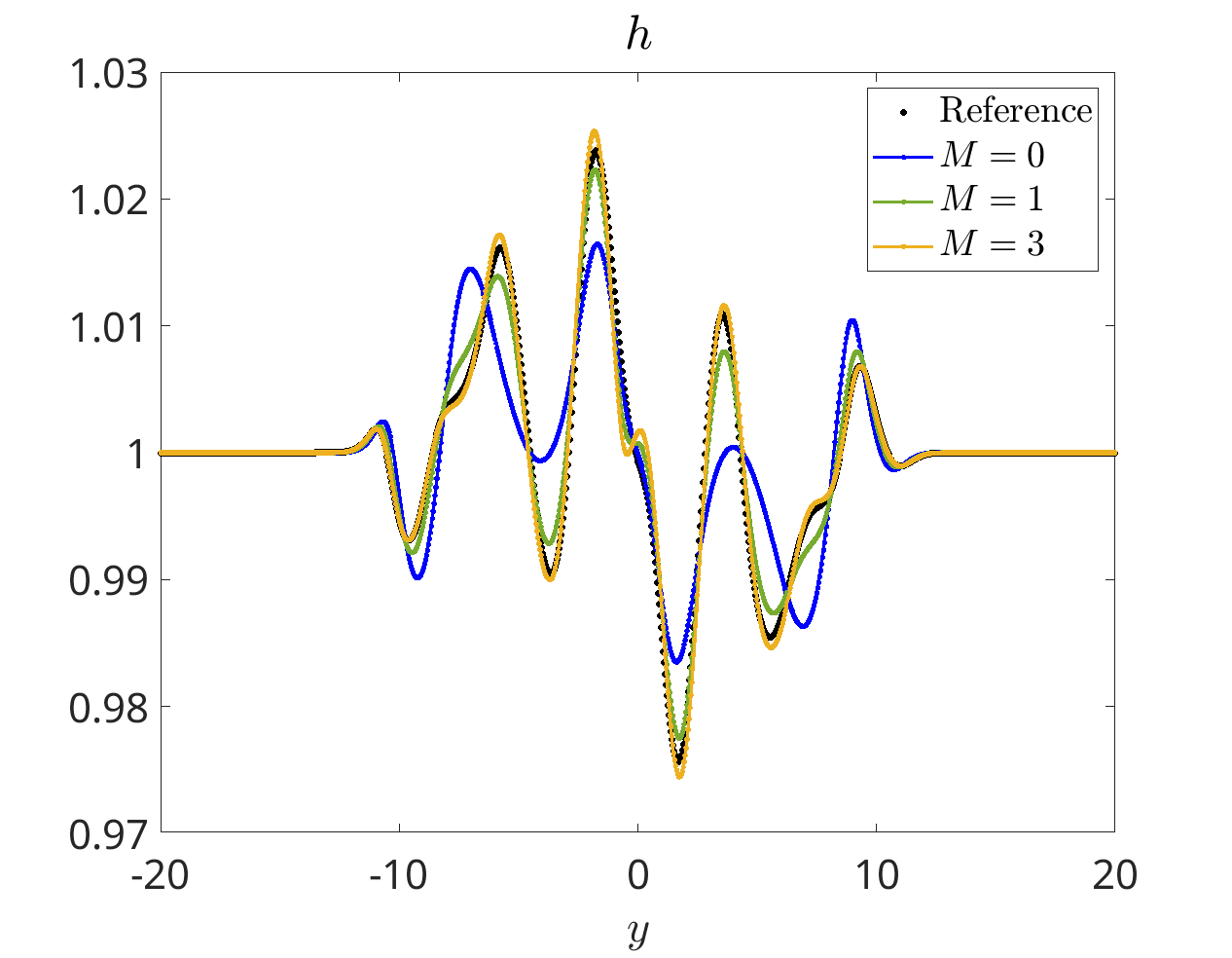}
\includegraphics[trim=0.9cm 0.1cm 1.6cm 0.2cm, clip, width=4.5cm]{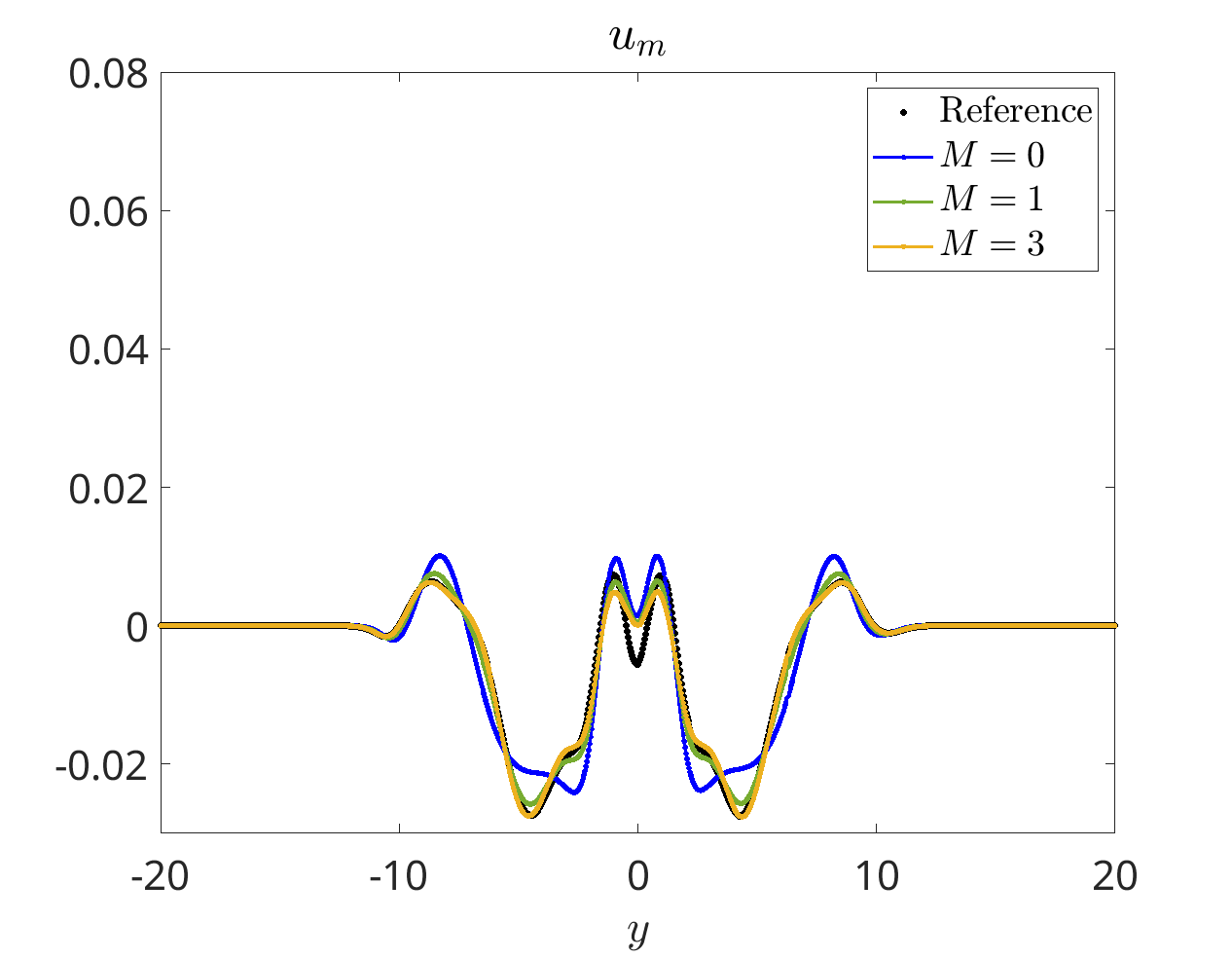}
\includegraphics[trim=0.9cm 0.1cm 1.6cm 0.2cm, clip, width=4.5cm]{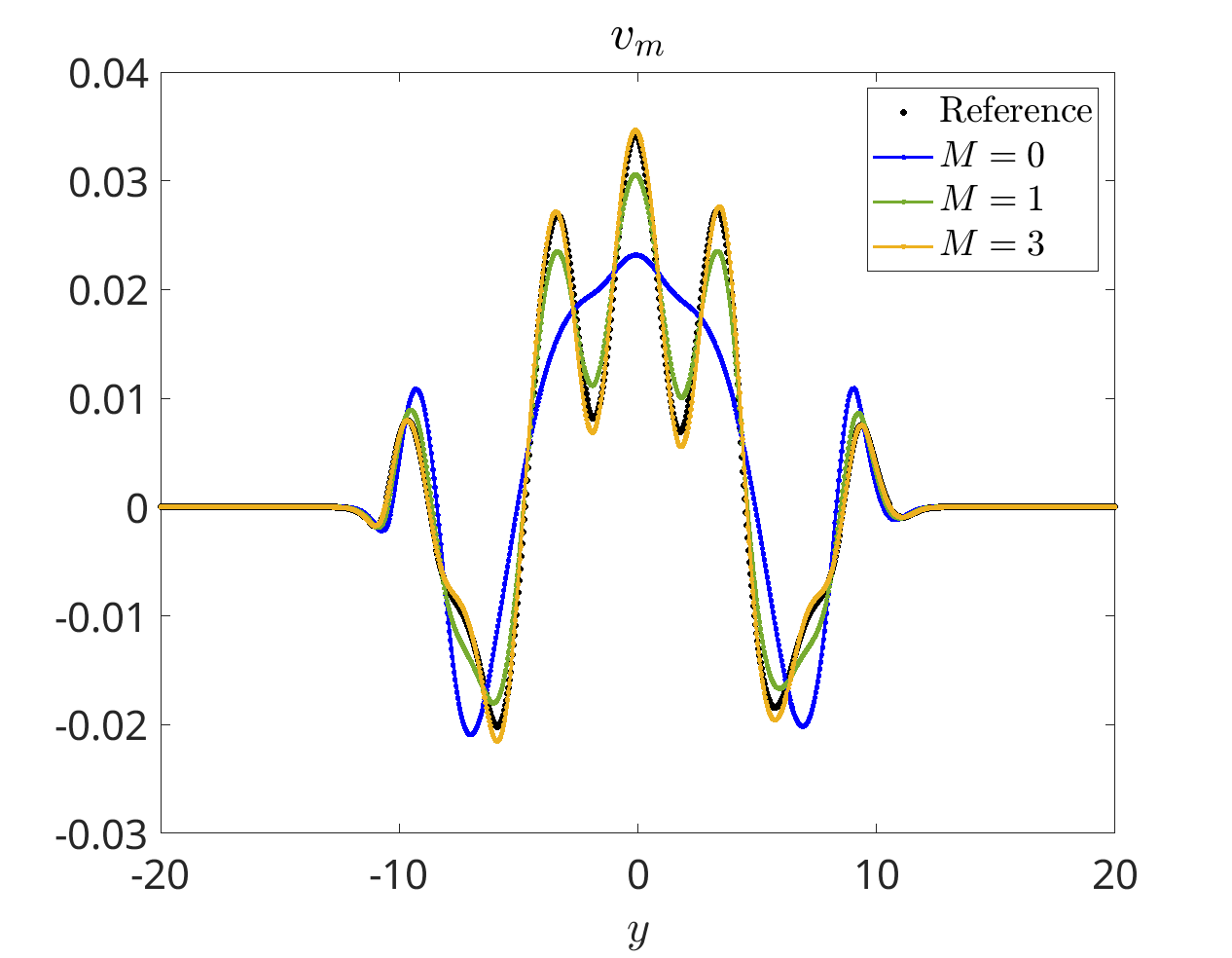}
\includegraphics[trim=0.9cm 0.1cm 1.6cm 0.2cm, clip, width=4.5cm]{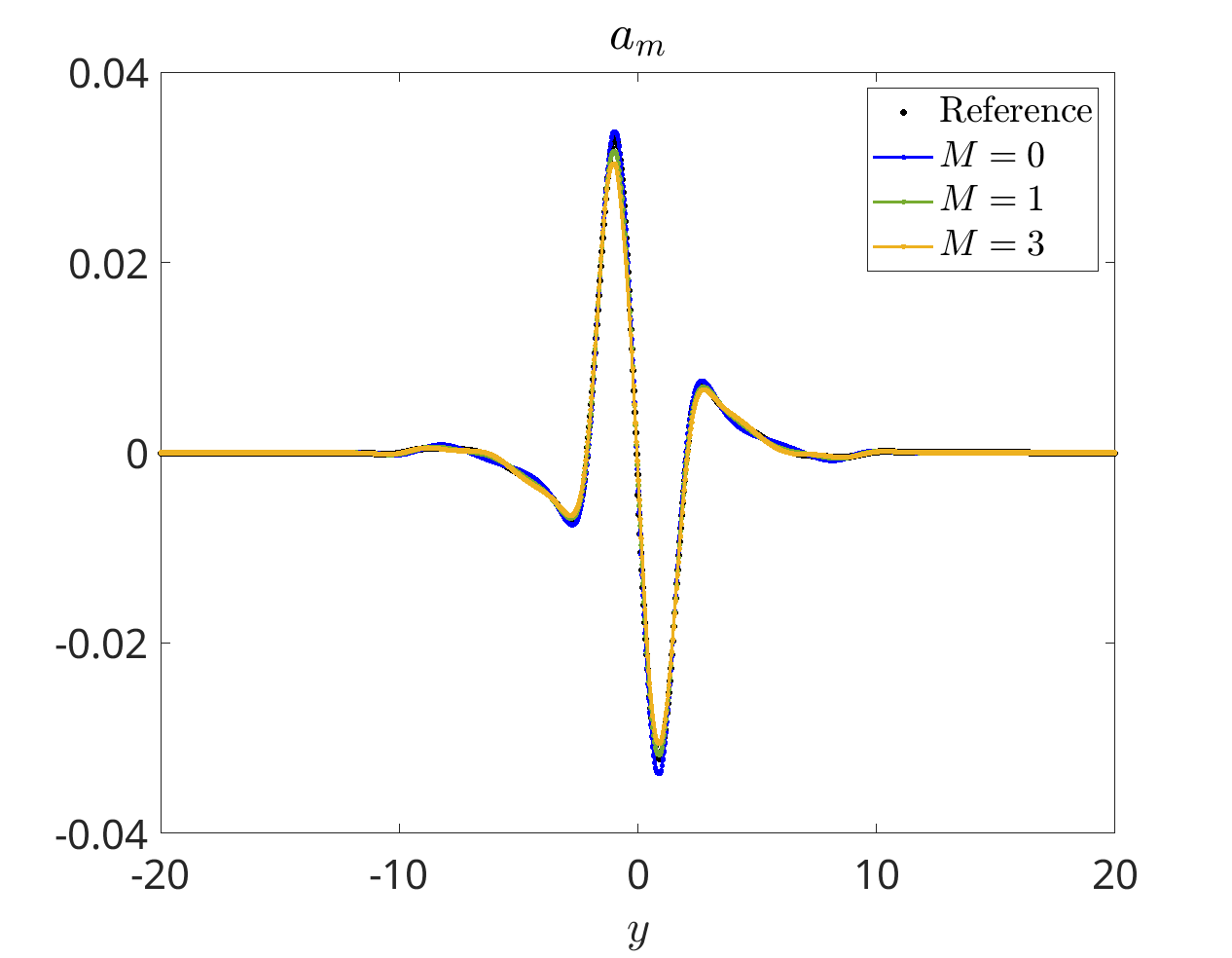}}
\caption{\sf Example 3: Solutions of $h$ (first column), $u_m$ (second column), $v_m$ (third column), and $a_m$ (last column) at times $t = 5$ (top row) and $t = 10$ (bottom row) of the reference system, MRSW ($M=0$), and MRSWME of order $M = 1$ and $M = 3$.}
\label{fig:Low_R_comp} 
\end{figure}

We present the solutions of $h$, $u_m$, $v_m$, and $a_m$ at times $t = 5$ and $t = 10$ for the reference system, MRSW, and MRSWME of order $M = 1$ and $M = 3$ in Figure \ref{fig:Low_R_comp}. 
Note that for the $M = 2$ MRSWME system, the initial conditions result in a non-hyperbolic state in which the complex portions of the eigenvalues are comparatively similar in magnitude to the real portions; solving such instances will be considered in future work. 
In this example, it is clear that outside of some amplitude difference, the standard MRSW captures the material wave in the center of the domain (in all variables except $v_m$), and the fast magnetogravity waves quite well. Ergo the small initial vertical perturbation in $v$ seemingly has little effect on these waves.
However, we do see significant differences in the Alfv\'en waves. 
Most clearly seen in the figure for $h$ and $v_m$, the inclusion of the initial sinusoidal vertical profile of $v$ causes these waves to slow, in addition to slightly changing the behavior of the material wave.
As expected, we see that the MRSWME captures the behavior of the Alfv\'en waves much more accurately in comparison to MRSW. 

\begin{figure}[ht!]
\centerline{
\includegraphics[trim=0.7cm 0.1cm 1.2cm 0.2cm, clip, width=6cm]{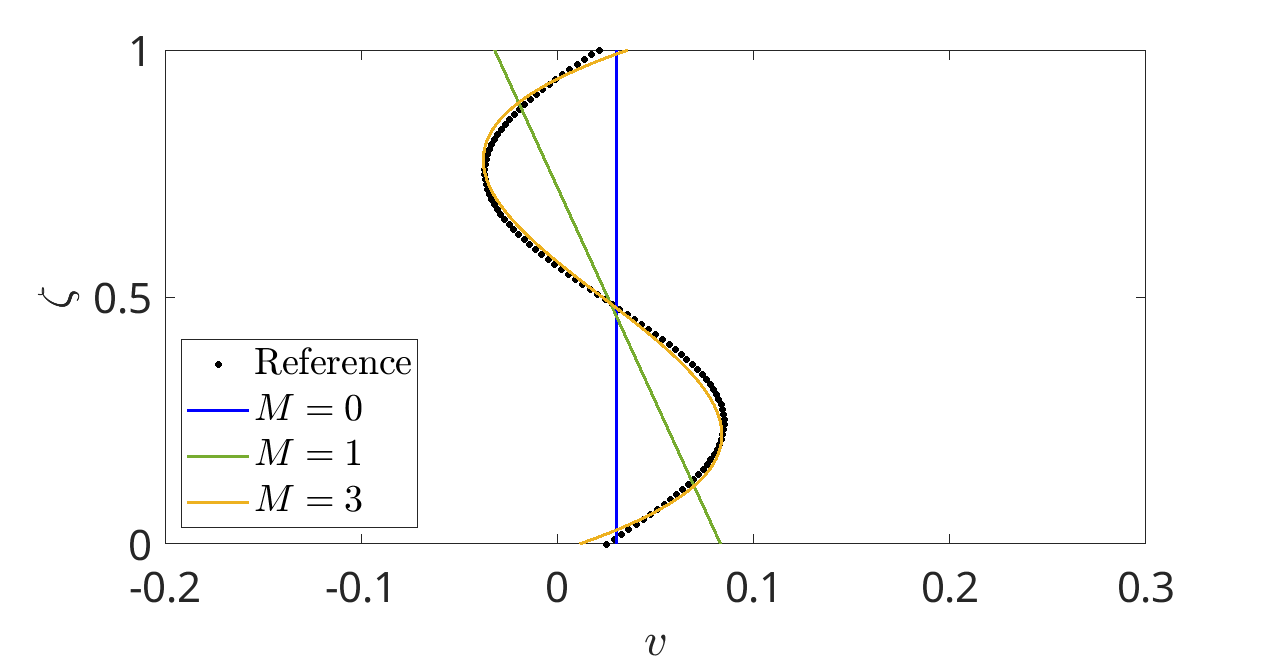}
\hspace{2pt}
\includegraphics[trim=0.7cm 0.1cm 1.2cm 0.2cm, clip, width=6cm]{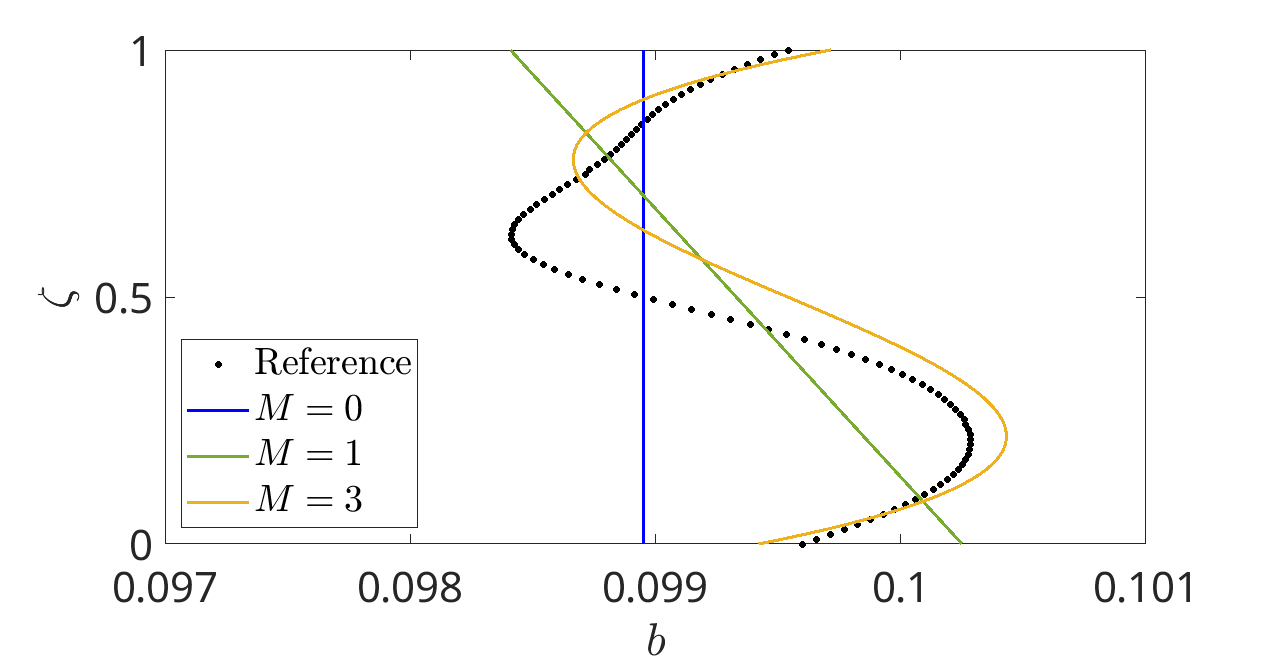}}
\vspace{7pt}
\centerline{
\includegraphics[trim=0.7cm 0.1cm 1.2cm 0.2cm, clip, width=6cm]{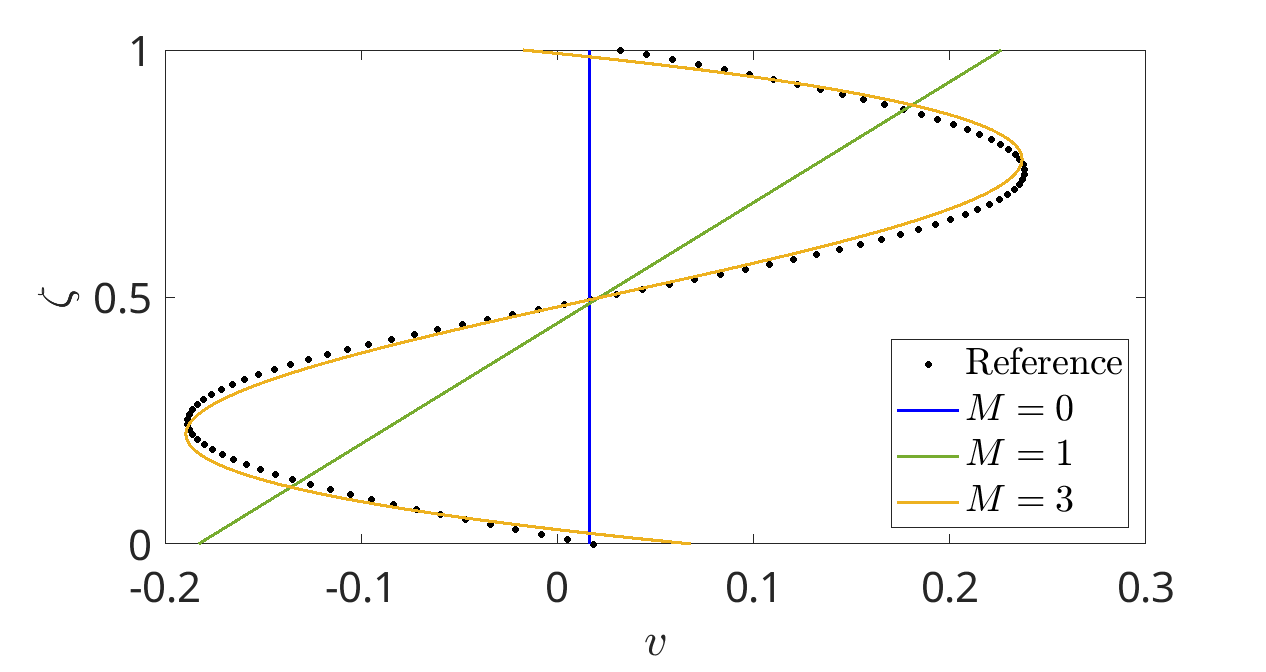}
\hspace{2pt}
\includegraphics[trim=0.7cm 0.1cm 1.2cm 0.2cm, clip, width=6cm]{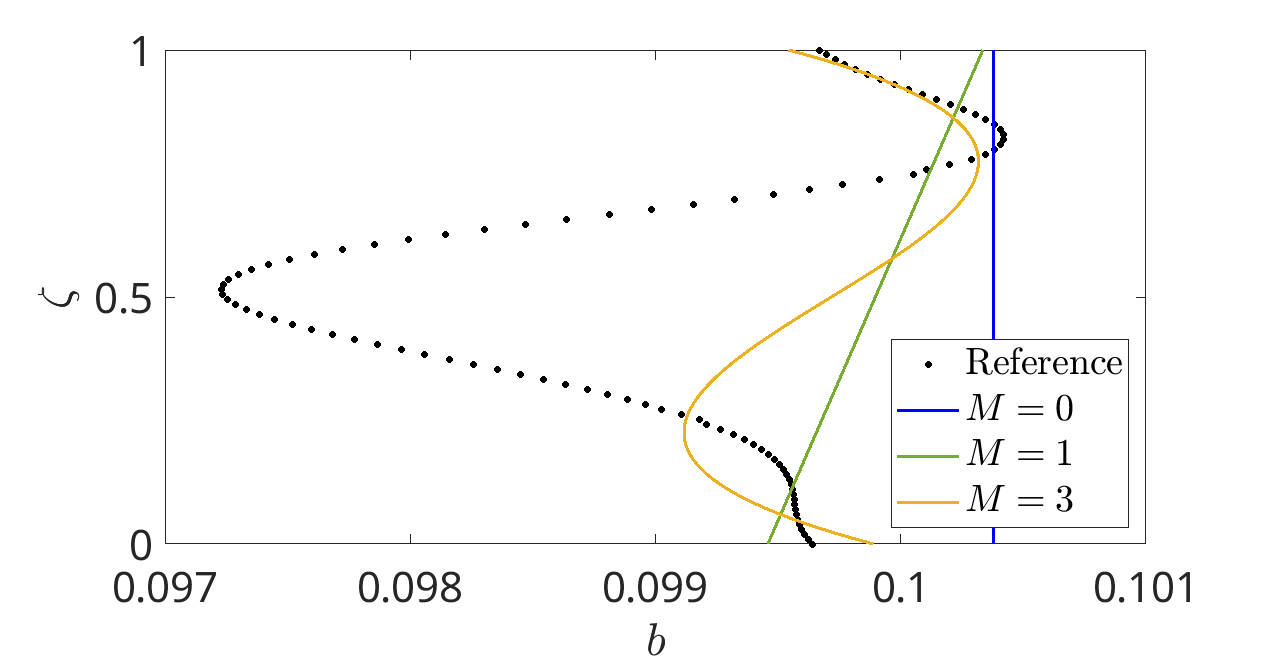}}
\caption{\sf Example 3: Vertical profiles at $y = -5$ of $v$ (left column) and $b$ (right column) at times $t = 5$ (top row) and $t = 10$ (bottom row) for the reference system, MRSW, and MRSWME of order $M = 1$ and $M = 3$.}
\label{fig:Low_R_prof} 
\end{figure}

This is again due the more accurate approximation of the vertical profiles; comparisons of the $v$ and $b$ profiles at slice $y = -5$ are presented in Figure \ref{fig:Low_R_prof}.
We notice that the initial sinusoidal profile of velocity remains as time advances in the reference solution, and the MRSWME captures the velocity profile quite accurately, especially in the third-order model.
Furthermore, this initial perturbation of the velocity vertical profile induces a vertical perturbation in the magnetic field that (i) is comparatively more complex in shape; and (ii) has a maximum relative perturbation of $\sim$5\% from the MRSW-computed $b_m$.
In the magnetic field, however, we see that MRSWME has some difficulty capturing this fully, even if more accurate than the MRSW approximation.
We expect this to improve as more moments are taken in the MRSWME.

%
%

\subsubsection*{Example 4---Magneto-geostrophic adjustment at high Rossby numbers}\label{ex4}
The last experiment is a magneto-geostrophic adjustment problem with high Rossby numbers (${\rm Ro} = {\rm Ro}_m = 1.1$), in which we additionally introduce a vertical profile in $v$ with mean $v_m = 0$.
Since the Rossby numbers are now larger than one, this problem lies outside of the dynamical regime near magneto-geostrophic equilibrium, resulting in shock formation in the outward-propagating waves. 
For this problem, we similarly take the initial conditions from \cite{CHERTOCK2024WB} and introduce a sinusoidal vertical perturbation in $v$; that is, we take

\begin{align*}
    &h(y,0) = 1, \qquad 
    u(y,\zeta, 0) = \frac{11}{10}\cdot\frac{(1+\tanh(4x+2))(1-\tanh(4x-2))}{(1+\tanh 2)^2},\\
    &v(y,\zeta,0) = \frac{1}{4}\sin(2\pi\zeta),
    \qquad 
    a(y,\zeta,0) = 0,\qquad
    b(y,\zeta,0) = 1.1,
\end{align*}
with a constant Coriolis parameter $f(y) = 1$ and a flat bottom topography $Z(y) = 0$ on the domain $[-20,20]$ subject to outflow boundary conditions.

For the 1-D MRSWME representation of this problem, the initial conditions read as follows: 
\begin{align*}
    &h(y,0) = 1, \qquad 
    u(y,\zeta, 0) = \frac{11}{10}\frac{(1+\tanh(4x+2))(1-\tanh(4x-2))}{(1+\tanh 2)^2},\\
    &v(y,\zeta,0) = \sum_{\ell = 1}^M\beta_\ell(y,0) \phi_\ell(\zeta),
    \qquad 
    a(y,\zeta,0) = 0,\qquad
    b(y,\zeta,0) = 1.1,
\end{align*}
where the initial $v$ vertical profile coefficients $\beta_\ell(y,0)$ are again initialized by applying the sinusoidal $v$ profile within \eqref{eq:moment}. 
Since the vertical profile in $v$ is identical to that in Example 3, \eqref{eq:moment} again results in the initial coefficients $\beta_\ell,\ \ell = 1,\dots,3,$ being those in \eqref{eq:beta_ic}.

\begin{figure}[ht!]
\centerline{
\includegraphics[trim=0.9cm 0.1cm 1.6cm 0.2cm, clip, width=4.5cm]{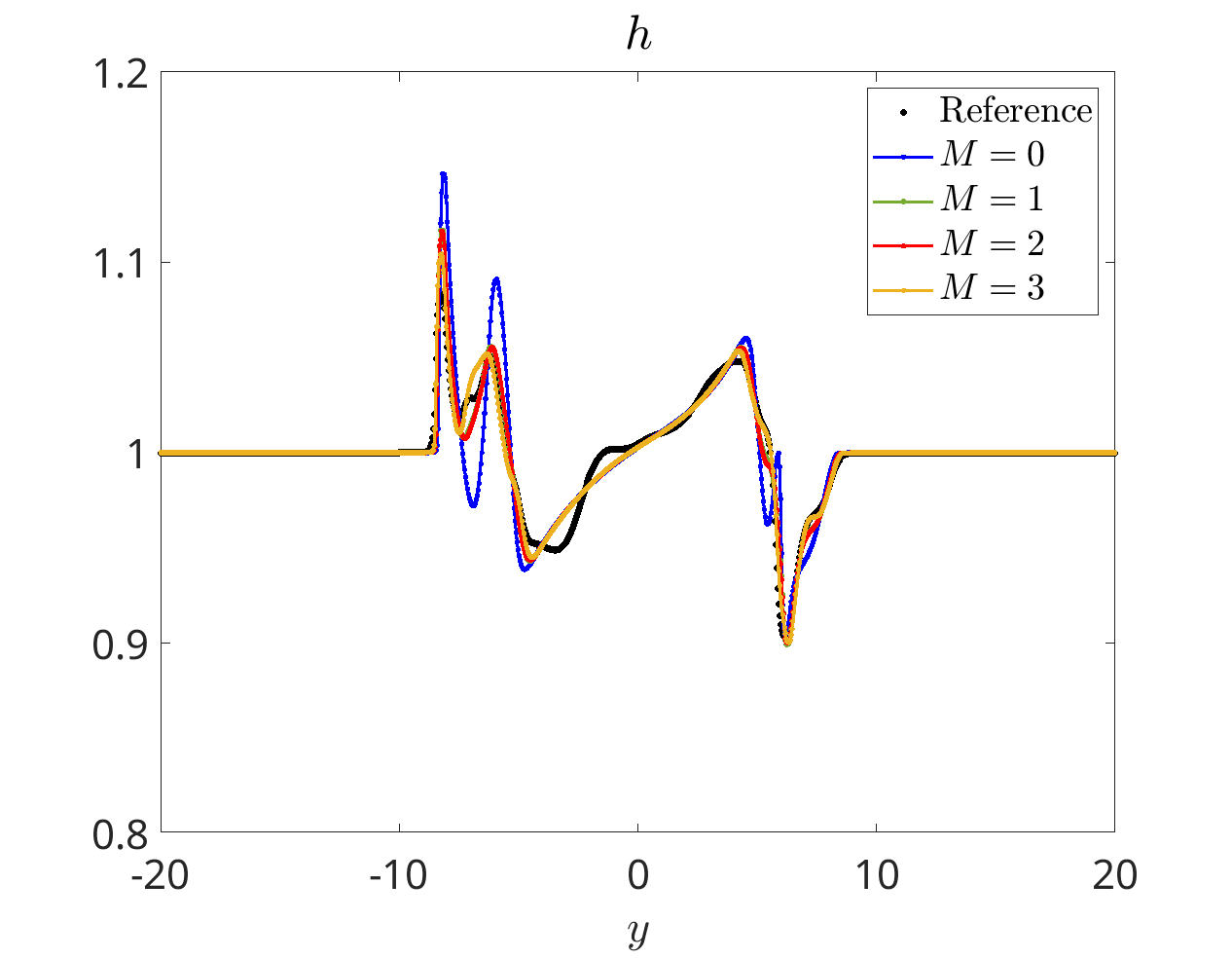}
\includegraphics[trim=0.9cm 0.1cm 1.6cm 0.2cm, clip, width=4.5cm]{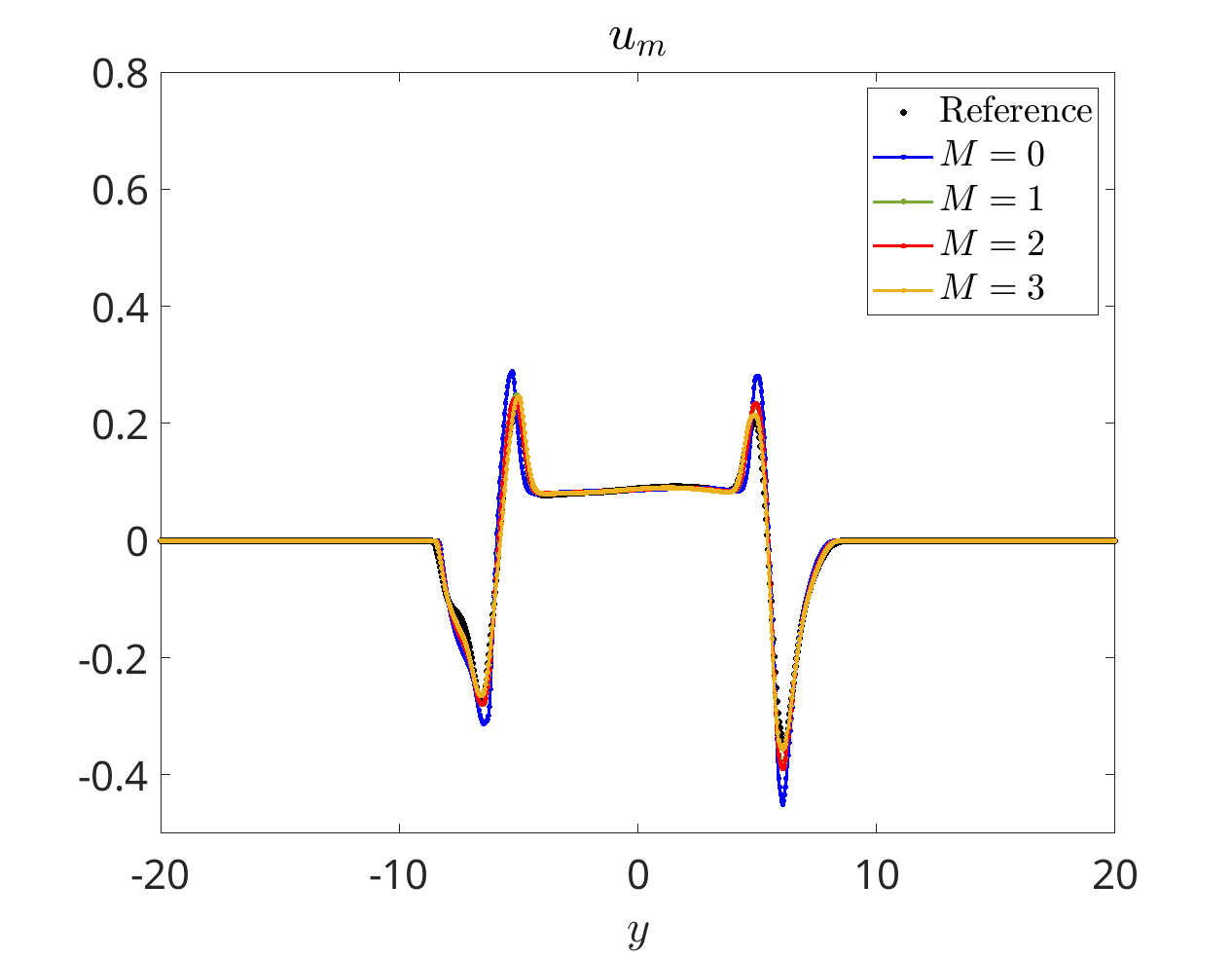}
\includegraphics[trim=0.9cm 0.1cm 1.6cm 0.2cm, clip, width=4.5cm]{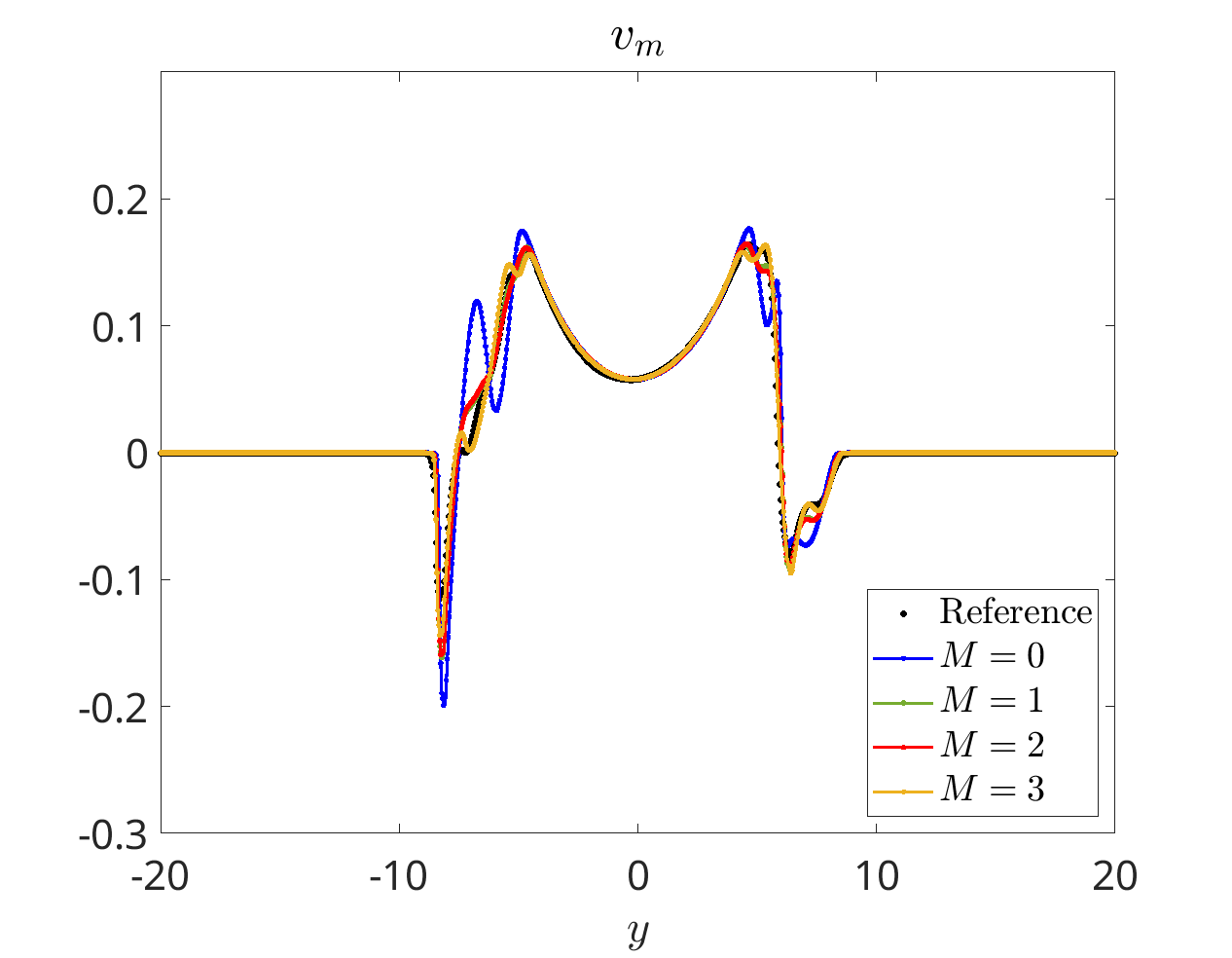}
\includegraphics[trim=0.9cm 0.1cm 1.6cm 0.2cm, clip, width=4.5cm]{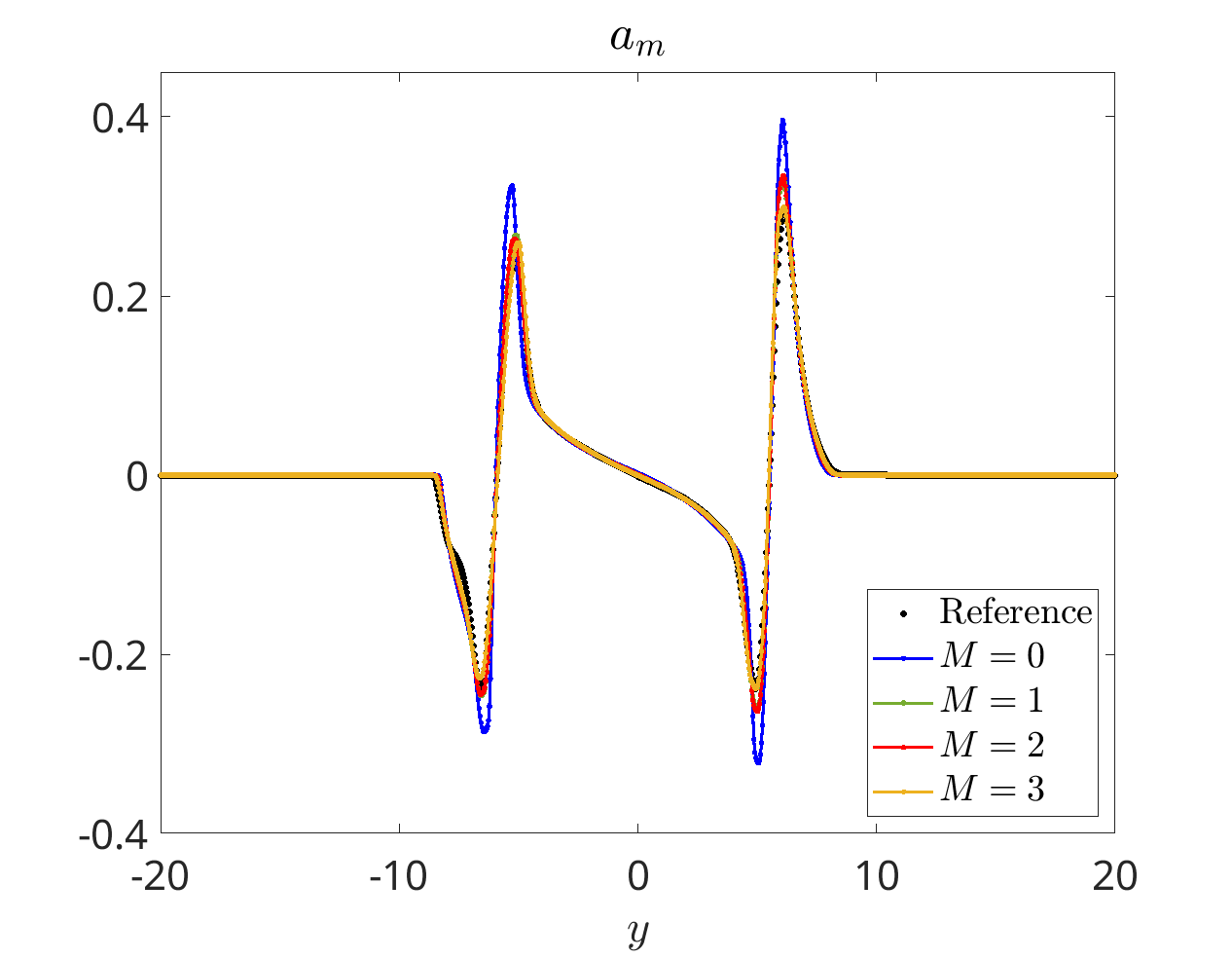}}
\vspace{7pt}
\centerline{
\includegraphics[trim=0.9cm 0.1cm 1.6cm 0.2cm, clip, width=4.5cm]{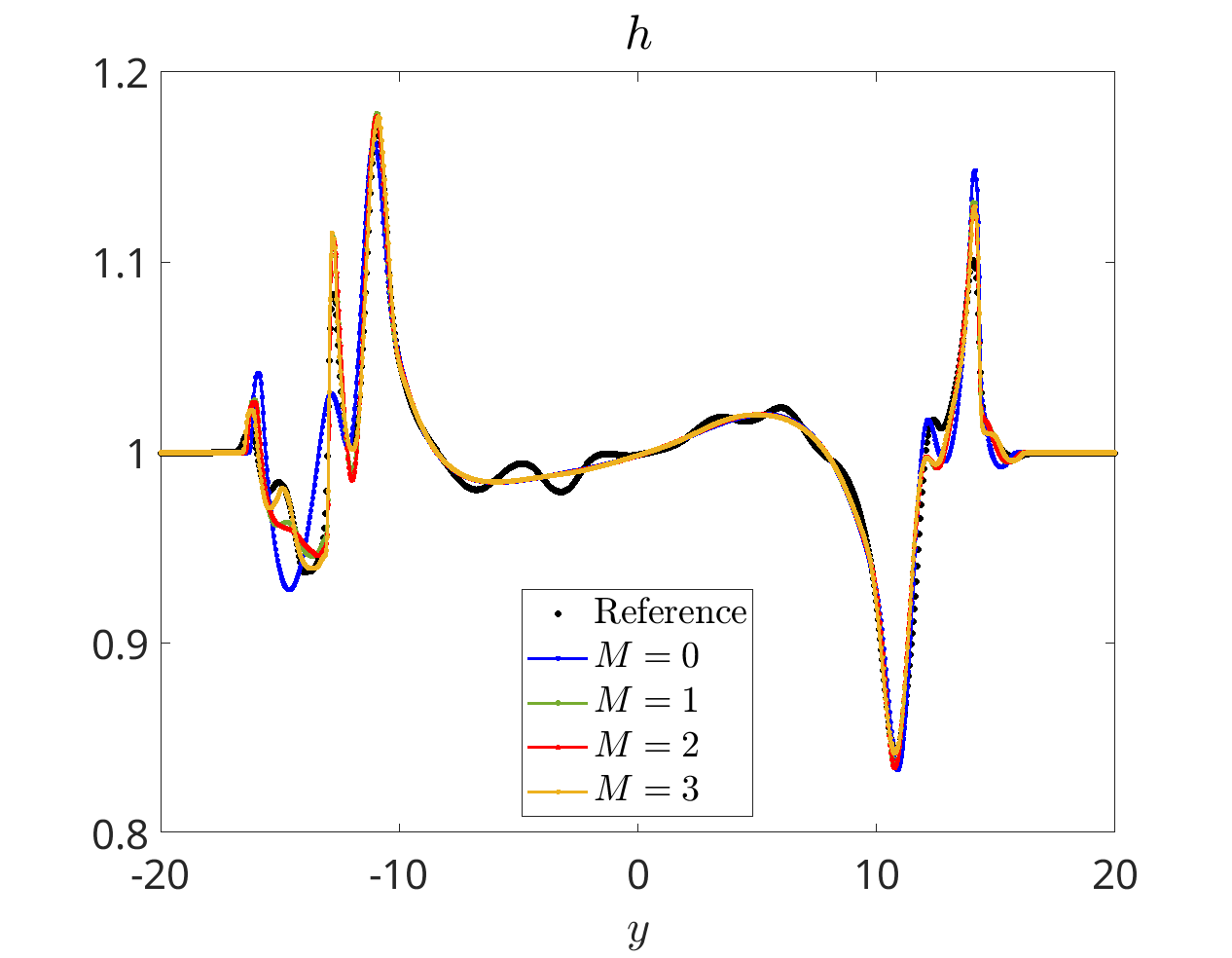}
\includegraphics[trim=0.9cm 0.1cm 1.6cm 0.2cm, clip, width=4.5cm]{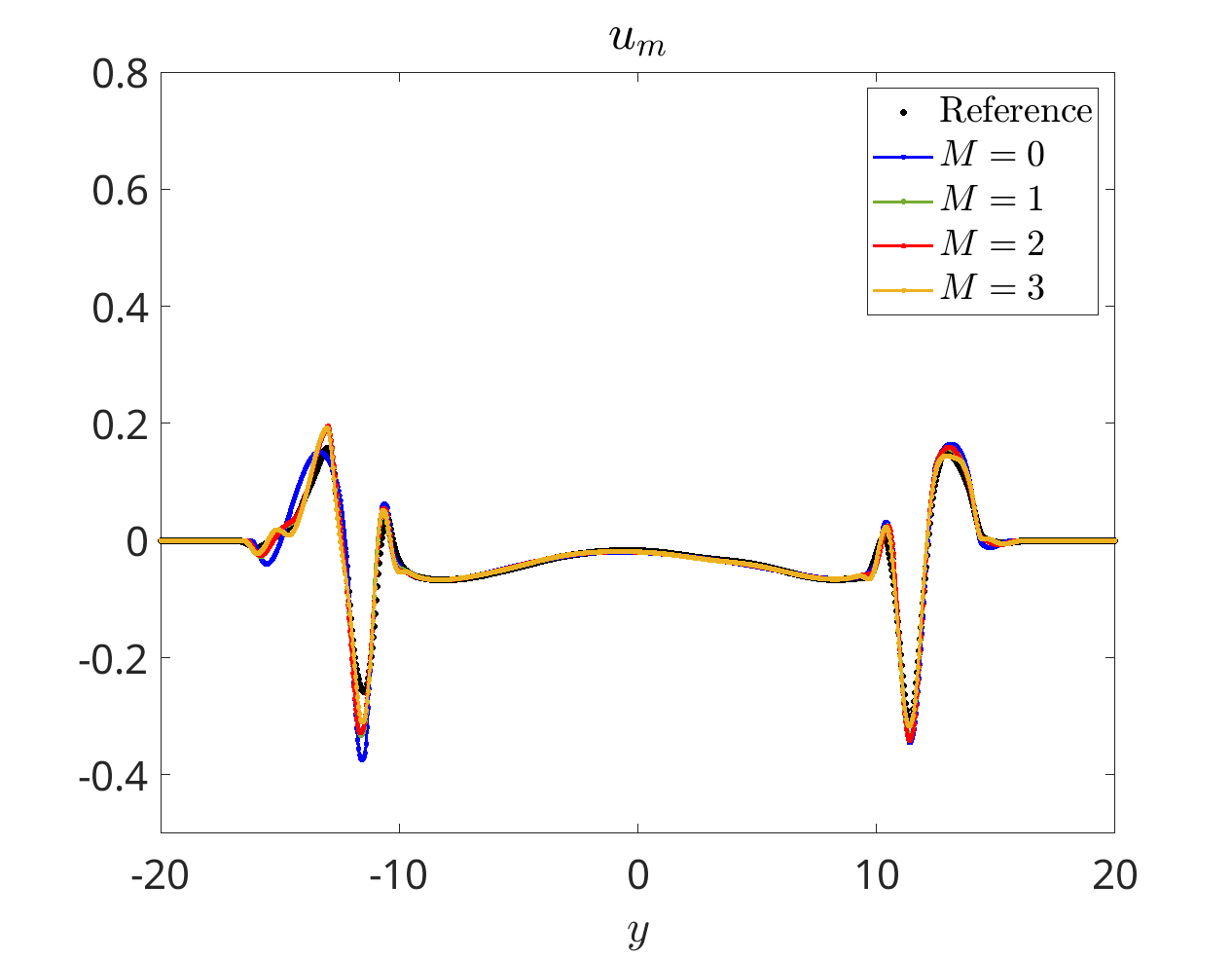}
\includegraphics[trim=0.9cm 0.1cm 1.6cm 0.2cm, clip, width=4.5cm]{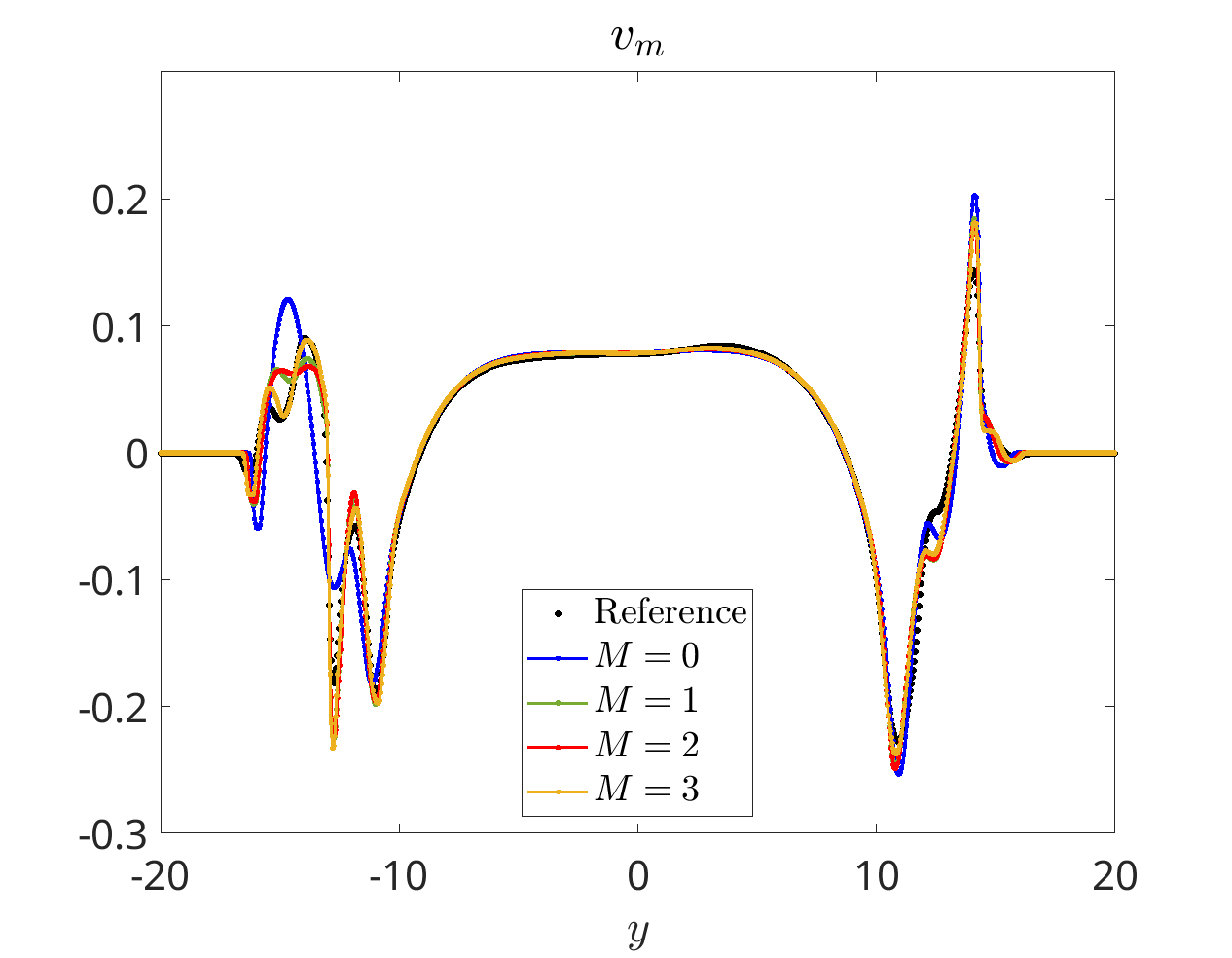}
\includegraphics[trim=0.9cm 0.1cm 1.6cm 0.2cm, clip, width=4.5cm]{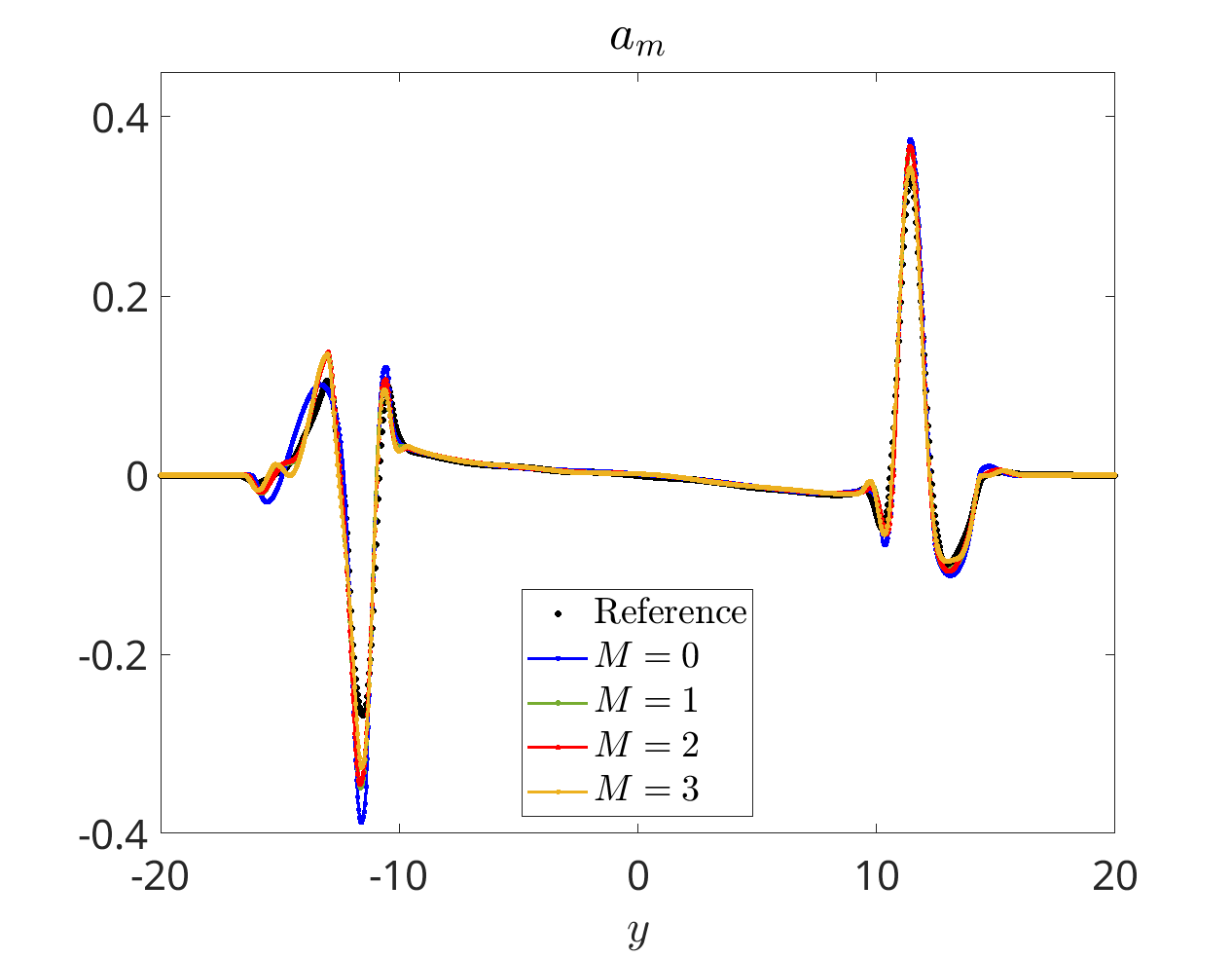}}
\caption{\sf Example 4: Solutions of $h$ (first column), $u_m$ (second column), $v_m$ (third column), and $a_m$ (last column) at times $t = 5$ (top row) and $t = 10$ (bottom row) of the reference system, MRSW, and MRSWME of order $M = 1, \dots,3$.}
\label{fig:Mid_R_comp} 
\end{figure}

The mean value solutions $h,\ u_m,\ v_m,$ and $a_m$ at times $t = 5$ and $t = 10$ for the reference system, MRSW ($M=0$), and MRSWME of order $M = 1, \dots, 3$ are presented in Figure \ref{fig:Mid_R_comp}. 
Much like that of the low Rossby number magneto-geostrophic adjustment, we see that the standard MRSW again captures the material wave and the fast magnetogravity waves in accuracy comparable to that of MRSWME, implying the vertical sinusoidal perturbation in $v$ has little influence on these waves.
Furthermore, the standard MRSW well approximates the expected discontinuities in $h$ and $v_m$ and the strong gradients in $u_m$ and $a_m$ of the reference solution, as does the MRSWME.
Where the sinusoidal vertical perturbation in $v$ has the most influence is again in the Alfv\'en waves, in which we see a significant approximation improvement when using the MRSWME.
Most clearly seen in the left-traveling wave in $h$ and $v_m$, we again see a decreased speed of the Alfv\'en waves. 
Unlike that of the low Rossby number example, though, we do not see significant changes near the center of the domain in $v_m$.

\begin{figure}[ht!]
\centerline{
\includegraphics[trim=0.7cm 0.1cm 1.2cm 0.2cm, clip, width=6cm]{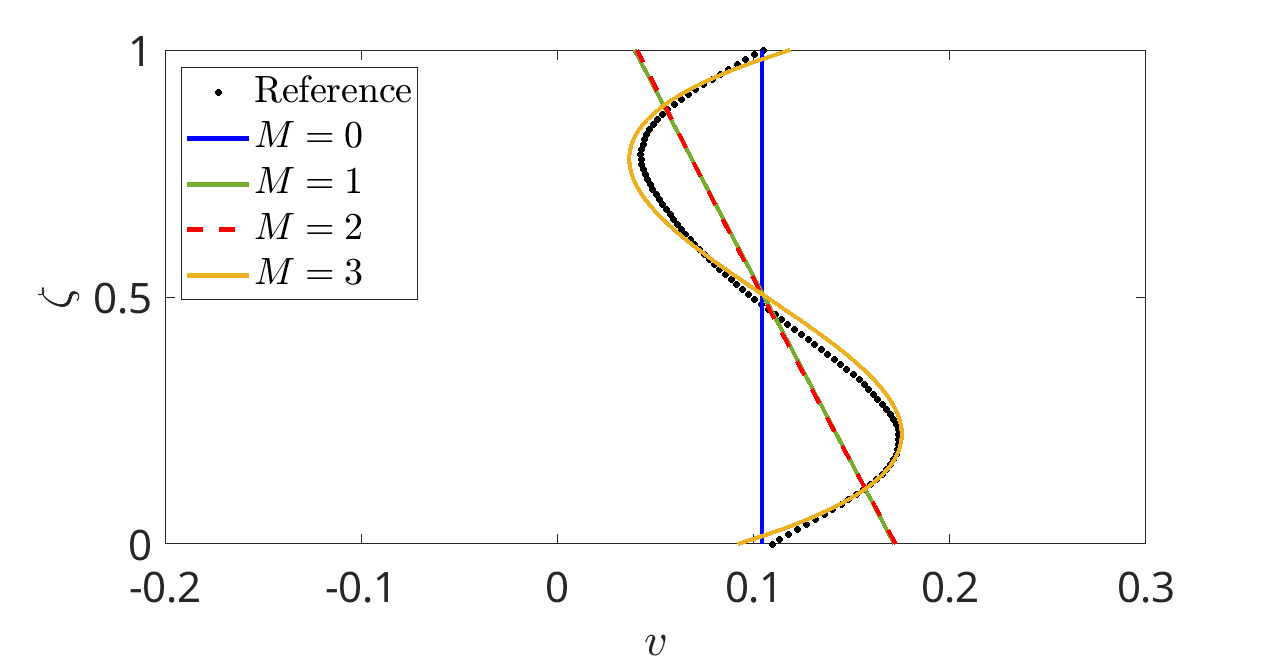}
\hspace{2pt}
\includegraphics[trim=0.7cm 0.1cm 1.2cm 0.2cm, clip, width=6cm]{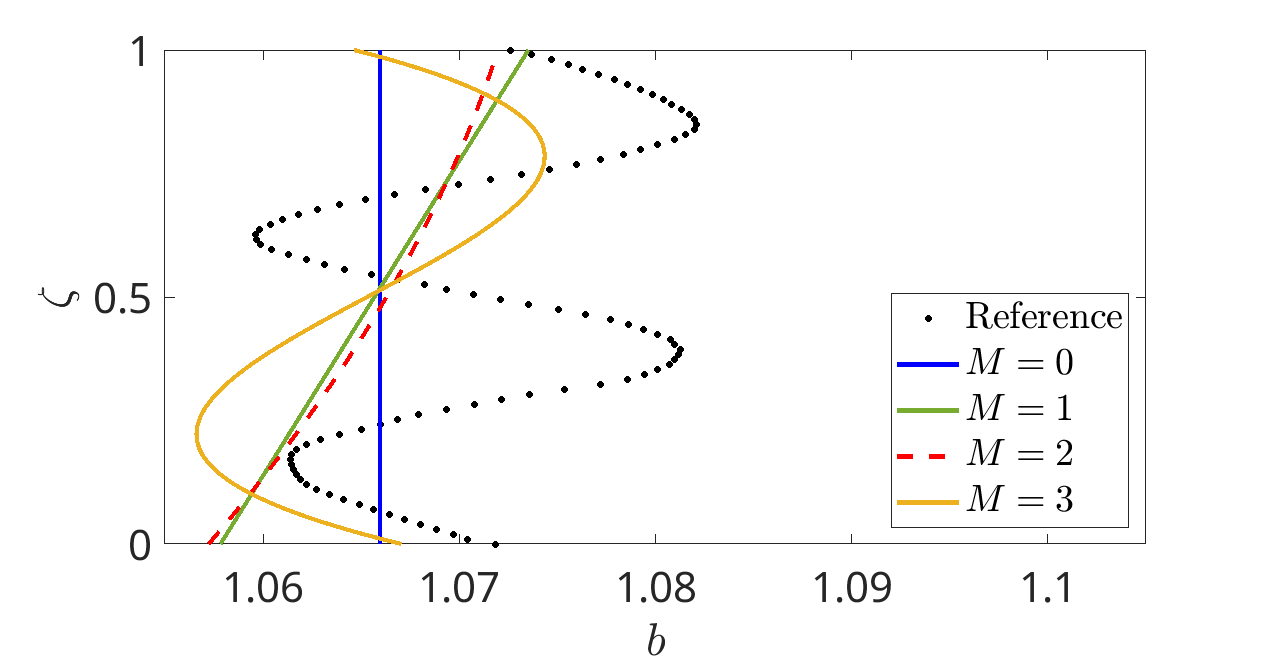}}
\vspace{7pt}
\centerline{
\includegraphics[trim=0.7cm 0.1cm 1.2cm 0.2cm, clip, width=6cm]{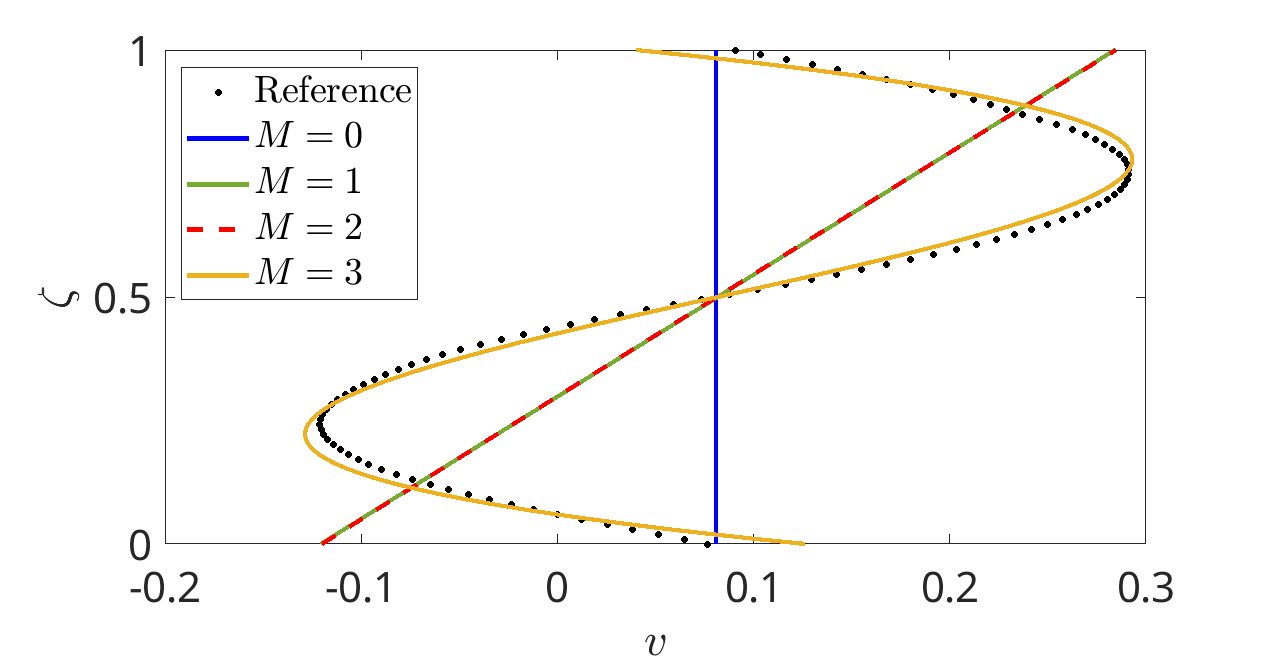}
\hspace{2pt}
\includegraphics[trim=0.7cm 0.1cm 1.2cm 0.2cm, clip, width=6cm]{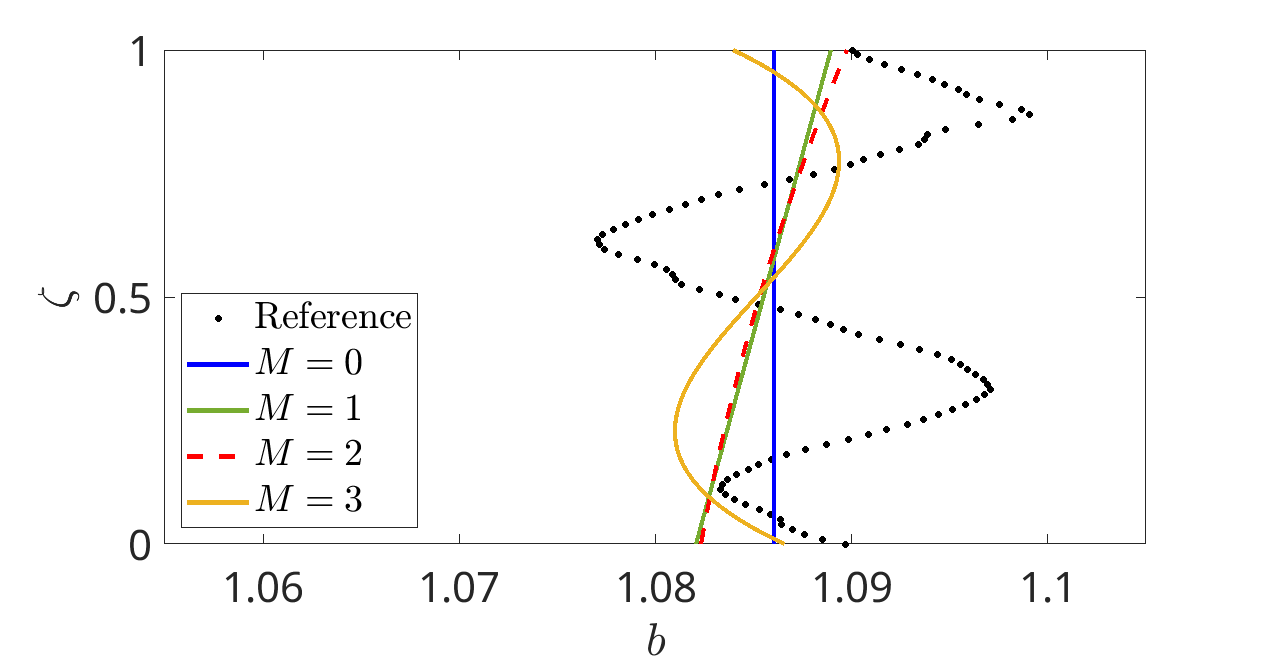}}
\caption{\sf Example 4: Vertical profiles at $y = -5$ of $v$ (left column) and $b$ (right column) at times $t = 5$ (top row) and $t = 10$ (bottom row) for the reference system, MRSW, and MRSWME of order $M = 1,\dots,3$.}
\label{fig:Mid_R_prof} 
\end{figure}

To further demonstrate the increase in the predictive power when using MRSWME, we again present vertical profiles of $v$ and $b$ at slice $y = -5$ in Figure \ref{fig:Mid_R_prof}.
much like that of the smooth magneto-geostrophic adjustment example, we see the sinusoidal structure in the $v$ profile remains as time advances. 
The MRSWME captures this structure accurately through time $t = 10$. 
We also see the profiles for the initial magnetic field $b$ again become more complex in shape as time advances.
While the MRSWME approximation of the profile is improved over the MRSW, it is clear that higher order moment models would need to be taken to capture these more accurately.

\section{Discussion}\label{sec7}
In this paper, we developed a novel 2-D reduced order model of the 3-D incompressible ideal MHD system in the context of magnetic rotating shallow flows. 
This is done by representing the vertical profiles of the horizontal velocities and magnetic fields by a polynomial expansion, in which the polynomial coefficients are now only spatially dependent on the horizontal components. 
By Galerkin projection onto polynomial test functions, i.e. `taking moments', of the the vertically resolved system, we obtain evolution equations for these polynomial coefficients, thus closing the new system.
This hierarchical moment approximation, referred as the MRSWME, maintains non-constant vertical profiles of the horizontal velocities and magnetic fields, ergo increasing the predictive power of these models in comparison to standard shallow fluid approaches (i.e., the MRSW system), while staying within the computationally efficient 2-D framework corresponding to depth-integration. 

Through this moment approximation, the new system gives rise to the Godunov-Powell source term---a nonconservative source dependent on the divergence-free condition of the magnetic field. 
While theoretically zero, keeping this term proves fruitful in the numerical treatment of the magnetic field. 
We take advantage of this naturally arising term to numerically treat the divergence-free condition of the magnetic field, which, after a direct extension of the method proposed in \cite{Chertock2024Locally}, results in a second-order PCCU scheme that is stable in the magnetic field treatment and not dependent on Riemann problem solvers.

Using this method, we present four different numerical examples demonstrating the importance of non-constant vertical profiles of the horizontal magnetic field.
In particular, we present two magneto-geostrophic adjustment problems where we introduce sinusoidal vertical perturbation in the velocity, in which the depth-averaged conditions coincide with the MRSW examples from \cite{CHERTOCK2024WB}.
Using the 2-D vertically-resolved MRSW reference system, we see that this perturbation has significant influence on the outward-propagating Alfv\'en waves, in addition to inducing a vertical profile on the magnetic field.
Furthermore, we show that the standard MRSW system does not properly capture these differences, while the MRSWME does so in an accurate manner.

While the MRSW moment approximation increases model accuracy, it unfortunately comes at the cost of losing global hyperbolicity of the system.
A well known challenge in moment models, this breakdown of hyperbolicity currently limits the capabilities of the presented numerical method. 
Numerical approaches in cases where hyperbolicity is lost will be considered and developed in future work together with hyperbolic regularizations. Further research could consider the numerical preservation of steady-states or an entropy-preserving numerical scheme based on the derivation of a suitable energy equation.

\section*{Author contributions}
\textbf{Julian Koellermeier:} Funding Acquisition, Project administration, Conceptualization, Methodology, Writing---original draft, Writing---review and editing.
\textbf{Michael Redle:} Project administration, Conceptualization, Methodology, Investigation, Formal Analysis, Software, Validation, Visualization, Writing---original draft, Writing---review and editing.
\textbf{Manuel Torrilhon:} Funding Acquisition, Conceptualization, Writing---review and editing.

\section*{Acknowledgments}
This publication is part of the project \textit{HiWAVE} with file number VI.Vidi.233.066 of the \textit{ENW Vidi} research programme, funded by the \textit{Dutch Research Council (NWO)} under the grant \url{https://doi.org/10.61686/CBVAB59929}. 
It is also funded in part by \textit{German Research Foundation (DFG)} Research Unit FOR5409 \textit{"Structure-Preserving Numerical Methods for Bulk- and Interface-Coupling of Heterogeneous Models (SNuBIC)"} (grant \#463312734).

\bibliographystyle{siam}
\bibliography{biblio}
\nocite{*}

\appendix

\section{Scaling and Dimensionless System}\label{secA1}
To non-dimensionalize the ideal incompressible MHD system in \eqref{eq:inc_mhd}, we use the following scaling for the spatial dimensions:
\begin{equation}\label{eq:sc_x}
\widehat{x} = \frac{x}{\ell_0},\quad \widehat{y} = \frac{y}{\ell_0},\quad \widehat{z} = \frac{z}{H}. 
\end{equation}
Since we are interested in free-surface shallow flow, the horizontal length scale $\ell_0$ is much larger than the vertical length scale $H$, and thus $H/\ell_0 = \eps \ll 1$. 
Furthermore, due to shallowness, the vertical velocity and magnetic field will be substantially smaller than their horizontal counterparts, and thus suggest the following scaling:
\begin{equation*}\label{eq:sc_ub}
\widehat{u} = \frac{u}{u_0},\quad \widehat{v} = \frac{v}{u_0},\quad \widehat{w} = \frac{w}{\eps u_0}, \quad
\widehat{a} = \frac{a}{u_0},\quad \widehat{b} = \frac{b}{u_0},\quad \widehat{c} = \frac{c}{\eps u_0},
\end{equation*}
where $u_0$ denotes the characteristic horizontal velocity.
Note again that $\bm b = (a,b,c)^\top$ of \eqref{eq:inc_mhd} is in units of velocity, and hence is scaled by characteristic velocities as well. 
Finally, we introduce a time scale given by the ratio of horizontal spatial and velocity scales, and assume the thermal pressure scales with the hydrostatic pressure with characteristic height $H$; that is, 
\begin{equation}\label{eq:sc_tp}
\widehat{t} = \frac{u_0}{\ell_0}t,\quad \widehat{p} = \frac{p}{\rho g H}. 
\end{equation}
Applying the scalings in \eqref{eq:sc_x}--\eqref{eq:sc_tp} to the ideal incompressible MHD system in \eqref{eq:inc_mhd}, the dimensionless component-wise system reads
\begin{equation*}\label{eq:scaled}
\begin{aligned}
    &\widehat{u}_{\widehat x} + \widehat{v}_{\widehat y} + \widehat{w}_{\widehat z} = 0,\\[4pt]
    & \widehat u_{\widehat t} 
    + \br{\widehat u^2 + \p{\frac{1}{\rm Fr^{2}}\widehat p + \frac{1}{2} \abs{\widehat{\bm b}}^2} - \widehat a^2}_{\widehat x} 
    + \br{\widehat u \widehat v - \widehat a \widehat b}_{\widehat y} 
    + \br{\widehat u \widehat w - \widehat a \widehat c}_{\widehat z} = \frac{1}{\rm Ro}\widehat{v}, \\[4pt]
    & \widehat v_{\widehat t} 
    + \br{\widehat u \widehat v - \widehat a\widehat b}_{\widehat x} 
    + \br{\widehat v^2 + \p{\frac{1}{\rm Fr^{2}}\widehat p + \frac{1}{2} \abs{\widehat{\bm b}}^2} - \widehat b^2}_{\widehat y} 
    + \br{\widehat v\widehat w - \widehat b\widehat c}_{\widehat z} 
    = -\frac{1}{\rm Ro}\widehat{u}, \\[4pt]
    &\eps \p{\widehat w_{\widehat t} 
    + \br{\widehat u \widehat w - \widehat a\widehat c}_{\widehat x}
    + \br{\widehat v \widehat w - \widehat b\widehat c}_{\widehat y}
    + \br{\widehat w^{2} - \widehat c^{2}}_{\widehat z}}
    = -\frac{1}{\eps} \p{\frac{1}{\rm Fr^{2}}\widehat p + \frac{1}{2} \abs{\widehat{\bm b}}^2}_{\widehat z}
    -\frac{1}{\eps}\frac{1}{\rm Fr^{2}}, \\[4pt]
    & \widehat a_{\widehat t} 
    + \p{\widehat a\widehat v-\widehat b\widehat u}_{\widehat y} 
    + \p{\widehat a\widehat w-\widehat c \widehat u}_{\widehat z} = 0,\\[4pt]
    & \widehat b_{\widehat t} 
    + \p{\widehat b\widehat u-\widehat a\widehat v}_{\widehat x} 
    + \p{\widehat b \widehat w-\widehat c\widehat v}_{\widehat z} = 0,\\[4pt]
    & \eps \p{\widehat c_{\widehat t} 
    + \p{\widehat c\widehat u-\widehat a\widehat w}_{\widehat x} 
    + \p{\widehat c \widehat v-\widehat b\widehat w}_{\widehat y}} = 0,\\[4pt]
    & \widehat a_{\widehat x} + \widehat b_{\widehat y} + \widehat c_{\widehat z} = 0,
\end{aligned}
\end{equation*}
where Fr$=u_0/gH$ and Ro$=u_0/Lf$ are the Froude and Rossby numbers, respectively. 
Again assuming that $\eps \ll 1$, the system can be reduced in that the equation for $\widehat{c}$ can be removed, and the equation for the $\widehat{w}$ can be reduced to obtain an explicit formula for the pressure profile. 
In doing so, we obtain the starting system for our derivation in Section \ref{sec3}:

\begin{equation}\label{eq:start}
\begin{aligned}
    &\widehat{u}_{\widehat x} + \widehat{v}_{\widehat y} + \widehat{w}_{\widehat z} = 0,\\[4pt]
    & \widehat u_{\widehat t} 
    + \br{\widehat u^2 + \p{\frac{1}{\rm Fr^{2}}\widehat p + \frac{1}{2} \abs{\widehat{\bm b}}^2} - \widehat a^2}_{\widehat x} 
    + \br{\widehat u \widehat v - \widehat a \widehat b}_{\widehat y} 
    + \br{\widehat u \widehat w - \widehat a \widehat c}_{\widehat z} = \frac{1}{\rm Ro}\widehat{v}, \\[4pt]
    & \widehat v_{\widehat t} 
    + \br{\widehat u \widehat v - \widehat a\widehat b}_{\widehat x} 
    + \br{\widehat v^2 + \p{\frac{1}{\rm Fr^{2}}\widehat p + \frac{1}{2} \abs{\widehat{\bm b}}^2} - \widehat b^2}_{\widehat y} 
    + \br{\widehat v\widehat w - \widehat b\widehat c}_{\widehat z} 
    = -\frac{1}{\rm Ro}\widehat{u}, \\[4pt]
    & \widehat a_{\widehat t} 
    + \p{\widehat a\widehat v-\widehat b\widehat u}_{\widehat y} 
    + \p{\widehat a\widehat w-\widehat c \widehat u}_{\widehat z} = 0,\\[4pt]
    & \widehat b_{\widehat t} 
    + \p{\widehat b\widehat u-\widehat a\widehat v}_{\widehat x} 
    + \p{\widehat b \widehat w-\widehat c\widehat v}_{\widehat z} = 0,\\[4pt]
    & \widehat a_{\widehat x} + \widehat b_{\widehat y} + \widehat c_{\widehat z} = 0,  
\end{aligned}
\end{equation}
where
\begin{equation}\label{eq:mhp_sc}
    \frac{1}{\rm Fr^{2}}\widehat p + \frac{1}{2} \abs{\widehat{\bm b}}^2 
    = \frac{1}{\rm Fr^{2}}\p{\widehat h(\widehat t, \widehat x, \widehat y) 
    + \widehat Z(\widehat t, \widehat x, \widehat y) - \widehat z },
\end{equation}
is referred to as the `magneto-hydrostatic pressure' balance, and is obtained via the assumption that the total pressure at the free surface $\widehat h + \widehat Z$ is equal to zero.

\section{Moments of vertical coupling terms}\label{secA2}
To complete the derivation of the moment evolution equations for the $x$-momentum, we must compute the following integrals:
$$ 
	\int_0^1 \phi_i \p{hu\omega}_\zeta \ d\zeta, \quad {\rm and}\quad  
	\int_0^1 \phi_i \p{ha\mc C}_\zeta \ d\zeta.
$$
Note that similar integrals must be computed in the derivation for the other moment equations, but we leave those out for brevity. 
Starting with the integration involving the velocity vertical coupling, we substitute in the expansion in \eqref{eq:expan} for $u$ and the vertical coupling definition in \eqref{eq:omega_new}, then integrate by parts:
\begin{equation}\label{eq:int_om_stp1}
\begin{aligned}
	\int_0^1 \phi_i \p{hu\omega}_\zeta \ d\zeta 
	= &- \int_0^1 \phi_i \br{\p{u_m + \sum_{n = 1}^M \alpha_n \phi_n}
	\p{\p{h \int_0^\zeta u - u_m \ d\hat\zeta}_x + 
	\p{h \int_0^\zeta v - v_m \ d\hat\zeta}_y}  }_\zeta\ d\zeta, \\[4pt]
	= & - \phi_i \p{u_m + \sum_{n = 1}^M \alpha_n \phi_n}
	\br{\p{h \int_0^\zeta u - u_m \ d\hat\zeta}_x + 
	\p{h \int_0^\zeta v - v_m \ d\hat\zeta}_y} 
	\Bigg|_{\zeta = 0}^{\zeta = 1}\\[4pt]
	& + \int_0^1 \phi_i' \p{u_m + \sum_{n = 1}^M \alpha_n \phi_n}
	\br{\p{h \int_0^\zeta u - u_m \ d\hat\zeta}_x + 
	\p{h \int_0^\zeta v - v_m \ d\hat\zeta}_y}\ d \zeta.
\end{aligned}
\end{equation}

The first remaining term in \eqref{eq:int_om_stp1} vanishes by the definition of $u_m$ and $v_m$.
We then expand the remaining terms to obtain
\begin{align*}
	\int_0^1 \phi_i \p{hu\omega}_\zeta \ d\zeta & 
	= u_m\br{\p{h \alpha_n}_x \sum_{n = 1}^M\p{ \int_0^1 \phi_i'
	\p{\int_0^\zeta \phi_n \d \hat\zeta}\d \zeta }
	+ \p{h \beta_n}_y \sum_{n = 1}^M \p{\int_0^1 \phi_i'
	\p{\int_0^\zeta \phi_n \d \hat\zeta}\d \zeta} }\\[4pt]
	& + \sum_{\ell = 1}^M\sum_{n = 1}^M \alpha_n\br{\p{h \alpha_\ell}_x 
	\p{ \int_0^1 \phi_i'
	\p{\int_0^\zeta \phi_\ell\d \hat\zeta}\phi_n}\d \zeta
	+\p{h \beta_\ell}_y \p{ \int_0^1 \phi_i'  
	\p{\int_0^\zeta \phi_\ell \d \hat\zeta} \phi_n\d \zeta}
	}.
\end{align*}
Using again integration by parts and Liebniz rule, one can show 
$$
	\int_0^1 \phi_i'\p{\int_0^\zeta \phi_n \d \hat\zeta}\d \zeta = 
	\begin{cases}
		\displaystyle-\frac{1}{2i+1} & {\rm if }\ i = n,\\[4pt]
		0 & {\rm otherwise}.
	\end{cases}
$$
Lastly, by defining,
\begin{equation}\label{eq:B_iln}
	B_{i\ell n} = (2i+1) \int_0^1 \phi_i'  
	\p{\int_0^\zeta \phi_\ell \ d \hat\zeta}\phi_n \ d \zeta,
\end{equation}
one recovers
\begin{equation}\label{eq:int_omega}
	\int_0^1 \phi_i \p{hu\omega}_\zeta \ d\zeta
	= -\frac{1}{2i+1}\p{u_m \br{\p{h\alpha_i}_x+\p{h\beta_i}_y}
	-\sum_{\ell=1}^M \sum_{n=1}^M B_{i\ell n} \alpha_n \br{\p{h\alpha_\ell}_x+\p{h\beta_\ell}_y}}.
\end{equation}
%

To compute the integral involving the magnetic field coupling term, we similarly apply the expansion of $a$ defined in \eqref{eq:expan} and the integral definition of $h\mc C$ in \eqref{eq:mcC_new}, then integrate by parts:
\begin{align*}
	\int_0^1 \phi_i &\p{ha\mc C}_\zeta \ d\zeta = 
	- \int_0^1 \phi_i \br{\p{a_m + \sum_{n=1}^M \gamma_n \phi_n}
	\p{\p{h \int_0^\zeta a \ d\hat\zeta}_x + 
	\p{h \int_0^\zeta b \ d\hat\zeta}_y}  }_\zeta\ d\zeta, \\[4pt]
	=& -\phi_i(1) \p{a_m + \sum_{n=1}^M \gamma_n \phi_n(1)}
	\br{\p{ha_m}_x + \p{hb_m}_y} \\[4pt]
	& + \int_0^1 \phi_i' \p{a_m + \sum_{n=1}^M \gamma_n \phi_n}
	\br{\p{h \int_0^\zeta a\ d\hat\zeta}_x + 
	\p{h \int_0^\zeta b \ d\hat\zeta}_y}\ d \zeta,
	\\[4pt]
 	=& -\phi_i(1) \p{a_m + \sum_{n=1}^M \gamma_n \phi_n(1)}
	\br{\p{ha_m}_x + \p{hb_m}_y} \\[4pt]
	& + \int_0^1 \phi_i' \p{a_m + \sum_{n=1}^M \gamma_n \phi_n}
	\br{\p{h \int_0^\zeta a_m \ d\hat\zeta}_x + 
	\p{h \int_0^\zeta b_m\ d\hat\zeta}_y}\ d \zeta,\\[4pt]
	& + \int_0^1 \phi_i' \p{a_m + \sum_{n=1}^M \gamma_n \phi_n}
	\br{\p{h \int_0^\zeta a-a_m \ d\hat\zeta}_x + 
	\p{h \int_0^\zeta b-b_m \ d\hat\zeta}_y}\ d \zeta.
\end{align*}
Notice that the last term remaining is identical to that of the velocity vertical coupling in \eqref{eq:int_om_stp1}, except velocities are replaced by magnetic fields. Hence, we will denote this term as $\boxed{\star}$ until the remaining terms are evaluated. 
Therefore:
\begin{align*}
	\int_0^1 \phi_i \p{ha\mc C}_\zeta \ d\zeta =&\ \boxed{ \star} -\phi_i(1) 
	\p{a_m + \sum_{n=1}^M \gamma_n \phi_n(1)}
	\br{\p{ha_m}_x + \p{hb_m}_y}\\[4pt]
	&+ \int_0^1 \p{\phi_i}_\zeta \p{a_m + \sum_{n=1}^M \gamma_n \phi_n}
	\br{\p{h \int_0^\zeta a_m \ d\hat\zeta}_x + 
	\p{h \int_0^\zeta b_m\ d\hat\zeta}_y}\ d \zeta,\\[4pt]
	=&\ \boxed{ \star} -\phi_i(1) \p{a_m + \sum_{n=1}^M \gamma_n \phi_n(1)}
	\br{\p{ha_m}_x + \p{hb_m}_y}\\[4pt]
	&+ \br{\p{ha_m}_x + \p{hb_m}_y}\br{a_m\int_0^1 \zeta \phi_i'
	\ d \zeta  + \sum_{n=1}^M \gamma_n \int_0^1 \zeta 
	\phi_i' \phi_n\ d \zeta}.
\end{align*}
Notice that all remaining terms outside of $\boxed{\star}$ contain the divergence-free condition of the magnetic field in \eqref{eq:newdivb}, and thus are all theoretically zero. 
However, as stated before, keeping these terms has proven beneficial in discretizations when having to manage the divergence-free condition, and therefore we will simplify these terms as well. 
One can simplify the remaining integrals using integration by parts and orthogonality of the Legendre polynomials:
\begin{align*}
	\int_0^1 \zeta \phi_i' \ d\zeta 
	& = \phi_i(1),\\[4pt]
	\int_0^1 \zeta \phi_i \phi_n\ d\zeta 
	& = \phi_i(1)\phi_n(1) - \frac{1}{2i+1} 
	- \int_0^1 \zeta\phi_i \phi_n'\ d \zeta.
\end{align*}
In turn, this gives
\begin{align*}
	\int_0^1 \phi_i \p{ha\mc C}_\zeta \ d\zeta =&\ \boxed{ \star} -\phi_i(1) 
	\p{a_m + \sum_{n=1}^M \gamma_n \phi_n(1)}
	\br{\p{ha_m}_x + \p{hb_m}_y}\\[4pt]
	&+\br{\p{ha_m}_x + \p{hb_m}_y}\br{a_m \phi_i(1)+ \sum_{n=1}^M \gamma_n
	\p{\phi_i(1)\phi_n(1) - \frac{1}{2i+1} 
	- \int_0^1 \zeta\phi_i \phi_n'\ d \zeta}}\\[4pt]
	=& \ \boxed{ \star} - \br{\p{ha_m}_x + \p{hb_m}_y}\br{\frac{1}{2i+1}
	\sum_{n=1}^M\gamma_n\p{1 + \Gamma_{in}}},
\end{align*}
where $\Gamma_{in} = (2i+1) \int_0^1\zeta\phi_i \phi_n'd \zeta$.
Substituting in $\boxed{\star}$, we obtain the integral term involving the magnetic field vertical coupling:
\begin{equation}\label{eq:int_mcC}
\begin{aligned}
	\int_0^1 \phi_i \p{ha\mc C}_\zeta \ d\zeta =& 
	-\frac{1}{2i+1}\Bigg(a_m \br{\p{h\gamma_i}_x+\p{h\eta_i}_y}
	-\sum_{\ell=1}^M \sum_{n=1}^M B_{i\ell n} \gamma_n \br{\p{h\gamma_\ell}_x
	+\p{h\eta_\ell}_y}\\
	&\hspace{6cm} + \br{\p{ha_m}_x + \p{hb_m}_y}
	\sum_{n=1}^M\gamma_\ell\p{1 + \Gamma_{i\ell}}\Bigg).
\end{aligned}
\end{equation}
As discussed in \S\ref{sec4.3}, the last term is the influence of the Godunov-Powell source term on the evolution equations for the vertical profile coefficients.

\end{document}